\documentclass[11pt]{amsart}
\usepackage{mathrsfs, amsmath, amssymb, array,graphics, booktabs,diagbox,listings,xcolor,textcomp,picinpar,wrapfig,cutwin,longtable,slashbox,amsthm,tikz,threeparttable,makecell,cases,longtable,pdflscape}

\usepackage[backref]{hyperref}
\usepackage{multicol}
\usepackage{multirow}
\usepackage[page,title,titletoc,header]{appendix}
\usepackage[font=small,labelfont=bf]{caption}
\DeclareCaptionFont{tiny}{\tiny}
\DeclareCaptionFont{normalsize}{\normalsize}

\usepackage{geometry}
\definecolor{listinggray}{gray}{0.9}

\definecolor{lbcolor}{rgb}{0.9,0.9,0.9}

\usepackage{algorithm}
\usepackage[noend]{algpseudocode}

\def\baselinestretch{1.2}

\newtheorem{theorem}{Theorem}[section]

\newtheorem{proposition}[theorem]{Proposition}

\theoremstyle{definition}
\newtheorem{definition}[theorem]{Definition}
\newtheorem{example}[theorem]{Example}

\theoremstyle{remark}
\newtheorem{remark}[theorem]{Remark}

\numberwithin{equation}{section}

\newcommand{\HH}{\mathbb{H}}

\newcommand{\ZZ}{\mathbb{Z}}

\begin{document}

\title{Finite-volume Hyperbolic Coxeter $4$-dimensional polytopes with $7$ facets}

\author{Jiming Ma}
\address{School of Mathematical Sciences \\Fudan University\\Shanghai 200433, China} \email{majiming@fudan.edu.cn}

\author{Fangting Zheng}
\address{Department of Mathematical Sciences \\ Xi'an Jiaotong-Liverpool University\\ Suzhou 200433,	China  }
\email{Fangting.Zheng@xjtlu.edu.cn}

\date{Dec. 23, 2024}
\thanks{Jiming Ma was partially supported by  NSFC  12171092. Fangting Zheng was supported by NSFC 12101504.}

\begin{abstract}
	In this paper, we obtain a complete classification of $331$ finite-volume hyperbolic Coxeter $4$-dimensional polytopes with $7$ facets.
\end{abstract}

\maketitle

\section{Introduction}

A Coxeter polytope is a convex polytope where each pair of intersecting bounding hyperplanes meets at a dihedral angle of $\frac{\pi}{k}$ for some integer $k\geq 2$. Let $\mathbb{X}^d$, where $d\geq 2$, represent the unit sphere $\mathbb{S}^d$, the Euclidean space $\mathbb{E}^d$, or the hyperbolic space $\mathbb{H}^d$. If $\Gamma\subset \text{Isom}(\mathbb{X}^d)$ is a finitely generated discrete reflection group, its fundamental domain is a Coxeter polytope $P$ in $\mathbb{X}^d$. Conversely, if $\Gamma=\Gamma(P)$ is generated by reflections in the bounding hyperplanes of a Coxeter polytope $P\subset\mathbb{X}^d$, then $\Gamma$ is a discrete group of isometries of $\mathbb{X}^d$ and $P$ serves as its fundamental domain. In other words, by iteratively reflecting a $d$-dimensional Coxeter polytope along its facets, we obtain a tessellation of $\mathbb{X}^d$. The group $\Gamma$ acts freely and transitively on this tessellation. Henceforth, we do not distinguish between Coxeter polytopes and finitely generated discrete reflection groups. A presentation for $\Gamma$ is $$\langle s_1, s_2, \ldots, s_m ~|~  s_i^2 = (s_is_j)^{k_{ij}}=1\rangle,$$ where $m$ is the number of facets of the polytope $P$. The index $i$ ranges from $1$ to $m$, and the pair $(i,j)$ varies among the incident facets, with a dihedral angle of $\frac{\pi}{k_{ij}}$.

There is a vast body of literature in this field. It was shown that 
any spherical Coxeter polytope (respectively, Euclidean) is a simplex (respectively, a simplex or a direct product of simplices). 
 For comprehensive lists of spherical and Euclidean Coxeter polytopes, readers can refer to, such as, \cite{Coxeter:1934} by Coxeter and \cite{Bourbaki: 1968} by  Bourbaki.

However, the classification of hyperbolic Coxeter polytopes remains far from being achieved. It has been proven by Vinberg \cite{Vinberg:1985} that no compact hyperbolic Coxeter polytope exists in dimensions $d\geq 30$, and by Prokhorov that finite volume hyperbolic Coxeter polytopes do not exist in dimensions $d \geq 996$ \cite{Prokhorov:1987}. These bounds, however, may not be sharp. Examples of compact hyperbolic Coxeter polytopes are known only up to dimension $8$ \cite{Bugaenko:1984,Bugaenko:1992}, while non-compact but finite volume hyperbolic Coxeter polytopes are known just up to dimension $21$ \cite{Vinberg:1972,VK:1978,Borcherds: 1998}. Under these dimension bounds, Allcock constructed infinitely many compact and non-compact but finite volume Coxeter polytopes in $\mathbb{H}^d$ for each $d\leq 6$ and $d\leq 19$, respectively \cite{Allcock: 2006}. Additionally, Bogachev-Douba-Raimbault exhibited that the garlands made by compact Markarov polytopes give infinitely many commensurability classes of compact Coxeter polytopes in $\mathbb{H}^4$ and $\mathbb{H}^5$ \cite{BDR: 2023}.  In terms of classification, complete results are only available in dimensions less than or equal to three. The classification of $2$-dimensional hyperbolic polytopes completed by Poincar\'{e} \cite{Poincare: 1882} has been utilized in the work of Klein and Poincare on discrete groups of isometries of the hyperbolic plane. Much later in 1970, Andreev proved a characterization theorem for the $3$-dimensional hyperbolic acute-angled polytopes \cite{Andreev1: 1970,Andreev2: 1970}. His result played a fundamental role in Thurston's work on the  geometrization of 3-dimensional Haken manifolds.

In higher dimensions, although a complete classification is not available, many interesting examples are presented in \cite{Makarov: 1965,Makarov: 1966,Vinberg: 1967,Makarov: 1968,Vinberg: 1969,Rusmanov: 1989,ImH: 1990,Allcock: 2006,Vinberg:2015,Alex:2023}. Additionally, enumerations are reported for cases where the differences between the numbers of facets $m$ and the dimensions $d$ of the corresponding polytopes are small numbers.  It is noteworthy that the quotient space of $\mathbb{H}^d$ by hyperbolic reflection groups with simpler presentations often yields hyperbolic orbifolds of small volume.  When $m-d=1$, Lann\'{e}r classified all compact hyperbolic Coxeter simplices \cite{Lanner: 1950}. The classification of all non-compact but finite volume hyperbolic simplices has been exhibited by several authors e.g., \cite{Bourbaki: 1968,Vinberg: 1967,Koszul:1967}. For $m-d=2$, Kaplinskaya described all finite volume hyperbolic Coxeter simplicial prisms \cite{Kaplinskaya: 1974}. Esselmann \cite{Esselmann: 1996} later enumerated other compact possibilities in this family, which are named \emph{Esselmann polytopes}. 
Tumarkin \cite{Tumarkin: n2} classified all other non-compact but finite volume hyperbolic Coxeter $d$-dimensional polytopes with
$d + 2$ facets, and the hyperbolic orbifold with minimal volume belongs to this family. In the case of $m-d=3$, Esselman has proven that compact hyperbolic Coxeter $d$-polytopes with $d+3$ facets  only exist when $d\leq 8$ \cite{Esselmann}. By expanding the techniques derived by Esselmann in \cite{Esselmann} and \cite{Esselmann: 1996}, Tumarkin completed the classification of compact hyperbolic Coxeter $d$-polytopes with $d+3$ facets \cite{Tumarkin: n3}. In the non-compact but finite volume case, Tumarkin has proven that such polytopes do not exist in dimensions greater than or equal to $17$ \cite{Tumarkin: n3fv,Tumarkin: n3fvs}; there is a unique such polytope in dimension  $16$.  Moreover, the author provided in the same papers the complete classification of a family of  pyramids over a product of three simplices, which exists only in dimensions of $4,5,\cdots,9$ and $13$. The classification for the non-compact but finite volume case has not been completed yet. Regarding this subfamily, Roberts provided a non-pyramidal list with exactly one non-simple vertex \cite{Roberts:15}. In this paper, we classify all the hyperbolic Coxeter $4$-polytope with $7$ facets. 
The main theorem is as follows:

\begin{theorem}\label{thm:main}
	There are exactly $331$ finite-volume hyperbolic Coxeter $4$-polytopes with $7$ facets, out of which $40$ are compact and $291$ are non-compact but of finite volume.
\end{theorem}

We would like to note that the compact list we present here coincides with the classification provided by Tumarkin in \cite{Tumarkin: n3} and the $4$-pyramids with $7$ facets are the same with those in \cite{Tumarkin: n3fv}. Additionally, we include the $140$ polytopes that were missing from Roberts's census in dimension $4$, thereby achieving a comprehensive coverage of $164$ non-pyramidal Coxeter $4$-polytopes with $7$ facets and one non-simple vertex. These polytopes are realized over the five combinatorial types  $P_{7}$, $P_{8}$, $P_{9}$, $P_{10}$, and $P_{11}$ in Table \ref{table: combinatoricSelction} and the combinatorial types $P_{9}$ and $P_{10}$ are involved in Roberts's original list. For further details,  please refer to Section \ref{section:validation}.  Additionally, for an up-to-date overview of the current knowledge on hyperbolic Coxeter polytopes, readers can consult Anna Felikson's webpage \cite{Annahomepage}.

The paper \cite{JT:2018} is the main inspiration for our recent work on enumerating hyperbolic Coxeter polytopes. In contrast to \cite{JT:2018}, we utilize a more versatile ``block-pasting" algorithm, initially introduced in \cite{MZ:2018}, rather than the attempted ``tracing back" method. We incorporate additional geometric constraints and optimize our program to significantly reduce computational load. Our algorithm efficiently enumerates hyperbolic Coxeter polytopes across various combinatorial types, not limited to just the $n$-cube. It is worth noting that the algorithm and procedure employed here are more intricate compared to those described in our previous works \cite{MZ:2022} and \cite{MZ:2023}, where the combinatorial structures were simpler due to compactness.

Last but not least, our primary motivation for studying hyperbolic Coxeter polytopes lies in their potential for constructing high-dimensional hyperbolic manifolds and orbifolds, some of which exhibit small or even minimal volumes. This scheme is beyond the scope of this paper. Interested readers can refer to works such as \cite{KM: 2013, Kellerhals: 2014} for related topics and further references.

The paper is organized as follows. We provide in Section $2$ some preliminaries about hyperbolic Coxeter polytopes. In Section $3$, we recall the 2-phases procedure and related terminologies introduced by Jacquemet and Tschantz \cite{Jacquemet2017,JT:2018} for numerating all hyperbolic Coxeter $n$-cubes. The $31$ combinatorial types of $4$-polytopes with $7$ facets are reported in Section $4$. We choose from them the $16$  polytopes with admissible vertex links and prepare censuses of dihedral angles around vertices accordingly for the next step of calculations. Enumeration of all the \emph{SELCper-potential matrices} is explained in Section $5$. We move on for the Gram matrices of actual hyperbolic Coxeter polytopes in Section $6$. Validations and complete lists of  result Coxeter polytopes  are shown in Section $7$.

\textbf{Acknowledgment}
We extend our sincere gratitude to both referees for their insightful comments and constructive suggestions. The first referee provided accurate and detailed feedback on many key points requiring revision, which significantly enhanced the precision of our statements.  For example, clarification was suggested regarding the optimal choice of $7$ and explanation about admissible vertex links. The meticulous review by the second referee, including a time-consuming near re-execution of our computational procedures to verify the results, is highly valued.  Authors are deeply moved and inspired by this dedication. This referee can even be considered the third author of this result, as a conscientious and independent review is crucial to the final presentation and credibility of a classification work involving computer-assisted traversal calculations. The oversight of the $16$ polytopes over the combinatorial type $P_4$ in our initial manuscript, as well as some mistakes in the volume and cusp number data presented in Table 12, due to errors in data conveyance, were pointed out by this referee. Additionally, the referee's advise of presenting the results in a tabular format instead of a graphical one has been implemented. This change not only enhances the clarity and machine readability of our findings but also reduces the risk of typographical errors. The corresponding Coxeter diagrams, as suggested by the first referee, are also included in the appendix of this arXiv version.

\section{preliminary} \label{section:cchp}
In this section, we recall some essential facts about Coxeter hyperbolic polytopes, including Gram matrices, Coxeter diagrams, characterization theorems, etc. Readers can refer to, for example, \cite{Vinberg:1993} for more details.  

\subsection{Hyperbolic space, hyperplane and convex polytope}
We first describe the Lobachevsky model of the $d$-dimensional hyperbolic space $\mathbb{H}^d$, which is the unique, up to isometry, complete simply-connected Riemannian manifold with constant sectional curvature $-1$. Let $\mathbb{E}^{d,1}$ be a ($d+1$)-dimensional Euclidean vector space equipped with a Lorentzian scalar product $\langle\cdot,\cdot \rangle$ of signature $(d,1)$. A vector in this space is categorized as time-like, light-like, or space-like if its square Lorentzian norm is negative, zero, or positive, respectively. In this quadratic space $\mathbb{E}^{d,1}$, we can define a model, also denoted by $\mathbb{H}^d$ for easing the notation, of the $d$-dimensional hyperbolic space to be
 $$\mathbb{H}^d=\{x=(x_1, ...,x_d,x_{d+1})\in \mathbb{E}^{d,1}:\langle x,x\rangle=-1, x_{d+1}>0\},$$ 
 with the distance function $d_{\mathbb{H}^d}(x,y)=\text{arccosh}(-\langle x,y\rangle)$. In the following, we always employ $\mathbb{H}^d$ to refer this model, unless otherwise indicated. The positive definiteness of the scalar product on the tangent space of every point in $\mathbb{H}^d$ can be confirmed by observing that the Lorentzian-orthogonal vector of the tangent space at every point in $\mathbb{H}^d$ is space-like.

Let $O(n,1)$ be the group of linear isomorphisms of $\mathbb{E}^{d,1}$ that preserve the Lorentzian scalar product and $O^+(n,1)$ is the index two subgroup that preserves the upper sheet $\mathbb{H}^d$. It can be shown that $\text{Isom}(\mathbb{H}^d)=O^+(n,1)$.

For an integer $1\leq k\leq d$, a $k$-dimensional affine geodesic subspace of $\mathbb{H}^d$ is the intersection of a $(k+1)$-dimensional vector space in $\mathbb{E}^{d,1}$ with $\mathbb{H}^d$, provided it is non-empty. This intersection is itself isometric to $\mathbb{H}^k$. In the case of $k=1$, the intersection gives rise to a geodesic. To be more specific, the geodesic originating from $p\in \mathbb{H}^d$ with unit velocity vector $v\in T_p \mathbb{H}^d$ can be parametrized as  $l(t)=\cosh(t)p+\sinh(t)v$.

Consider geodesic half-lines $l_i:[0,\infty)\rightarrow \mathbb{H}^d$ with constant unit speed. The set $\partial \mathbb{H}^d$, termed as the \emph{boundary of $\mathbb{H}^d$}, represents the collection of all geodesic half-lines, subjecting to the following equivalence relation: $l_i\sim l_j$ if and only if $\sup\limits_{t\in [0,\infty)}\{d_{\mathbb{H}^d}(l_i(t),l_j(t))\}<\infty$. Equivalently, the boundary $\partial\mathbb{H}^d$ can be identified as the set of Euclidean norm $1$ light-like vectors with $x_{n+1}>0$. 
The hyperbolic space $\mathbb{H}^d$ can be compactified by incorporating its boundary, resulting in the extended hyperbolic space $\overline{\mathbb{H}^d} = \mathbb{H}^d \cup \partial \mathbb{H}^d$. The points of the boundary $\partial \mathbb{H}^d$ are called the \emph{ideal points}.

The affine subspaces of $\mathbb{H}^d$ of dimension $d-1$ are \emph{hyperplanes}. In particular, every hyperplane of $\mathbb{H}^d$ can be represented as $$H_e=\{x\in \mathbb{H}^d~|~\langle x,e\rangle=0\},$$ where $e$ is a space-like vector of Lorentzian norm 1. The half-spaces bounded by $H_e$ are denoted by $H_e^+$ and $H_e^-$, where
\begin{equation}
	H_e^-=\{x\in \mathbb{H}^d~|~\langle x,e\rangle\leq 0\}. \label{1}
\end{equation} 

The \emph{mutual disposition of hyperplanes} $H_{e}$ and $H_{f}$ can be described in terms of the corresponding two space-like vectors $e$ and $f$ as follows: 
\begin{itemize}
	\item The hyperplanes $H_{e}$ and $H_{f}$ intersect if $\vert\langle e,f \rangle\vert<1$. The value of the dihedral angle of $H_{e}^-\cap H_{f}^-$, denoted by $\angle H_e H_f$, can be obtained via the formula $$\cos \angle H_e H_f=-\langle e,f\rangle;$$
	\item The hyperplanes $H_{e}$ and $H_{f}$ are parallel if  $\vert\langle e,f\rangle\vert=1$;
	\item The hyperplanes $H_{e}$ and $H_{f}$ diverge if $\vert\langle e,f\rangle\vert>1$. The distance $\rho(H_e,H_f)$ between $H_e$ and $H_f$, when $H_e^+\subset H_f^-$ and $H_f^+\subset H_e^-$, is determined by  $$\cosh \rho(H_e,H_f)=-\langle e,f \rangle.$$
\end{itemize}


We say a hyperplane $H_e$ \emph{supports} a closed bounded convex set $S$ if $H_e\cap S \ne \empty 0$ and $S$
lies in one of the two closed half-spaces bounded by $H_e$. If a hyperplane $H_e$ supports $S$, then $H_e \cap S$ is called a \emph{face} of $S$. 

\begin{definition}
	A $d$-dimensional convex non-degenerated hyperbolic polytope is a subset of the form 
	
	\begin{equation}
		P=\mathop{\cap}\limits_{i\in \mathcal{I}} H_{i}^-\in \mathbb{H}^d, \label{2}
	\end{equation}
	where $H_i^-$ is the negative half-space bounded by the hyperplane $H_i$ in $\mathbb{H}^d$,
	under the following assumptions:
	\begin{itemize}
		\item  $P$ contains a non-empty open subset of $\mathbb{H}^d$;
		\item  Every bounded subset $S$ of $P$ intersects only finitely many $H_{i}$.
	\end{itemize}
\end{definition}

A convex polytope of the form (\ref{2}) is termed an \emph{acute-angled} polytope if, for distinct $i$ and $j$, either $\angle H_iH_j\leq \frac{\pi}{2}$\footnote{We follow the convention in [\cite{Vinberg:1993}, Ch. 6] to include right angles for the acute-angled polytopes, which may be atypical in other contexts.} or $H_i^+\cap H_j^+=\emptyset$. It is evident that Coxeter polytopes satisfy this condition. We denote $e_i$ as the unit space-like normal vector to $H_i$, meaning $e_i$ is Lorentzian-orthogonal to $H_i$ and points away from $P$.

In the following, a $d$-dimensional convex polytope $P$ is referred to as a \emph{$d$-polytope}. A $j$-dimensional face is called a $j$-face of $P$. Specifically, a $(d-1)$-face is termed a \emph{facet} of $P$, and a $0$-face is an ordinary vertex of $P$. We assume that each hyperplane $H_i$ intersects $P$ at its facet, meaning the hyperplane $H_i$ is uniquely determined by $P$ and is known as a \emph{bounding hyperplane} of the polytope $P$. Henceforth, we always assume that the intersection in the definition of a polytope only involves bounding hyperplanes.

Let $\overline{P}$ denote the closure of a \emph{finite-volume} polytope $P$ in $\overline{\mathbb{H}^d}$. The term  ``finite-volume" signifies that $P$ is a subset of $\mathbb{H}^d$ with a finite hyperbolic volume. The closure process involves the inclusion of ideal points corresponding to representatives of equivalent classes of geodesic half-lines contained in $P$. The vertices of $\overline{P}$ lying on boundary $\partial \mathbb{H}^d$ and interior $\mathbb{H}^d$ are referred to as \emph{ideal} and \emph{ordinary} vertices of $P$, respectively. A finite-volume hyperbolic polytope $P$ is \emph{compact} if all its vertices are ordinary. Otherwise, it is non-compact but of finite volume and has at least one ideal vertex.

\subsection{Gram matrices, Perron-Frobenius Theorem, and Coxeter diagrams}\label{hcp}

Most of the content in this subsection is well-known by  peers in this field. We present them here for the convenience of the readers. 
In particular, Theorems \ref{thm:signature} and \ref{Vinberg:thm3.1} are extremely important throughout this paper. 

For a hyperbolic acute-angled $d$-polytope $P=\mathop{\cap}\limits_{k\in \mathcal{I}} H_{k}^-$, we consider its Gram matrix $G(P)$ whose ($ij$)-entries are the products $\langle e_i,e_j\rangle$, where the vectors are from the set $\{e_k\in \mathbb{E}^{d,1}\vert i\in\mathcal{I}\}$ that determining hyperplanes $H_k^-$s. Many of the combinatorial, metrical, and arithmetic properties of $P$ can be inferred from $G(P)$. For instance, the entry $\langle e_i,e_j\rangle$ characterizes the configuration of bounding hyperplanes as follows:

\begin{center}
	$\langle e_i,e_j\rangle=
	\left\{
	\begin{array}{ccl}1  \hspace{0.7cm} &\mbox~~~~~{{\rm if}}& j=i, \\
		-\cos\alpha_{ij}&\mbox ~~~~~{{\rm if}}&  H_i~ \text{and}~ H_j~ \text{intersect} \text{ at~a~dihedral~ angle~}\alpha_{ij},\\
		-\cosh \rho_{ij}&\mbox ~~~~~{{\rm if}}&  H_i~ \text{and}~ H_j~ \text{divergent~at~a~distance}~\rho_{ij},\\
		-1  \hspace{0.7cm} &\mbox~~~~~{{\rm if}}& H_i~\text{and}~H_j~\text{are~parallel}.
	\end{array}\right.
	$
\end{center}

For a Coxeter polytope, it is convenient to represent it by a weighted graph, which is known as the \emph{Coxeter graph} and is denoted by $\Sigma=\Sigma(P)$. Every node $i$ in $\Sigma$ corresponds a bounding hyperplane $H_i$ of $P$. Two nodes $i_1$ and $i_2$ are joined by an edge with a weight $2\leq k_{ij} < \infty$ or the symbol $\infty$, if $H_i$ and $H_j$ intersect in $\mathbb{H}^n$ with angle $\frac{\pi}{k_{ij}}$ or diverge. If two hyperplanes $H_i$ and $H_j$ have a common perpendicular of length $\rho_{ij}>0$ in $\mathbb{H}^n$, the nodes $i_1$ and $i_2$ are joined by a dotted edge, sometimes labeled $\cosh \rho_{ij}$. In the following, an edge of weight $2$ is omitted. Edges with weights $3$, $4$, $5$, and $6$ are drawn as simple, double, triple, and quadruple edges, respectively.  Edges with weight $\infty$ are drawn thick. Moreover, the \emph{order}  of diagram $\Sigma$, denoted as $|\Sigma|$, refers to the total number of nodes it contains. The \emph{signature} and \emph{rank} of diagram $\Sigma$ correspond to the signature and rank, respectively, of the matrix $G(\Sigma)$.

A \emph{permutation} of a square matrix refers to a rearrangement of its rows combined with an identical permutation of its columns. A square matrix, denoted as $M$, is said to be a direct sum of matrices $M_1, M_2, \cdots, M_n$ if, through a suitable permutation, it can be transformed into the form
\begin{center}
	$ \begin{pmatrix}
		M_1&&&&0 \\
		&M_2&&&\\
		&& \cdot&&\\
		&&&\cdot& \\
		0&&&&M_n
	\end{pmatrix}_. $
	
\end{center}

\noindent A matrix $M$ that cannot be represented as a direct sum of two matrices is said to be \emph{indecomposable}\footnote{It is also referred to as ``irreducible" in some references.}. Every matrix can be represented uniquely as a direct sum of indecomposable matrices, which are called (indecomposable) components. 
We say a polytope is \emph{indecomposable} if its Gram matrix $G(P)$ is indecomposable. Note that a non-degenerate hyperbolic polytope $P$ is decomposable if it has a proper face $F$ that is orthogonal to every hyperplane which does not contain it. The orthogonal projection onto the face $F$ establishes a fibration of $P$ into polyhedral cones, ensuring that all vertices of $P$ reside within $F$.  Consequently, every convex hyperbolic polytope of finite volume is indecomposable.

\begin{figure}[h]
	\scalebox{0.26}[0.26]{\includegraphics {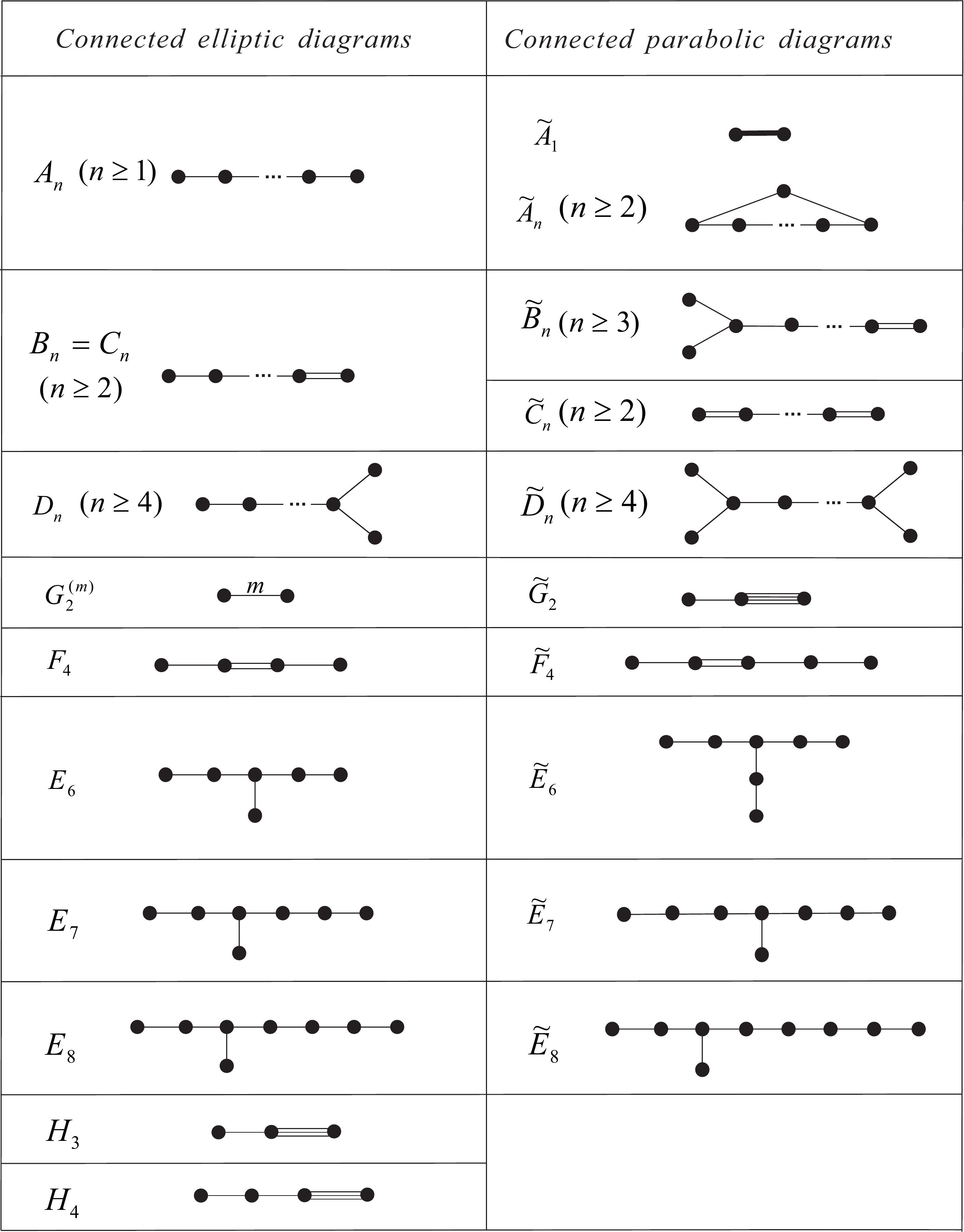}}
	\caption{Connected elliptic (left) and connected parabolic (right) Coxeter diagrams.}
	\label{figure:coxeter}
\end{figure}

In 1907, Perron found a remarkable property of eigenvalues and eigenvectors of matrices with positive entries. Frobenius later generalized it by investigating the spectral properties of indecomposable non-negative matrices.

\begin{theorem}[Perron-Frobenius, \cite{G:1959}]\label{thm:PF}
	An indecomposable real matrix $A=(a_{ij})$ with non-negative entries always possesses a maximal eigenvalue $r$, which is positive, simple, and the corresponding eigenvector has all positive coordinates.
\end{theorem}

The Gram matrices $G(P)$ of an indecomposable Coxeter polytope are symmetric matrices with non-positive entries off the diagonal, and all diagonal elements are equal to $1$. By applying Theorem \ref{thm:PF} to $I-G(P)$, where $I$ represents the identity matrix, we can identify a unique smallest eigenvalue of $G(P)$, denoted by $\lambda$. Depending on the sign of $\lambda$, $G(P)$ can be classified as positive definite, semi-positive definite, or indefinite when $\lambda$ is greater than, equal to, or less than zero, respectively. In the case of being semi-definite, once again by utilizing Theorem \ref{thm:PF}, the deficiency of an indecomposable semi-positive definite matrix $G(P)$ does not exceed $1$, and any proper submatrix of $G(P)$ is positive definite. For a Coxeter polytope $P$, its Coxeter diagram $\Sigma(P)$ is said to be \emph{elliptic} if $G(P)$ is positive definite; $\Sigma (P)$ is called \emph{parabolic} if every indecomposable component of $G(P)$ is degenerate and all proper subdiagrams of each indecomposable component are elliptic. The elliptic and connected parabolic diagrams are proven to be the Coxeter diagrams of spherical and Euclidean Coxeter simplices, respectively. They are classified  by Coxeter \cite{Coxeter:1934} as shown in Figure \ref{figure:coxeter}.

\begin{figure}[h]
	\scalebox{0.35}[0.35]{\includegraphics {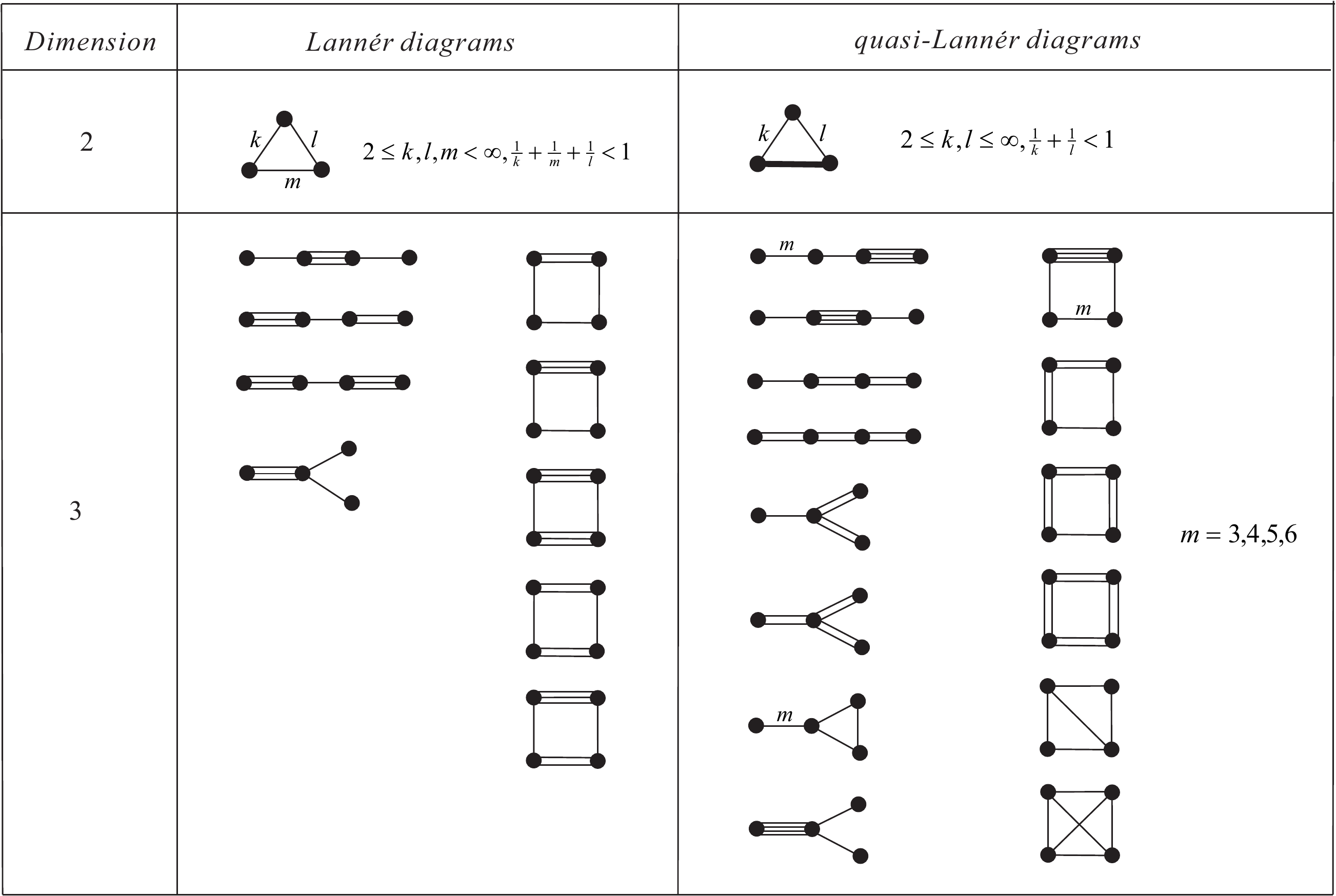}}
	\caption{The Lann\'{e}r and quasi-Lann\'{e}r diagrams correspond to $2$- or $3$-dimensional hyperbolic Coxeter simplices, respectively.}
	\label{figure:lanner}
\end{figure}

A connected diagram $\Sigma$ is referred to as a \emph{Lann\'{e}r diagram} if it is neither elliptic nor parabolic, and every proper subdiagram of $\Sigma$ is elliptic. These diagrams represent compact hyperbolic Coxeter simplices. The complete list of such diagrams was first reported by Lann\'{e}r \cite{Lanner: 1950}, and the lists for $2$- and $3$-dimensional cases are provided in Figure \ref{figure:lanner}. 
Similarly, non-compact hyperbolic simplices of finite volume can also be enumerated. They are known as quasi-Lann\'{e}r diagrams, and characterized by connected non-elliptic and non-parabolic Coxeter graphs where proper subgraphs are either elliptic or connected parabolic graphs.  The $2$- and $3$-dimensional lists of these simplices are presented in Figure \ref{figure:lanner}.

Although the full list of hyperbolic Coxeter polytopes remains incomplete, some powerful algebraic restrictions are known \cite{Vinberg:1985s}:

\begin{theorem} \label{thm:signature}
	(\cite{Vinberg:1985s}, Th. 2.1). Let $G=(g_{ij})$ be an indecomposable symmetric matrix
	of signature $(d,1)$, where $g_{ii}=1$ and $g_{ij}\leq 0$ if $i\ne j$. Then there exists a unique, up to isometry, convex hyperbolic polytope $P\subset\mathbb{H}^d$, whose Gram matrix coincides with $G$.
\end{theorem}

The counterpart theorem in the spherical and Euclidean settings is analogous, with the only difference being that the signature condition is revised to $(d+1,0)$ and $(d,0)$, respectively. This distinction also implies that we cannot expect a comparatively straightforward characterization of the potential combinatorial types of acute-angled polytopes in hyperbolic space.

\begin{theorem} \label{Vinberg:thm3.1} 
	(\cite{Vinberg:1985s}, Th. 3.1, Th. 3.2)
	Let $P=\mathop{\cap}\limits_{i\in I} H_i^- \in \mathbb{H}^d$ be an  acute-angled polytope and $G=G(P)$ be the Gram matrix. Denote $G_J$ the principal submatrix of G formed from the rows and columns whose indices belong to $J\subset I$. Then, 
	\begin{enumerate}
		\item The intersection $\mathop{\cap}\limits_{j\in J}H_j^-\in \mathbb{H}^d,J\subset I$, is a face $F$ of $P$ if and only if the matrix $G_J$ is positive definite;
		\item    The intersection $\overline{\mathop{\cap}\limits_{j\in J}H_j^-},J\subset I$, is an ideal vertex of $P$ if and only if the matrix  $G_J$ is a parabolic matrix of rank $d-1$.
	\end{enumerate}
\end{theorem}

\begin{theorem} \label{Vinberg:thm4.2} 
	(\cite{Vinberg:1985s}, Th. 4.2)
	Let $P=\mathop{\cap}\limits_{i\in I} H_i^- \in \mathbb{H}^d$ be a convex polytope and $G=G(P)$ be its Gram matrix. 
	Denote $G_J$ the principal submatrix of G formed from the rows and columns whose indices belong to $J\subset I$. Let $\mathcal{F}$ (resp. $\widetilde{\mathcal{F}}$) be the collection of elliptic submatrices $G_J$s (resp. elliptic and parabolic submatrices $G_J$s ). Partially order $\mathcal{F}$ (resp. $\widetilde{\mathcal{F}}$) by submatrix relation.
	Then, $P$ is compact (resp. finite-volume) if and only if any of the  following conditions holds:
	
	\begin{enumerate}
		\item The polytope $P$ contains at least one ordinary vertex (resp. ordinary or ideal vertex). For every ordinary vertex (resp. vertex) of $P$ and any edge of $P$ emanating from it, there is precisely one other ordinary vertex (resp. vertex) of $P$ on that edge.
		\item  The partially ordered set   $\mathcal{F}$ (resp. $\widetilde{\mathcal{F}}$) is isomorphic to the poset of some $d$-dimensional abstract combinatorial polytope.
		
	\end{enumerate}
\end{theorem}

Therefore, in order to classify all hyperbolic Coxeter $4$-polytopes with $7$ facets, a two-phase approach is employed. Initially, we systematically enumerate all potential matrices over admissible $4$-polytopes with $7$ facets, ensuring they adhere to elliptic conditions around every ordinary vertex and parabolic conditions around every ideal vertex. To enhance the efficiency of our computational endeavors, we strategically introduce additional geometric constraints during this phase, which serve to refine our search criteria without compromising the comprehensiveness of our classification. Subsequently, we conduct a targeted programmed search for matrices with an indecomposable signature  $(d,1)$ Gram matrix, thereby accomplishing our classification.

Last but not least, the volume of an even dimensional hyperbolic Coxeter polytope can be computed via its elliptic subgraphs. And we present some basic information here. Almost all the terminology and fact hold for general finitely generated groups, but we only state for finitely generated discrete reflection groups for our purpose. We have also identified all the arithmetic non-compact finite-volume Coxeter polytopes within our list, which are distinctly labeled with ``A" in the ``Arithmetic" column of Table \ref{table:resultall}. Those without any label are non-arithmetic. Furthermore, the ones labeled with ``C" are compact, and their arithmeticity remains to be checked.

\begin{definition}(Poincar\'{e} series)
	Let $\Gamma$  be a finitely generated group with a generating set $S$ corresponding to a $d$-dimensional hyperbolic polytope $P$. The Poincar\'{e} series of $\Gamma$ is the formal power series $$f(\Gamma,S)=f_S(x)=\mathop{\Sigma}\limits_{g\in\Gamma}x^{l_s(g)},$$ where $l_g$ is the length function of $\Gamma$ with respect to S, i.e. $$l_g(s)=min\{l\in \mathbb{N}~|~ \exists s_1,\cdots,s_l\in S\bigcup S^{-1}~\text{such that}~ g=s_1s_2\cdots s_l\}.$$
\end{definition}

\begin{proposition}\label{prop:volume} (\cite{KP: 2011})
	Given a $d$-dimensional hyperbolic polytope $P$ where $d$ is even, and $\Gamma=\Gamma(P)$ is the corresponding finitely generated discrete reflection group, then
	
	$$\text{volume}(P)=(-1)^{\frac{d}{2}}\cdot\frac{\pi^{\frac{d}{2}}\cdot2^d\cdot \frac{d}{2}!}{d!}\cdot \chi(\Gamma),$$
	
\noindent where $\chi(\Gamma)$ is the orbifold Euler characteristic of the $\mathbb{H}/\Gamma$.

Moreover, let $\mathcal{E}=\{\Sigma^{'}| \Sigma^{'} \text{is an elliptic subgraphs of }\Sigma(P)\}$, and  $\Gamma^{'}$ the Coxeter subgroup corresponds to $\Sigma^{'}$ , by Steinberg's formula, we have

$$\chi(\Gamma)=\sum\limits_{\Sigma^{'}\in \mathcal{E}}\frac{(-1)^{|\Sigma^{'}|}}{f_{S^{'}}(1)},$$ 

\noindent where $f_{S^{'}}(1)$ is the value of Poincar\'{e} series at identity, i.e. the order of the subgroup $\Gamma^{'}$, and $S^{'}\subset S$ is the generating set of the subgroup of $\Gamma^{'}$.

\end{proposition}

By means of the values in Table \ref{table:ellorder} \cite{Guglielmetti:2017}, it is straightforward to compute the hyperbolic volume of an even dimensional hyperbolic Coxeter polytope.

\begin{table}[H]
	{\footnotesize
			\begin{tabular}{cccccccccc}
				\Xcline{1-10}{1.2pt}
				
				$A_n$&$B_n=C_n$&$D_n$&$G_2^{(m)}$&$F_4$& $E_6$&$E_7$&$E_8$&$H_3$&$H_4$\\
				\hline
				$(n+1)!$&$2^n\cdot n!$&$2^{n-1}\cdot n!$&$2m$&1,152&51,840&2,903,040&696,729,600&120&14,400\\
				\Xcline{1-10}{1.2pt}

			\end{tabular}
		}
		
		\hspace*{0.5cm}
		\caption{Orders of finite Coxeter groups corresponding to elliptic diagrams.}
		\label{table:ellorder}
	\end{table} 

Next, we provide a concise definition of an arithmetic hyperbolic Coxeter polytope. Discrete subgroups of finite covolume are termed \emph{lattices}. Readers can refer to \cite{Morris:2015} for a more general definition of arithmetic subgroup of a connected semi-simple Lie group, where a non-orientable lattice is arithmetic if it orientable double cover is arithmetic.  Let $k$ be a totally real number field with the ring of integers $\mathcal{O}_k$, and let $q$ be a quadratic form of signature $(d,1)$ defined over $k$. We require that for every non-identity Galois embedding $\sigma :k\rightarrow \mathbb{R}$, the conjugate form $q^\sigma$ is of signature $(d+1,0)$.  Denote by  $O_q(k\otimes_{\mathbb{Q}}\mathbb{R})$ the set $\{A\in GL(n+1,k\otimes \mathbb{R})~|~\forall x,y \in k^{n+1}\otimes_\mathbb{Q}\mathbb{R},~q(x,y)=q(Ax,Ay)\}$, and let $O_q^+(k\otimes_{\mathbb{Q}}\mathbb{R})$ be  its index two subgroup that preserve one sheet of the hyperboloid $q(x,x)=-1$. The group $\widetilde{\Gamma}=O_q^+(\mathcal{O}_k)$ of the integral automorphisms of $q$ is a lattice of $O_q^+(k\otimes_{\mathbb{Q}}\mathbb{R})\cong O^+(n,1)=\text{Isom}(\mathbb{H}^d)$. The discreteness of G follows from the discreteness of $\mathcal{O}_k$ in $\mathbb{R}^{[k:\mathbb{Q}]}$ and compactness of $O^+_{q^\sigma}(k\otimes_{\mathbb{Q}}\mathbb{R})$ for $\sigma\ne \sigma_{id}$. By the Borel-Harish-Chandra theorem, $\widetilde{\Gamma}$ is of finite covolume.

Let $H$ be a group. Two subgroups $H_1, H_2 < H$ are commensurable if and only if there exists an element $h \in  H$ such that $H_1\cap h^{-1}H_2h$ has finite index in both $H_1$ and $h^{-1}H_2h$. The subgroups $\Gamma$ of the isometry group of $d$-dimensional hyperbolic space that are commensurable with the group $\widetilde{\Gamma}$, obtained in the manner previously described, are known as \emph{arithmetic hyperbolic lattices of the simplest type}. 
We extend this terminology to the corresponding quotient space $\mathbb{H}^n/\Gamma$. In general, there are three types of arithmetically defined hyperbolic lattices. Vinberg showed that all finitely generated hyperbolic discrete reflection groups are the arithmetic lattice of the simplest type \cite{Vinberg: 1967}. The corresponding fundamental domains are refer to as \emph{arithmetic hyperbolic Coxeter polytope}.

Furthermore, Vinberg has formulated an efficient criterion to determine the arithmeticity. Consider a hyperbolic Coxeter polytope $P$ with its Gram matrix denoted as $G(P) = (g_{ij})_{1\leq i,j\leq N}$. Define 
$K(P)$ as the number field $\mathbb{Q}(\{g_{ij}\}_{1\leq i,j\leq N})$
, and let $k(P)$ be the field generated by all cyclic products of the matrix entries of $G(P)$, that is,  $g_{i_1i_2}g_{i_2i_3}\cdots g_{i_{l-1}i_l}g_{i_li_1}$ for any subset $\{i_1,i_2,\cdots,i_l\}\subset\{1,2,\cdots,N\}$. The criterion for arithmeticity is as follows:

\begin{theorem}\label{Vinberg:thm2} 
	(\cite{Vinberg: 1967}, Th. 2)
Let $P$  be a hyperbolic Coxeter polytope. The associated finitely generated discrete reflection group $\Gamma(P)$ is arithmetic if and only if:
\begin{enumerate}
	\item $K(P)$ is totally real;
	\item for any embedding $\sigma:K(P)\rightarrow \mathbb{R}$, such that $\sigma|_{k(P)}\ne Id$, the matrix $G^{\sigma}(P)$, obtained by applying $\sigma$ to all the coefficients of $G(P)$, is semi-positive definite;
	\item all the cycles of $2\cdot G(P)$ are algebraic integers in $k(P)$.
\end{enumerate}
\end{theorem}

Specifically, if a hyperbolic Coxeter polytope $P$ is non-compact but of finite volume, then $k(P)=\mathbb{Q}$. In this case,  $\Gamma(P)$ is arithmetic if and only if all the cycles of the matrix $2\cdot G$ are rational integers. We have checked the arithmeticity of all non-compact but finite-volume hyperbolic Coxeter polytope listed in our study, as depicted in Table \ref{table:resultall}.

\section{Potential hyperbolic Coxeter matrices}\label{section:potential}

Almost all of the terminologies and theorems in this section are proposed by Jacquemet and Tschantz, with only minor adjustments adopted. We recall them here for reference, and readers can refer to \cite{JT:2018} for more details.

\subsection{Coxeter matrices}\label{subsection:coxeterMatrix}

The \emph{Coxeter matrix} of a hyperbolic Coxeter polytope $P$ is a symmetric matrix $M=(m_{ij})_{1\leq i,j\leq N}$ with entries in $\mathbb{N}\cup\{\infty\}$ such that
\[m_{ij}=\left\{\begin{array}{cl} 
	1,&\text{if } j=i,\\ 
	k_{ij}, &\text{if }H_i\text{ and }H_j\text{ intersect in }\HH^n\text{ with angle }\frac{\pi}{k_{ij}},~k_{ij}\in\ZZ_{\geq2},\\
	0, &\text{if }H_i\text{ and }H_j\text{ intersect in }\partial\HH^n,\\
	\infty, & \text{otherwise}.
\end{array}\right.\]
Note that, compared with Gram matrix, the Coxeter matrix does not involve the specific information of the distances of the disjoint pairs.

\begin{remark}
	In the subsequent discussions, we refer to \textit{the Coxeter matrix $M$ of a graph $\Sigma$} as the Coxeter matrix $M$ of the Coxeter polyhedron $P$ such that $\Sigma=\Sigma(P)$.\\
\end{remark}

\subsection{Partial matrices}\label{subsection:partial}

\begin{definition}
	Let $\Omega=\{n\in \mathbb{Z}\,|n=0~\text{or}~ \,n\geq 2\}\cup\{\infty\}$ and let $\bigstar$ be a symbol representing an undetermined real value. A \textit{partial matrix of size $m\geq 1$} is a symmetric $m\times m$ matrix $M$ whose diagonal entries are $1$, and whose non-diagonal entries belong to $\Omega\cup\{\bigstar\}$.
\end{definition}

\begin{definition}
	Let $M$ be an arbitrary $m\times m$ matrix, and $s=(s_1,s_2,\cdots,s_k)$, $1\leq s_1<s_2<\cdots<s_k\leq m$. Let $M^{s}$ be the $k\times k$ submatrix of $M$ with $(i,j)$-entry $m_{s_i,s_j}$. 
\end{definition}

\begin{definition}
	We say that a partial matrix $M=(p_{ij})_{1\leq i,j,\leq m}$ is a \emph{potential matrix} for a given polytope $P$ if 
	
	$\bullet$ There are no entries with the value $\bigstar$;
	
	$\bullet$ There are entries $\infty$ in positions of $M$ that correspond to disjoint pair;
	
	$\bullet$ For every tuple $s$ of indices of facets meeting at an ordinary vertex $v$ of $P$, the matrix, obtained by replacing value $k_{ij}$ with $-\cos\frac{\pi}{k_{ij}}$ of submatrix $M^s$, is elliptic;  For every tuple $s$ of indices of facets meeting at an ideal vertex $v$ of $P$, the matrix, obtained by replacing value $k_{ij}$ with $-\cos\frac{\pi}{k_{ij}}$ and value $0$ with $-1$ of submatrix $M^s$, is parabolic.
	
\end{definition}

For brevity, we use a \emph{potential vector} $$C=(p_{12},p_{13},\cdots p_{1m},p_{23},p_{24},\cdots,p_{2m},\cdots p_{ij},\cdots p_{m-1,m})$$ to denote the potential matrix, where $p_{ij}$'s are from the corresponding potential matrix $M=(p_{ij})_{m\times m}$. The potential matrix  and potential vector $C$  can be constructed one from each other easily. In general, an arbitrary Coxeter matrix corresponds to a Coxeter vector following the same manner. For brevity, we mainly use the language of \emph{vectors} to explain the methodology and  report the results. A Coxeter vector is \emph{connected} if the corresponding Coxeter diagram is connected. It is worth remarking that for a given Coxeter diagram, the corresponding (potential / Gram) matrix and vector are not unique in the sense that they are determined under a given labeling system of the facets and may vary when the system changes. In Section \ref{chapter:algorithm}, we apply a permutation group to the nodes of the diagram and remove the duplicates to obtain all of the distinct desired vectors

For each rank $r\geq 2$, there are infinitely many finite Coxeter groups, because of the infinite 1-parameter family of all dihedral groups, whose graphs consist of two nodes joined by an edge of weight $k\geq 2$. However, a simple but useful truncation can be utilized:

\begin{proposition}
	There are finitely many finite Coxeter groups of rank $r$ with Coxeter matrix entries at most seven.
\end{proposition}

We remark that in order to minimize the number of angle variables involved in the second phase of the signature obstruction computation in \emph{Mathematica}, we have deliberately selected the integer $7$ as the bound. Because $7$ is the smallest weight not appearing in connected elliptic or parabolic diagrams of order at least $3$. 
Thus, it suffices to enumerate potential matrices with entries at most seven, and the other candidates can be obtained from substituting integers greater than seven with the value seven. In the second phase, we will have more angle variables. 
In the following, we always use the terms ``Coxeter matrix" or ``potential matrix" to refer to the one with integer entries less than or equal to seven, unless otherwise mentioned.

\section{ combinatorial type of 4-polytopes with 7 facets}\label{section:4d7f}
The vertex link $lk(v)$ of a vertex $v$ in a polytope $P$ is the intersection of the
polytope $P$ with a sufficiently small sphere centered at $v$. This sphere should not intersect any facet of $P$ that is not incident to $v$. In particular geometric context, the terminology might be adjusted accordingly. For instance, we will use horosphere rather than sphere to form the vertex link if the vertex of a hyperbolic polytope is ideal. 

The classification of combinatorial types of $d$-polytopes with $m$ facets, where $m\leq d+3$, has been accomplished using Gale diagrams. If the vertex link of a vertex $v$ in an $n$-polytope is an ($n-1$)-simplex, then $v$ is referred to as a \emph{simple} vertex. Otherwise, the vertex is \emph{non-simple}. An $n$-polytope is termed \emph{simple} if every vertex is simple. For simple $d$-polytopes with $m$ facets, there is a counting formula by Perles, see for example \cite{G:1967}, and also \cite{BD:1998}. We learn the data of $4$-polytopes with $7$ facets as shown in Table \ref{table: combinatoricSelction} from \cite{FMM:13}. In this text, the numbers  $0$, $1$, $\cdots$, $6$ are used to label the seven facets, with each line signifying a combinatorial type. Each bracket, termed a \emph{facet bracket}, corresponds to a vertex that is intersected by facets with indicated labels. For instance, in the polytope $P_1$, there are $9$ vertices. The first vertex is formed by the intersection of the facets $F_1$, $F_2$, $F_3$, $F_4$, $F_5$, and $F_6$

\begin{table}[h]
\renewcommand\arraystretch{1.5}
{\tiny
	\begin{tabular}{c|c|l}
		\Xcline{1-3}{1.2pt}
label&admissible&facet brackets\\
\hline				
{\color{red}1}&$P_1$&[{\color{green}6},5,{\color{pink}4},{\color{pink}3},2,{\color{green}1}] [6,5,4,0] [6,5,3,0] [6,4,2,0] [6,3,2,0] [5,4,1,0] [5,3,1,0] [4,2,1,0] [3,2,1,0]\\
\hline
2&&[6,5,4,3,2,1] [6,5,4,0] [6,5,3,0] [6,4,3,0] [5,4,2,0] [5,3,2,0] [4,3,1,0] [4,2,1,0] [3,2,1,0]\\
\hline
3&&[6,5,4,3,2,1] [6,5,4,3,0] [6,5,2,0] [6,4,2,0] [5,3,1,0] [5,2,1,0] [4,3,1,0] [4,2,1,0]\\
\hline
4&&[6,5,4,3,2,1] [6,5,4,3,0] [6,5,2,0] [6,4,2,0] [5,3,2,0] [4,3,1,0] [4,2,1,0] [3,2,1,0]\\
\hline
5&&[6,5,4,3,2,1] [6,5,4,3,0] [6,5,2,1,0] [4,3,2,1,0] [6,4,2,0] [5,3,1,0]\\
\hline
6&&[6,5,4,3,2,1] [6,5,4,3,0] [6,5,2,1,0] [6,4,2,0] [5,3,1,0] [4,3,2,0] [3,2,1,0]\\
\hline
7&&[6,5,4,3,2,1] [6,5,4,3,2,0] [6,5,1,0] [6,4,1,0] [5,3,1,0] [4,2,1,0] [3,2,1,0]\\
\hline
{\color{red}8}&$P_2$&[6,{\color{green}5},4,3,{\color{green}2}] [6,{\color{green}5},4,1,{\color{green}0}] [6,{\color{green}3},2,1,{\color{green}0}] [6,5,3,1] [6,4,2,0] [5,4,3,1] [4,3,2,1] [4,2,1,0]\\
\hline
{\color{red}9}&$P_3$&[{\color{green}6},5,4,3,{\color{green}2}] [6,{\color{green}5},4,1,{\color{green}0}] [6,5,3,1] [6,4,3,0] [6,3,1,0] [5,4,2,1] [5,3,2,1] [4,3,2,0] [4,2,1,0] [3,2,1,0]\\
\hline
{\color{red}10}&$P_4$&[{\color{green}6},5,4,3,{\color{green}2}] [6,{\color{green}5},4,1,{\color{green}0}] [6,5,3,1] [6,4,3,0] [6,3,1,0] [5,4,2,1] [5,3,2,1] [4,3,2,1] [4,3,1,0]\\
\hline
{\color{red}11}&$P_5$&[{\color{green}6},5,4,3,{\color{green}2}] [{\color{green}6},5,4,1,{\color{green}0}] [6,5,3,1] [6,4,3,1] [5,4,2,0] [5,3,2,1] [5,2,1,0] [4,3,2,0] [4,3,1,0] [3,2,1,0]\\
\hline
{\color{red}12}&$P_6$&[{\color{green}6},5,4,3,{\color{green}2}] [{\color{green}6},5,4,1,{\color{green}0}] [6,5,3,1] [6,4,3,1] [5,4,2,0] [5,3,2,1] [5,2,1,0] [4,3,2,1] [4,2,1,0]\\
\hline
{\color{red}13}&$P_7$&[{\color{green}6},5,4,3,{\color{green}2}] [6,5,4,1] [6,5,3,1] [6,4,3,0] [6,4,1,0] [6,3,1,0] [5,4,2,1] [5,3,2,0] [5,3,1,0] [5,2,1,0] [4,3,2,0] [4,2,1,0]\\
\hline
{\color{red}14}&$P_8$&[{\color{green}6},5,4,3,{\color{green}2}] [6,5,4,1] [6,5,3,1] [6,4,3,0] [6,4,1,0] [6,3,1,0] [5,4,2,1] [5,3,2,1] [4,3,2,0] [4,2,1,0] [3,2,1,0]\\
\hline
{\color{red}15}&$P_9$&[{\color{green}6},5,4,3,{\color{green}2}] [6,5,4,1] [6,5,3,1] [6,4,3,1] [5,4,2,0] [5,4,1,0] [5,3,2,0] [5,3,1,0] [4,3,2,0] [4,3,1,0]\\
\hline
{\color{red}16}&$P_{10}$&[{\color{green}6},5,4,3,{\color{green}2}] [6,5,4,1] [6,5,3,1] [6,4,3,1] [5,4,2,1] [5,3,2,0] [5,3,1,0] [5,2,1,0] [4,3,2,0] [4,3,1,0] [4,2,1,0]\\
\hline
{\color{red}17}&$P_{11}$&[{\color{green}6},5,4,3,{\color{green}2}] [6,5,4,1] [6,5,3,1] [6,4,3,1] [5,4,2,1] [5,3,2,1] [4,3,2,0] [4,3,1,0] [4,2,1,0] [3,2,1,0]\\
\hline
18&&[6,5,4,3,2] [6,5,4,3,1] [6,5,2,0] [6,5,1,0] [6,4,2,0] [6,4,1,0] [5,3,2,0] [5,3,1,0] [4,3,2,0] [4,3,1,0]\\
\hline
19&&[6,5,4,3,2] [6,5,4,3,1] [6,5,2,1,0] [4,3,2,1,0] [6,4,2,1] [5,3,2,0] [5,3,1,0]\\
\hline
20&&[6,5,4,3,2] [6,5,4,3,1] [6,5,2,1,0] [6,4,2,0] [6,4,1,0] [5,3,2,0] [5,3,1,0] [4,3,2,0] [4,3,1,0]\\
\hline
21&&[6,5,4,3,2] [6,5,4,3,1] [6,5,2,1,0] [6,4,2,1,0] [5,3,2,0] [5,3,1,0] [4,3,2,0] [4,3,1,0]\\
\hline
22&&[6,5,4,3,2] [6,5,4,3,1] [6,5,2,1,0] [6,4,2,1] [5,3,2,0] [5,3,1,0] [4,3,2,0] [4,3,1,0] [4,2,1,0]\\
\hline
23&&[6,5,4,3,2] [6,5,4,3,1] [6,5,2,1,0] [6,4,2,1] [5,3,2,0] [5,3,1,0] [4,3,2,1] [3,2,1,0]\\
\hline
24&&[6,5,4,3,2] [6,5,4,3,1] [6,5,2,1] [6,4,2,0] [6,4,1,0] [6,2,1,0] [5,3,2,0] [5,3,1,0] [5,2,1,0] [4,3,2,0] [4,3,1,0]\\
\hline
25&&[6,5,4,3,2] [6,5,4,3,1] [6,5,2,1] [6,4,2,1] [5,3,2,0] [5,3,1,0] [5,2,1,0] [4,3,2,0] [4,3,1,0] [4,2,1,0]\\
\hline
26&&[6,5,4,3,2] [6,5,4,3,1] [6,5,2,1] [6,4,2,1] [5,3,2,1] [4,3,2,0] [4,3,1,0] [4,2,1,0] [3,2,1,0]\\
\hline
{\color{red}27}&$P_{12}$&[6,5,4,3] [6,5,4,2] [6,5,3,2] [6,4,3,1] [6,4,2,1] [6,3,2,0] [6,3,1,0] [6,2,1,0] [5,4,3,1] [5,4,2,0] [5,4,1,0] [5,3,2,0] [5,3,1,0] [4,2,1,0]\\
\hline
{\color{red}28}&$P_{13}$&[6,5,4,3] [6,5,4,2] [6,5,3,2] [6,4,3,1] [6,4,2,1] [6,3,2,1] [5,4,3,0] [5,4,2,0] [5,3,2,0] [4,3,1,0] [4,2,1,0] [3,2,1,0]\\
\hline
{\color{red}29}&$P_{14}$&[6,5,4,3] [6,5,4,2] [6,5,3,2] [6,4,3,1] [6,4,2,1] [6,3,2,1] [5,4,3,1] [5,4,2,0] [5,4,1,0] [5,3,2,0] [5,3,1,0] [4,2,1,0] [3,2,1,0]\\
\hline
{\color{red}30}&$P_{15}$&[6,5,4,3] [6,5,4,2] [6,5,3,2] [6,4,3,2] [5,4,3,1] [5,4,2,1] [5,3,2,0] [5,3,1,0] [5,2,1,0] [4,3,2,0] [4,3,1,0] [4,2,1,0]\\
\hline
{\color{red}31}&$P_{16}$&[6,5,4,3] [6,5,4,2] [6,5,3,2] [6,4,3,2] [5,4,3,1] [5,4,2,1] [5,3,2,1] [4,3,2,0] [4,3,1,0] [4,2,1,0] [3,2,1,0]\\

				\Xcline{1-3}{1.2pt}
				
			\end{tabular}
		}
		
		\hspace*{0.5cm}
		\caption{Out of the 31 4-dimensional polytopes with 7 facets, there are 16 admissible polytopes that can potentially be realized as combinatorial types of Coxeter hyperbolic polytopes. Note that the polytope $P_1$ is the only pyramid in this family which is a pyramid over the product of three $1$-simplices.}
		\label{table: combinatoricSelction}
	\end{table}

	According to Theorem \ref{Vinberg:thm3.1}, the subgraphs formed by nodes that correspond to facets intersecting at ideal points of $4$-dimensional hyperbolic Coxeter polytopes yield rank $3$ parabolic graphs. Such diagrams and basic related combinatorial data are listed in Table \ref{table:rank3Euclidean}. 
	We say a polytope is \emph{admissible} if its vertex links are either $3$-simplices, simplicial $3$-prisms, or $3$-cubes. All these three polyhedra are simple, namely every vertices is intersected by $3$ facets. Among the total of $31$ $4$-polytopes with $7$ facets, there are $16$ admissible polytopes as shown in Table \ref{table: combinatoricSelction}. For each vertex whose vertex link constitutes a simplicial $3$-prism, the two facets that comprise a parallel pair are marked in green; for each vertex with a vertex link that is a $3$-cube, the three pairs of parallel facets are colored by  green, black, and pink, respectively.
	
		\begin{table}[H]
		{\footnotesize
			\begin{tabular}{c|c|c|c}
				\Xcline{1-4}{1.2pt}
				\hline
				rank $3$ Euclidean Coxeter diagrams & \rule{0pt}{12pt} $\widetilde{A}_3$ / $\widetilde{B}_3$ / $\widetilde{C}_3$ & $\widetilde{A}_2\sqcup\widetilde{A}_1$ / $\widetilde{B}_2\sqcup\widetilde{A}_1$ / $\widetilde{G}_2\sqcup\widetilde{A}_1$ &
				$\widetilde{A}_1\sqcup\widetilde{A}_1\sqcup\widetilde{A}_1$\\
				\hline
				combinatorial type & $3$-simplex& simplicial $3$-prism & $3$-cube \\
				\hline
				\# facets &4 & 5 & 6\\
				\hline
				\# vertices & 4& 6 & 8\\
				
				\Xcline{1-4}{1.2pt}
			\end{tabular}
		}
		
		\caption{Rank $3$ Euclidean Coxeter graphs.}
		\label{table:rank3Euclidean}
	\end{table}

Now, let us exemplify the process of selecting the 16 polytopes as follows:

\begin{enumerate}
	\item In a hyperbolic Coxeter $4$-polytope $P$ with $7$ facets, if there is a vertex  intersected by six facets, the combinatorial type of $P$ can only be that of a cone over a $3$-cube. In this case, precisely one vertex has a vertex link which is a $3$-cube, while the vertex links of all other vertices are $3$-simplices. 
	Moreover, the vertex link of the unique non-simplex vertex must consist of precisely $8$ vertices, which are all simple. Among the first 7 combinatorial types listed in Table \ref{table: combinatoricSelction}, featuring one facet bracket of length $6$, only the first two exhibit a unique non-simplex vertex. In either of these two cases, the vertex link contains exactly $8$ vertices, all of which are simple. Next, we must ascertain whether the vertex link is a $3$-cube. First, we identify the set of vertices that intersect with the facet bracket of the non-simple vertex at precisely three shared entries, and record these common entries. For the polytope $P_1$, the combinatorial type of the vertex link $lk(v_1)$ of the first vertex $v_1$, represented by $[6,5,4,3,2,1]$, can be shown via facet brackets as follows:
	
	\begin{center}
		$lk(v_1)_{P_1}=\{[6,5,4] ,[6,5,3] ,[6,4,2], [6,3,2] ,[5,4,1] ,[5,3,1], [4,2,1] [3,2,1]\}.$
	\end{center}
	Each facet bracket in $lk(v_1)_{P_1}$ yields a vertex of the vertex link, which corresponds in a natural manner to an edge in $P_1$, as depicted in the low-dimensional illustration, Figure \ref{figure:link} (1).
	Similarly, for the polytope $P_2$, we have the vertex link of the first vertex to be
	
	\begin{center} $lk(v_1)_{P_2}=\{[6,5,4],[6,5,3],[6,4,3],[5,4,2],[5,3,2],[4,3,1],[4,2,1],[3,2,1]\}.$
	\end{center}
	 
	 \noindent It is easy to check that only the vertex link $lk(v_1)_{P_1}$ is combinatorial equivalent to a $3$-cube as shown in Figure \ref{figure:link} (2). 
	 \begin{center}
			\begin{figure}[h]
			\scalebox{0.3}[0.3]{\includegraphics {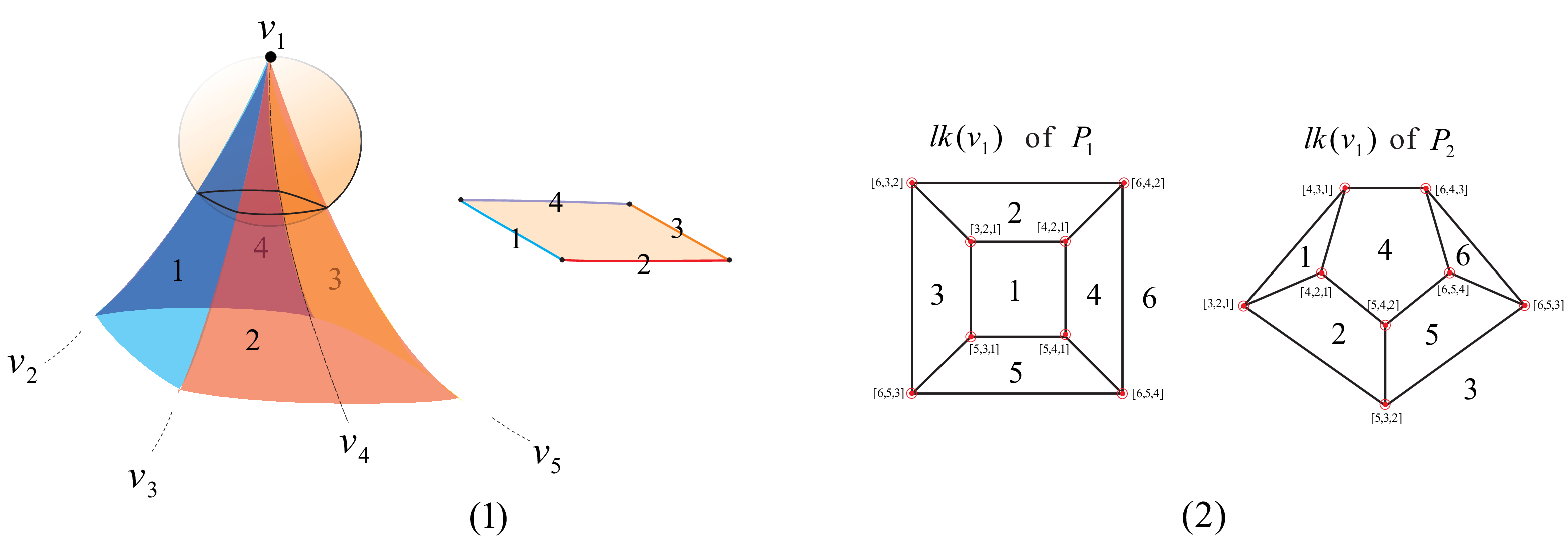}}
			\caption{Vertex links of polytopes.}\label{figure:link}
		\end{figure}
	\end{center}
	\item For a vertex $v$ with facet bracket of length $5$, we want to examine whether its vertex link $lk(v)$ is a simplicial $3$-prism.  Following the approach outlined previously,  we identify the set of vertices whose facet brackets intersect with the facet bracket of $v$ at exactly three shared entries and record the common entries. Subsequently, we exclude vertex links that correspond to cones over 2-cubes, thereby locating all the desired simplicial $3$-prisms.
		
\end{enumerate}

Furthermore, based on the combinatorial information in Table \ref{table: combinatoricSelction}, we can derive the following (1)--(7) \textbf{\emph{data}} for each admissible polytope $P_k$. To illustrate, we provide an example of the \textbf{\emph{data}} for $P_{11}$ in Table \ref{table:p11data}.

\begin{enumerate}
	\item The set $d_k$ consists of pairs of facets that are disjoint.
	\item The set $l_3$ (respectively, $l_4$) comprises sets of three (respectively, four) facets, where the intersections exhibit the combinatorial structure of a triangle (respectively, a tetrahedron).
	\item The set $s_3$ comprises sets of three facets, where the intersections are not edges of $P_k$, and no disjoint pairs are included.
	\item The set $e_3$ consists of sets of three facets of $P_k$, from which the triples in ideal vertices corresponding to Euclidean triangles are excluded, and no disjoint pairs are included.
	\item The set $s_4$ consists of sets of four facets of $P_k$, where the intersections are not vertices of $P_k$, and no disjoint pairs are included.
	\item The set $e_4$ consists of sets of four facets, where the intersections are not vertices of $P_k$. Additionally, it excludes the tetrads in ideal vertices corresponding to Euclidean cubes and does not include any disjoint pairs.
	\item The set $s_5$ (respectively, $s_6$) comprises sets of five (respectively, six) facets, where no disjoint pairs are included.
	\item The set $e_5$ (respectively, $e_6$) consists of sets of five (respectively, six) facets, excluding the 5-tuples (respectively, 6-tuples) corresponding to ideal vertices with vertex links of $3$-prisms (respectively, $3$-cubes). Additionally, no disjoint pairs are included.
	
\end{enumerate}

\begin{table}[h]
	{\scriptsize
			\begin{tabular}{|c|c|l|}
				\Xcline{1-3}{1.2pt}
				\multicolumn{3}{|c|}{\textbf{$P_{11}$ }}\\
				\hline
				\multirow{1}{*}{Vert}& \multirow{1}{*}{10}&  $\{\{{\color{green}6},5,4,3,{\color{green}2}\}, \{6,5,4,1\}, \{6,5,3,1\}, \{6,4,3,1\}, \{5,4,2,1\}, \{5,3,2,1\}, \{4,3,2,0\}, \{4,3,1,0\}, \{4,2,1,0\}, \{3,2,1,0\}\} $\\
							
				\hline
				$d_{11}$ & 2 & $\{ \{0, 5\}, \{0, 6\}\}$ \\
				\hline
				$l_3$& 4 & $\{\{1, 2, 6\}, \{2, 3, 6\}, \{2, 4, 6\}, \{2, 5, 6\}\}$ \\
				\hline
				$l_4$& 4 & $\{\{1, 2, 3, 4\},\{1, 3, 4, 5\}, \{2, 3, 4, 5\}, \{3, 4, 5, 6\}\}$ \\
				\hline
				$s_3$ &5& $\{\{1, 2, 6\}, \{2, 3, 6\}, \{2, 4, 6\}, \{2, 5, 6\}, \{3, 4, 5\}\}$ \\		 
				\hline
				\multirow{3}{*}{$e_3$}& \multirow{3}{*}{$25$}&  $\{\{0, 1, 2\},
				\{0, 1, 3\},
				\{0, 1, 4\},
				\{0, 2, 3\},
				\{0, 2, 4\},
				\{0, 3, 4\},
				\{1, 2, 3\},
				\{1, 2, 4\},
				\{1, 2, 5\},
				\{1, 2, 6\},
				\{1, 3, 4\},	\{1, 3, 5\}
			$\\
				&&$
				\{1, 3, 6\},
				\{1, 4, 5\},
				\{1, 4, 6\},
				\{1, 5, 6\},
				\{2, 3, 4\},
				\{2, 3, 5\},
				\{2, 3, 6\},
				\{2, 4, 5\},
				\{2, 4, 6\},
				\{2, 5, 6\},
				\{3, 4, 6\},\{3, 5, 6\},$\\
				&&$ \{4, 5, 6\}\}$\\
				
				\hline
				\multirow{1}{*}{$s_4$}& \multirow{1}{*}{$10$}&  
				$\{\{1, 2, 3, 4\},
				\{1, 2, 3, 6\},
				\{1, 2, 4, 6\},
				\{1, 2, 5, 6\},
				\{1, 3, 4, 5\},
				\{2, 3, 4, 5\},
				\{2, 3, 4, 6\},
				\{2, 3, 5, 6\},
				\{2, 4, 5, 6\}
				\{3, 4, 5, 6\}\}  $\\
				\hline
				
				\multirow{1}{*}{$e_4$}& \multirow{1}{*}{$10$}&  
				$\{\{1, 2, 3, 4\},
				\{1, 2, 3, 6\},
				\{1, 2, 4, 6\},
				\{1, 2, 5, 6\},
				\{1, 3, 4, 5\},
				\{2, 3, 4, 5\},
				\{2, 3, 4, 6\},
				\{2, 3, 5, 6\},
				\{2, 4, 5, 6\}
				\{3, 4, 5, 6\}\}  $\\
				\hline
				
				$s_5$& 7 &$\{\{0, 1, 2, 3, 4\},
				\{1, 2, 3, 4, 5\},
				\{1, 2, 3, 4, 6\},
				\{1, 2, 3, 5, 6\},
				\{1, 2, 4, 5, 6\},
				\{1, 3, 4, 5, 6\},
				\{2, 3, 4, 5, 6\}\}$\\
				\hline
				$e_5$& 6 &$\{\{0, 1, 2, 3, 4\},
				\{1, 2, 3, 4, 5\},
				\{1, 2, 3, 4, 6\},
				\{1, 2, 3, 5, 6\},
				\{1, 2, 4, 5, 6\},
				\{1, 3, 4, 5, 6\}\}$\\
				\hline
				$s_6$ &1& $\{\{1, 2, 3, 4, 5, 6\}\}$\\
				\hline
				$e_6$ &1& $\{\{1, 2, 3, 4, 5, 6\}\}$\\
				\hline
				
			\end{tabular}
		}
		
		\caption{Combinatorics of the combinatorial polytope $P_{11}$.}
		\label{table:p11data}
	\end{table}

			\section{Block-pasting algorithms for enumerating all the candidate matrices over certain combinatorial types. } \label{chapter:algorithm}
			
			We now use the \emph{block-pasting algorithm} to determine all of the potential matrices for  the $16$   combinatorial types reported in Section \ref{section:4d7f}. Recall that the entries have only finite options, i.e., $\{0,1,2,3,\cdots, 7\}\cup \{\infty\}$. Compared to the backtracking search algorithm raised in \cite{JT:2018}, ``block-pasting" algorithm is more efficient and universal. Generally speaking, the backtracking search algorithm uses the method of ``a series circuit", where the potential matrices are produced one by one. Whereas, the block-pasting algorithm adopts the idea of ``a parallel circuit", incorporating some breadth-first strategy compared to the purely depth-first backtracking search algorithm, where different parts of a potential matrix are generated concurrently and subsequently assembled.

			For each vertex $v_i$ of a 4-dimensional hyperbolic Coxeter polytope $P_k$, we define a ``chunk" $k_i$ consists of an ordered set of facets that intersect at vertex $v_i$ with the subscripts in ascending order. The chunk is called ``ordinary" or ``ideal" if the vertex is ordinary or ideal. For example, for the polytope $P_{1}$ discussed above, there are $8$ ordinary chunks and $1$ ideal chunk as it has $8$ simple vertices and $1$ non-simple vertex with the vertex link a 3-cube. We may also use $k_i$ to denote the ordered set of subscripts, i.e., $k_i$ is referred to as a set of integers of length four, five or six.
			
			Since the vertex links of finite-volume hyperbolic $4$-dimensional polytopes are either simplex, simplicial 3-prism or 3-cube, each chunk possesses $6$, $10$, or $15$ dihedral angles\footnote{Each pair of parallel facets is conventionally possesses  a dihedral angle of $0$.}. For every chunk $k_i$, we define an \emph{i-label} set $e_i$ to be the ordered set $\{10a+b~|~\{a,b\}\in E_i\}$, where $E_i$, named the \emph{i-index} set, is the ordered set of labels of pairs of facets and sorted lexicographically. For example, For the polytope $P_7$ shown in Table \ref{table: combinatoricSelction}, the facets $F_6$, $F_5$, $F_4$, $F_3$, and $F_2$ intersect at the first ideal vertex with vertex link a $3$-dimensional prism, where $F_6$ and $F_2$ are parallel. Then, we have 
			$$k_1=\{F_2,F_3,F_4,F_5,F_6\} ~(\text{or}~ \{2,3,4,5,6\}),$$ 
			$$E_1=
			\{\{2,3\},\{2,4\},\{2,5\},\{2,6\},\{3,4\},\{3,5\},\{3,6\},\{4,5\},\{4,6\},\{5,6\}\},$$ $$e_1=
			\{23,24,25,26,34,35,36,45,46,56\}.$$
			
\begin{figure}[h]
				\scalebox{0.35}[0.35]{\includegraphics {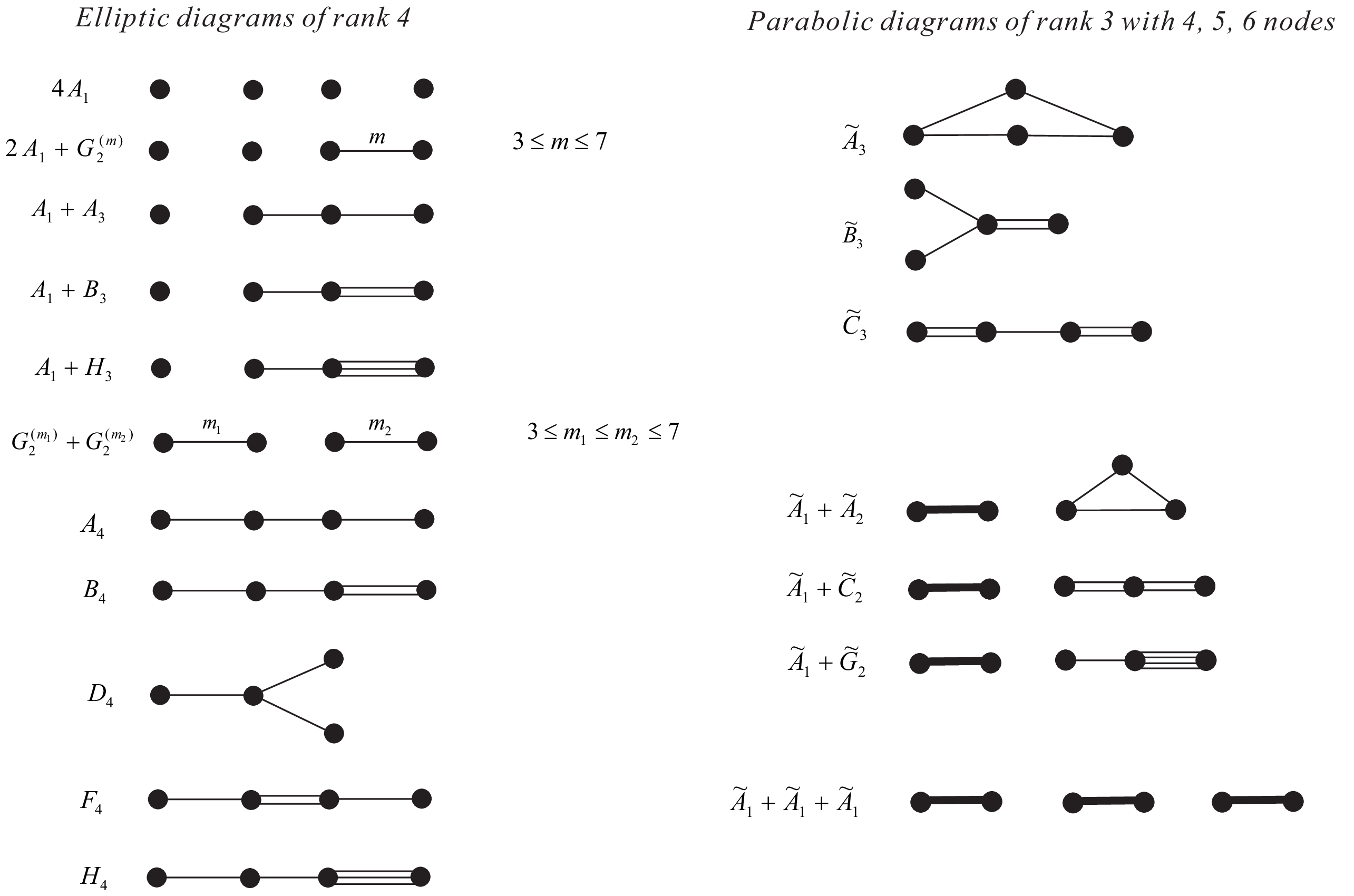}}
				\caption{Spherical Coxeter diagrams of rank $4$ and parabolic Coxeter diagram of rank $3$ with nodes $4$, $5$, and $6$. Only diagrams with integer entries less than or equal to seven are under consideration.}\label{figure:elliptic4}
			\end{figure}

			Next, we list all of the Coxeter vectors of rank $4$ elliptic Coxeter diagrams and rank $3$ parabolic Coxeter diagrams with $4$, $5$, or $6$ nodes as shown in Figure \ref{figure:elliptic4}. Note that the finiteness of the diagrams is a result of our convention to consider only those diagrams with integer weights less than or equal to seven.
			
\begin{figure}[h]
				\scalebox{0.3}[0.3]{\includegraphics {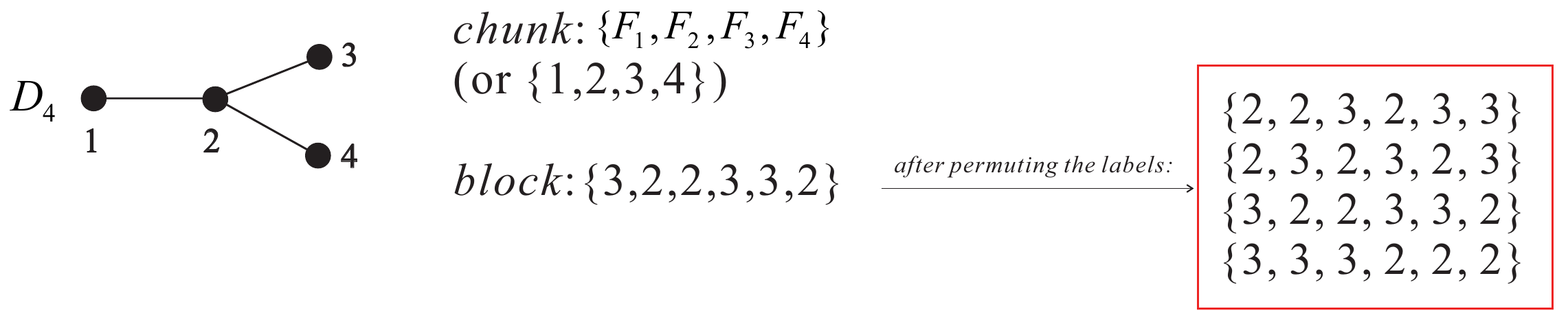}}
				\caption{Prepare the pre-block.}\label{figure:d4}
			\end{figure}

			We apply permutation groups on four, five, and six letters, denoted by $S_4$, $S_5$, and $S_6$ respectively, to the labels of the nodes in the Coxeter diagrams shown in Figure \ref{figure:elliptic4}. This process generates all possible vectors by varying the order of the facets. For instance, there are four vectors associated with the single diagram $D_4$, as depicted in Figure \ref{figure:d4}. The Coxeter diagrams in Figure \ref{figure:elliptic4} yield $269$, $10$, and $1$ distinct vectors corresponding to the order $4$, $5$, and $6$ of diagrams respectively. These vectors are then arranged into sets $\mathcal{S}^4$, $\mathcal{S}^5$, and $\mathcal{S}^6$, referred to as the $4$-preblock, $5$-preblock, and $6$-preblock respectively. These sets can be represented as $269\times 6$, $10\times 10$, and $1\times 15$ matrices in a straightforward manner. In the following, we do not differentiate between these two perspectives and may refer to $\mathcal{S}^*$, where $*=4, 5, 6$, as either a set or a matrix depending on the context.
			
			We proceed by generating a dataframe $B_i$, referred to as the $i$-block, for each chunk $k_i$ of a given polytope $P_k$, where $1\leq i \leq |V_k|$ and $|V_k|$ is the number of vertices of $P_k$. Note that the size of the dataframe depends on whether the vertex link corresponds to a $3$-simplex, a $3$-dimensional prism, or a $3$-dimensional cube, of which the dataframe will be of size $269\times 6$, $10\times 10$, or $1\times 15$ accordingly. 
			
			Firstly we evaluate $B_i$ by suitable $\mathcal{S}^*$ and take the ordered set $e_i$ defined above as the column names of $B_i$. For example, for $e_1=\{12,14,15,24,25,45\}$, the columns of $B_i$ are referred to as $(12)$-, $(14)$-, $(15)$-, $(24)$-, $(25)$-, $(45)$-columns. 
			
			For each admissible polytopes $P_k$, the set of disjoint pairs of facets, denoted by $d_k$ is shown in Table \ref{table:disjoint}. We may also use $d_k$ to represent the label set $\{10a+b~|~{a,b}\in d_k\}$ if there is no ambiguity.
			
\begin{table}[h]
{\tiny
	\begin{tabular}{c|ccc|ccc|cc|c}
	\Xcline{1-10}{1.2pt}
	
	$P_{k}$&1&2&3&4&6&8&9&11&16\\
	\hline
	Disjoint pairs $d_k$ &\{\}&\{\}&\{\}&\{\{0,2\}\}&\{\{0,3\}\}&\{\{0,5\}\}&\{\{0,6\},\{1,2\}\}&\{\{0,5\},\{0,6\}\}&\{\{0,5\},\{0,6\},\{1,6\}\}\\
	\Xcline{1-10}{1.2pt}
	&5&7&12&10&14&&13&15&\\
	\hline
	&\{\}&\{\}&\{\}&\{\{0,6\}\}&\{\{0,6\}\}&&\{\{1,5\},\{0,6\}\}&\{\{0,6\},\{1,6\}\}&\\
	\Xcline{1-10}{1.2pt}
	\end{tabular}
     \caption{ The set $d_k$ of disjoint pairs of facets for each admissible polytope $P_k$. }\label{table:disjoint}
	}
\end{table}
			
			Denote $L$ to be a vector of length $21$ as follows:
			$$L=\{01,02,03,...,06{\color{red},}~12,13,...,16{\color{red},}~23,24,...,26{\color{red},}~34,35,36{\color{red},}~45,46{\color{red},}~56\}.$$
			\noindent Next, we extend each $269\times 6$, $10\times 10$, or $1\times 15$ dataframe to a $k\times 21$ dataframe, where $k=269,~10,~1$, with column names denoted by $L$. This is done by placing each $(ij)$-column in the corresponding labeled column position, filling in the symbol of infinity for the $(ij)$-columns where $\{ij\}\in d_k$, and filling in the value zero for the remaining positions. We continue to use the same notation $B_i$ for the extended dataframe.

		For example, in the case of the polytope $P_8$, the set $d_8$ is $\{0,5\}$, and the vertex link of the first ideal vertex is a $3$-prism with a pair of parallel facets labeled $\{2,6\}$, as shown in Table \ref{table: combinatoricSelction}. Table \ref{table:block1} displays the first block $B_1$ of the polytope $P_8$.
		\begin{table}[h]
			\footnotesize
			\begin{tabular}{c|ccccccccccccccccccccc}
				\hline
				& \textbf{01}& \textbf{02} & \textbf{03} & \textbf{04} & \textbf{05} & \textbf{06} & 
				 \textbf{12} &  \textbf{13} & \textbf{14} & \textbf{15} & \textbf{16} & 
				 \textbf{23} & \textbf{24} & \textbf{25} & \textbf{26} & 
				  \textbf{34} & \textbf{35} & \textbf{36} & 
				  \textbf{45} & \textbf{46} & 
				  \textbf{56} 
				\\
				\hline
				\textbf{1}& *&*&*&*&$\infty$&*&*&*&*&*&*&2&2&2&0&6&3&2&2&2&2\\
				\hline
				\textbf{2} &*&*&*&*&$\infty$&*&*&*&*&*&*&2&2&2&0&4&4&2&2&2&2
				\\
				\hline
				\textbf{3} &*&*&*&*&$\infty$&*&*&*&*&*&*&2&2&2&0&3&6&2&2&2&2
				\\
				\hline
				\textbf{4} &*&*&*&*&$\infty$&*&*&*&*&*&*&2&2&2&0&6&2&2&3&2&2
				\\
				\hline
				\textbf{5} &*&*&*&*&$\infty$&*&*&*&*&*&*&2&2&2&0&3&3&2&3&2&2
				\\
				\hline
				\textbf{6}&*&*&*&*&$\infty$&*&*&*&*&*&*&2&2&2&0&2&6&2&3&2&2
				\\
				\hline
				\textbf{7} &*&*&*&*&$\infty$&*&*&*&*&*&*&2&2&2&0&4&2&2&4&2&2
				\\
				\hline
				\textbf{8} &*&*&*&*&$\infty$&*&*&*&*&*&*&2&2&2&0&2&4&2&4&2&2
				\\
				\hline
				\textbf{9} &*&*&*&*&$\infty$&*&*&*&*&*&*&2&2&2&0&3&2&2&6&2&2
				\\
				\hline
				\textbf{10} &*&*&*&*&$\infty$&*&*&*&*&*&*&2&2&2&0&2&3&2&6&2&2
				\\
				\hline
			\end{tabular}
			\caption{ The first block $B_1$ of the polytope $P_8$. }\label{table:block1}
		\end{table}

			After preparing all of the blocks $B_i$ for a given polytope $P_k$, we proceed to paste them up. More precisely, when pasting $B_1$ and $B_2$, a row from $B_1$ is matched up with a row of $B_2$ where every two entries specified by the same index $i$, where $i\in e_1\cap e_2$, have the same values. The index set $e_1\cap e_2$ is called a \emph{linking key} for the pasting. The resulting new row is actually the sum of these two rows in the non-key positions; the values are retained in the key positions. The dataframe of the new data is denoted by $B_1\cup^*B_2$.
			
			We use the following example to explain this process. Suppose
			
			$B_1=\{x_1,x_2\}=\{(2,2,3,2,3,4,\infty,*,*,*\cdots ,*), (2,3,4,4,2,2,\infty,*,*,*,\cdots,*)\}$,
			
			$B_2=\{y_1,y_2,y_3\}=\{(2,2,3,2,*,*,\infty,4,3,*,\cdots,*),(2,3,4,4,*,*,\infty,2,2,*,\cdots,*)$,
			
			\hspace{3.5cm}$(3,3,2,2,*,*,\infty,3,3,*,\cdots,*)\} .$

			In this example, $x_1$ and $x_2$  have the same values with $y_1$ and $y_2$ on the $(01)$-, $(02)$-, $(03)$- and $(04)$- positions, respectively. In other words, the linking key here is $\{01,02,03,04\}$. Thus, $y_1$ and $y_2$ can paste to $x_1$ and $x_2$, respectively to form the Coxeter vectors $$(2,2,3,2,3,4,\infty,4,3,*,\cdots,*)~\text{and}~(2,3,4,4,2,2,\infty,2,2,*,\cdots,*), \text{respectively}.$$
			\noindent In contrast, $y_3$ cannot be pasted to any element of $B_1$ as there are no vectors with entry $3$ on the $(01)$-position. Therefore, $$B_1\cup^* B_2=\{(2,2,3,2,3,4,\infty,4,3,*,\cdots,*),(2,3,4,4,2,2,\infty,2,2,*,\cdots,*)\}.$$

			We then move on to paste the sets $B_1 \cup^*B_2$ and $B_3$. We follow the same procedure with an updated index set. The linking key is now $e_1\cup e_2\cap e_3$. We conduct this procedure until we finish pasting the final set $B_{\vert V_k\vert}$. The set of linking keys used in this procedure is $$\{e_1\cap e_2, e_1\cup e_2\cap e_3,\cdots,e_1\cup e_2\cup\cdots\cup e_{i-1}\cap e_i,\cdots ,e_1\cup e_2\cup \cdots \cup e_{\vert V\vert -1}\cap e_{\vert V\vert}\}.$$
			
			During the process of pasting, in order to manage the peak volume of resulting vectors and ensure computational feasibility, we have imposed additional restrictions to refine the program. The refined approach introduces additional necessary conditions, apart from the vertex elliptic or parabolic restriction, with the objective of reducing the number of vectors involved in the block-pasting process. 
			
			Firstly, we collect data sets $\mathcal{L}_3$, $\mathcal{L}_4$, $\mathcal{S}_3$, $\mathcal{S}_4$, $\mathcal{S}_5$, $\mathcal{S}_6$, $\mathcal{E}_3$, $\mathcal{E}_4$,
			$\mathcal{E}_5$, and $\mathcal{E}_6$ as claimed in Table \ref{table:library}, by the following two steps: 
			\begin{enumerate}
				\item Prepare Coxeter diagrams with the appropriate rank and node as described in Table \ref{table:library}, and record the Coxeter vectors using an arbitrary system of node labeling.
				\item Apply the permutation group $S_n$ to the node labels and generate the desired dataset, which includes all distinct Coxeter vectors obtained from different labeling systems.
			\end{enumerate}
			
			
			\begin{table}[h]
				{\tiny
					\begin{tabular}{c|c|c|c}
						\Xcline{1-4}{1.2pt}
						\multirow{3}{*}{\textbf{Types of Coxeter diagrams}}
						& \multirow{3}{*}{\textbf{\# Coxeter diagrams }}& \textbf{\# distinct Coxeter } 
						& \multirow{3}{*}{\textbf{data sets}}\\
						&& \textbf{ Vectors after } & \\  
						&& \textbf{ permutation on nodes} & \\ 
						\hline	
								
					Coxeter diagrams of
							hyperbolic $2$-simplices
						& 44 (Lann\'{e}r)& 
						\multirow{2}{*}{175+124=299}
						& 
						\multirow{2}{*}{$\mathcal{L}_3$}\\
						with weights no more than 7 &27 ( quasi-Lann\'{e}r)&&\\
						 				 
						 \hline	
						 
						 \multirow{2}{*}{Coxeter diagrams of
						 	hyperbolic 3-simplices}
						 & 9 (Lann\'{e}r)& 
						 \multirow{2}{*}{108+284=392}
						 & 
						 \multirow{2}{*}{$\mathcal{L}_4$}\\
						 
						 &23 (quasi-Lann\'{e}r)&&\\

						\hline
						rank $3$ elliptic Coxeter diagrams & 9& 31&$\mathcal{S}_3$ \\
						\hline
						rank $4$ elliptic Coxeter diagrams & 29& 242 & $\mathcal{S}_4$ \\
						\hline
						rank $5$ elliptic Coxeter diagrams &47 & 1946 &$\mathcal{S}_5$\\
						\hline
						rank $6$ elliptic Coxeter diagrams & 117& 20206 &$\mathcal{S}_6$\\
						\hline
						parabolic Coxeter diagrams with $3$ nodes (rank $2$)&3 & 10 &$\mathcal{E}_3$ \\
						\hline
						parabolic parabolic Coxeter diagrams with 4 nodes (rank 2(1+1) or 3)& 1+3=4 & 3+27=30 &$\mathcal{E}_4$\\	
						\hline
						parabolic parabolic Coxeter diagrams with 5 nodes (rank 3(1+2) or 4) &3+5=8&100+257=357&
						$\mathcal{E}_5$\\
						\hline
						parabolic Coxeter diagrams with node 6 &\multirow{2}{*}{1+9(3+6)+4=14}&\multirow{2}{*}{15+(405+1000)+870=2290}&\multirow{2}{*}{$\mathcal{E}_6$}\\
						(rank 3(1+1+1), 4((1+3) or (2+2)), or  5)&&& \\
						\Xcline{1-4}{1.2pt}
					\end{tabular}
				}
				
				\caption{Data sets used to reduce the computational load.}
				\label{table:library}
			\end{table}


			Next, we modify the block-pasting algorithm by using additional metric restrictions. More precisely, remarks \ref{remark:1}--\ref{remark:3}, which are practically reformulated from Theorem \ref{Vinberg:thm3.1}, must be satisfied. 
			
			\begin{remark}(``(quasi-)Lann\'{e}r-condition") \label{remark:1}
				The connected Coxeter vector of the three/six dihedral angles formed by the three/four facets with the labels indicated by the data in $l_3$/$l_4$ is {\color{red} IN} $\mathcal{L}_3$/$\mathcal{L}_4$.
			\end{remark}
			
			\begin{remark}(``spherical-condition") \label{remark:2}
				The Coxeter vector of the three/six/ten/fifteen dihedral angles formed by the three/four/five/six facets with the labels indicated by the data in $s_3$/$s_4$/$s_5$/$s_6$ is {\color{red} NOT IN} $\mathcal{S}_3$/$\mathcal{S}_4$/$\mathcal{S}_5$/$\mathcal{S}_6$.
			\end{remark}
			
			\begin{remark}(``Euclidean-condition") \label{remark:3}
				The Coxeter vector of the three/six/ten/fifteen dihedral angles formed by the three/four/five/six facets with the labels indicated by the data in $e_3$/$e_4$/$e_5$/$e_6$ is {\color{red} NOT IN} $\mathcal{E}_3$/$\mathcal{E}_4$/$\mathcal{E}_5$/$\mathcal{E}_6$.
			\end{remark}
			

			The {\color{red}``IN"} and {\color{red}``NOT IN"} tests are called the ``saving" and the ``killing" conditions, respectively. The  ``saving" conditions are significantly more efficient than the ``killing" conditions due to the greater restrictiveness of determining ``qualified" vectors compared to identifying ``non-qualified" vectors.

			We now program these conditions and insert them into appropriate layers during the pasting to reduce the computational load. Here the ``appropriate layer" means the layer where the dihedral angles are non-zero for the first time. For example, for $\{1,2,3\}\in e_3$, we find that after the $j$-th block pasting, the data in columns ($12$-,$13$-,$23$-) of the dataframe $B_1\cup^*B_2\cdots \cup^*B_j$ become non-zero. Therefore, the $e3$-condition for $\{1,2,3\}$ is inserted immediately after the $j$-th block pasting. 
			
		Once the pasting is complete, we proceed to identify and classify all connected Coxeter vectors. We apply symmetry equivalency to classify Coxeter vectors and select representatives from each group. Specifically, two Coxeter vectors, $c_1$ and $c_2$, are considered equivalent, denoted as $c_1 \sim c_2$, if $c_1 = \lambda \cdot c_2$. Here, we regard $c_2$ as a row of data entries with column names represented by $L$. The left action of $\lambda$ on $c_2$ is rearranging the entries of $c_2$ based on the new column names $L_{\lambda}$ obtained through a permutation from the symmetry group $S_7$ operating on the set $\{0,1,2,3,4,5,6\}$.
		
		For example, suppose $\lambda=(01)(234)(56)$, and $c_1$ is as claimed in Table \ref{table:symmetry}, we have $c_1=\lambda\cdot c_2$ where $c_2$ is shown in Table \ref{table:symmetry}.
		\begin{table}[h]
			\footnotesize
			\begin{tabular}{c|ccccccccccccccccccccc}
				\hline
				$L$ & \textbf{01}& \textbf{02} & \textbf{03} & \textbf{04} & \textbf{05} & \textbf{06} & 
				\textbf{12} &  \textbf{13} & \textbf{14} & \textbf{15} & \textbf{16} & 
				\textbf{23} & \textbf{24} & \textbf{25} & \textbf{26} & 
				\textbf{34} & \textbf{35} & \textbf{36} & 
				\textbf{45} & \textbf{46} & 
				\textbf{56} 
				\\
				\hline
				$c_1$& 2& 2& 2& 2& $\infty$& $\infty$& 3& 2& 5& 2& 2& 2& 2& 2& 0& 3& 6& 2& 2& 2& 2\\
				\hline
				$L_\lambda$ & \textbf{01}& \textbf{13} & \textbf{14} & \textbf{12} & \textbf{16} & \textbf{15} & 
				\textbf{03} &  \textbf{04} & \textbf{02} & \textbf{06} & \textbf{05} & 
				\textbf{34} & \textbf{23} & \textbf{36} & \textbf{35} & 
				\textbf{24} & \textbf{46} & \textbf{45} & 
				\textbf{26} & \textbf{25} & 
				\textbf{56} 
				\\
				\hline
				$c_2$& 2& 5& 3& 2& 2& 2& 2& 2& 2& $\infty$& $\infty$& 2& 3&  2& 2 & 2&0&2& 2&6&2\\
				\hline
			\end{tabular}
		\hspace*{0.5cm}
		\caption{Example of symmetry equivalency.}
		\end{table}\label{table:symmetry}

		The matrices (or vectors) after all these conditions, i.e., metric restrictions, connectivity test and symmetry equivalence, are called  \emph{``SELCper"-potential matrices (or vectors)} of certain combinatorial types. The abbreviation \textbf{SELCper} stands for a method that applies \textbf{s}pherical, \textbf{E}uclidean, \textbf{L}anner, and \textbf{c}onnectivity constraints, and results are representatives of equivalency when factoring out the \textbf{per}mutaiton group action on labelling This approach has been Python-programmed. The machine is equipped with Windows 11 Home. Its processor is the 11th Gen Intel(R) Core(TM) i7-1185G7, with four
		computing cores of a 3.0 GHz clockspeed, and the RAM is 32GB. The statistics of the results are reported in Table \ref{table:gall}. The \emph{refined $f$-vector} $((a_1,a_2,a_3),a_4,a_5,a_6,a_7)$ of a polytope $P_k$ provides that the $f$-vector of $P_k$ is $(a_1+a_2+a_3,a_4,a_5,a_6,a_7)$ and it specifies the counts of vertices with vertex links of $3$-cubes, simplicial $3$-prisms, and simplices as $a_1$, $a_2$, and $a_3$, respectively. It is worthy to remark that we have conducted many experiments concerning the gluing part, examining various factors such as the sequence of the gluing process and the possible implementation of more advanced gluing techniques. However, these refinements did not significantly reduce our running time. The current computational time is considered acceptable. Therefore, the blocks are just pasted in a sequential manner. We have documented the running times in Table 9. The notation $s$, $m$, and $h$ on the rows labeled \emph{time} correspond to second, minute, and hour, respectively.

			\begin{table}[h]
				{\scriptsize
					\begin{tabular}{c|c|c|c|c|c|c|c|c|c|c}
						\Xcline{1-11}{1.2pt}
						~~~~\textbf{label}~~~~ &~~ ~~\#~$d_k$~~~~ &refined $f$-vector& ~~~~ \# \textbf{SELCper} ~~~~&time& &~~~~\textbf{label}~~~~&~~~~$d_k$~~~~&refined $f$-vector&~~~~\# \textbf{SELCper} ~~~~&time\\
						\Xcline{1-5}{0.8pt}\Xcline{7-11}{0.8pt}
						1&	0& ((1,0,8),20,18,7,1)	&13	&45.1s&&       9& 2&((0,1,9),21,18,7,1)&	548&7.6m\\
						\Xcline{1-5}{0.8pt}\Xcline{7-11}{0.8pt}
					    2&	0&	((0,3,5),19,18,7,1) &2  &1.6s&&		10&	1&((0,1,10),23,19,7,1) &	289&1.3m\\
					    \Xcline{1-5}{0.8pt}\Xcline{7-11}{0.8pt}
					    3&	0&	((0,2,8),22,19,7,1) &47 &17.9s&&		11&	2& ((0,1,9),21,18,7,1)&	13&5.8s\\
					    \Xcline{1-5}{0.8pt}\Xcline{7-11}{0.8pt}
					    4&	1& ((0,2,7),20,18,7,1)  &16 &6.8s&&		12&	0&((0,0,14),28,21,7,1) & 1,347&6.0h	\\
					    \Xcline{1-5}{0.8pt}\Xcline{7-11}{0.8pt}
					    5&	0&	((0,2,8),22,19,7,1))&18 &9.5s&&		13&	2& ((0,0,12),24,19,7,1)&	10,651&17.6h\\
					    \Xcline{1-5}{0.8pt}\Xcline{7-11}{0.8pt}
					    6&	1&	((0,2,7),20,18,7,1) &37 &41.4s&&		14&	1&((0,0,13),26,20,7,1)&4,992&1.7h\\
					    \Xcline{1-5}{0.8pt}\Xcline{7-11}{0.8pt}
					    7&	0&	((0,1,11),25,20,7,1)&303&3.3m&&		15&	2&((0,0,12),24,19,7,1)&583&7.9m\\
					    \Xcline{1-5}{0.8pt}\Xcline{7-11}{0.8pt}
					    8&	1&	    ((0,1,10),23,19,7,1) &1,541&9.2m&&		16&	3&((0,0,11),22,18,7,1)&	81&1.7m\\
						
						\Xcline{1-11}{1.2pt}
					\end{tabular}
				}
				
				\hspace*{0.5cm}
				\caption{Results of ``SELCper"-potential matrices. Recall that $d_k$ in the table means the number of disjoint pairs as defined before.}
				\label{table:gall}
			\end{table}

			\section{Signature Constraints of hyperbolic Coxeter \texorpdfstring{$n$}-polytopes}\label{section:signature}
			
			After preparing all of the SELCper matrices, we proceed to calculate the signatures of the potential Coxeter  matrices to determine if they lead to the Gram matrix $G$ of an actual hyperbolic Coxeter polytope. 
			
			Firstly, we modify every SELCper matrix $M$ as follows: 
			
			\begin{enumerate}
				\item replace $\infty$s by length unknowns $-x_i$;
				\item replace $0$ by $-1$;
				\item replace $2$, $3$, $4$, $5$, and $6$ by $0$, $-\frac{1}{2}$, $-\frac{l}{2}$, $-\frac{m}{2}$, $-\frac{n}{2}$, where $$l^2-2=0,~l>0,~m^2-m-1=0,~m>0,~n^2-3=0,~n>0;$$
				\item replace $7$s by angle unknowns of $-\frac{y_i}{2}$.
			\end{enumerate} 
			
			By Theorem \ref{thm:signature}, the resulting Gram matrix must have the signature $(4, 1)$. This implies that the determinant of every $6\times 6$ minor of each modified $7\times 7$ SELCper matrix is zero. Therefore, we have the following system of equations and inequality on $x_i$, $l$, $m$, $n$, and $y_i$ to further restrict and lead to the Gram matrices of the desired polytopes:
			$$(6.1) ~~~
			\begin{cases}
				2\det (M_i)=0,~\text{for any of the~}\tbinom{7}{6}=7 ~6\times 6~\text{minor}~M_i~\text{of}~ M. \\
				1.8<y_i<2, ~\text{for~ all}~y_i.\\
				x_i>1,~\text{for~ all}~x_i.\\
				l^2-2=0,~l>0,~m^2-m-1=0,~m>0,~n^2-3=0,~n>0.
			\end{cases}$$
			
			The above conditions are initially stated by Jacquemet and Tschantz in \cite{JT:2018}. Due to practical constraints in \emph{Mathematica}, we denote $2cos(\frac{\pi}{4}),~2cos(\frac{\pi}{5}),~2cos(\frac{\pi}{6})$ by $\frac{l}{2},~\frac{m}{2},~\frac{n}{2}$, rather than $l,~m,~n$ and set $2\det(M_i)=0$ rather than $\det(M_i)=0$. Specific reasons for doing so can be found in \cite{JT:2018}. Moreover, we first find the \emph{Gr\"{o}bner bases} of the polynomials involved, i.e., $2\det(M_i),~l^2-2,~m^2-m-1,~ n^2-3$, before solving the system. This might help to quickly pass over the cases that have no solution. Last but not least, we need to check whether the signature is indeed $(4,1)$ and whether $  \pi/\arccos(\frac{y_i}{2})$ is an integer among the pre-result set. In practice, to expedite the computation process, we have conducted numerous attempts and have identified the following efficient four-step strategy. All intermediate data is presented in Table \ref{table:Mathmaticaprocedure}.
			
			\begin{enumerate}
			
			\item \emph{One equation killing}
			
			\noindent For cases where the number of SELCper matrices exceeds 100, we select equations corresponding to $6\times 6$ minors that exclude length unknowns. We use each equation together with the inequalities corresponding to the unknowns left in the minor as a condition set and solve under a time constraint of 1s.\\

			\noindent The results are categorized into ``out set" (where a non-empty solution is found), ``left set" (where the solution search is terminated due to time constraints), and ``break set" (where no solution exists). We proceed to the second round with the SELCper matrices whose results are in either  ``out set" or ``left set".\\
			
			\item \emph{Seven equations killing}
			
			\noindent We apply the condition system (6.1) to the the SELCper matrices that pass the first round, with a 10s time constraints, engaging all 7 equations. The results also divided into ``out / left / break" sets. In practice, no cases remain unsolved in the ``left set", and we retain the ``out set" for the third round. For cases with 100 or fewer SELCper matrices, we bypass the first step and implement this step directly.\\
			
				\item \emph{Integrality Analysis}
			
			\noindent What we expect about the solution of an angle unknown $y_i$ is of the form $ 2\cos\frac{\pi}{n}$, where $n\in\mathbb{Z}_{\geq 7}$. That is, $  \pi/\arccos(\frac{y_i}{2})$ must be an integer greater than or equal to 7.

			For SELCper matrices that yield a single solution, we verify whether the solution for the angle unknown $y_i$, if it exists, conforms to the form $ 2\cos\frac{\pi}{n}$. If no angle unknowns are present, the solution automatically passes this check;

			If a solution from a SELCper matrix in \emph{Mathematica} is flagged by the warning message
			 {\scriptsize \emph{``Solve::svars:The equations may not give solutions for all variables"}}, we refer to this type of solution as a \emph{vague solution}.  In such cases, we proceed by analyzing the defining equations that pertain only to angles, where the upper bound for at least one of the angle variable $y_i$ is strictly less than $2$. With this constraint, we can exhaust all potential solutions for the angle unknowns due to the integrality restriction.
				
		 Refer to the Example \ref{exampleP15} about the calculation of polytope $P_{15}$ for an illustration. \\
			
			\item \emph{Signature checking}
			
			\noindent The final step entails verifying that the signature is indeed $(4,1)$.
					
		\end{enumerate}

			\begin{example}\label{exampleP15}
				In the computation for the second round concerning polytope $P_{15}$, we find that $4$ SELCper matrices yield vague solution sets, while $3$ present singleton solutions. For example, the Coxeter vector $\{2, 2, 2, 7, 2, \infty, 7, 2, 2, 2, \infty, 3, 2, 2, 2, 2, 5, 2, 3, 2, 2\}$ results in the following system (6.2) of conditions and is one of the four cases with vague solution sets.
		$$(6.2) ~~~
		\begin{cases}
			16 m^2 x_2^2+4 m^2 y_2^2-16 m^2-36 x_2^2-12 y_2^2+36=0\\
			16 m^2 x_1^2+4 m^2 y_1^2-16 m^2-36 x_1^2-12y_1^2+36=0\\
			16 m^2 x_1^2-4 m^2 x_2^2 y_1^2+16 m^2 x_2^2+4 m^2 y_1^2-16 m^2-48 x_1^2+16 x_2^2 y_1^2-48x_2^2-16 y_1^2+48=0\\
			12 x_1^2 y_2^2-48 x_1^2+16 x_2^2 y_1^2-48 x_2^2+4 y_1^2 y_2^2-16y_1^2-12 y_2^2+48=0\\
			-4 m^2 x_1^2 y_2^2+16 m^2 x_1^2+16 m^2 x_2^2+4 m^2 y_2^2-16 m^2+16 x_1^2 y_2^2-48x_1^2-48 x_2^2-16 y_2^2+48=0\\
			16 x_1^2 y_2^2-48 x_1^2+12 x_2^2 y_1^2-48 x_2^2+4 y_1^2 y_2^2-12 y_1^2-16 y_2^2+48=0\\
			-m^2 y_1^2 y_2^2+4 m^2 y_1^2+4 m^2 y_2^2-16 m^2+4 y_1^2 y_2^2-12y_1^2-12 y_2^2+36=0\\
			m^2-m-1=0\\
			x_1,x_2>1\\
			1.8<y_1,y_2<2\\
			m>0		
		\end{cases}$$

			\noindent The Gr\"{o}bner bases of the system (6.2) is as shown in (6.3) and a vague message showed up when we solve it in \emph{Mathematica}.
				\noindent
			
				$$(6.3) ~~~
				\begin{cases}
					24-12y_1^2-12y_2^2+5 y_1^2 y_2^2 - 28 m + 4 y_1^2 m +4 y_2^2 m = 0\\
					-11 + 11 x_2^2 + 2 y_2^2 - 3 y_2^2 m = 0\\
					-11 + 
					11 x_1^2 + 2 y_1^2 - 3 y_1^2 m = 0\\
					m^2-m-1=0\\
					x_1,x_2>1\\
					1.8<y_1,y_2<2\\
					m>0		
				\end{cases}$$				\noindent We then derive from the system a unique defining relation for $y_2$, expressed as \noindent
			{
				$$y_2=f(y_1), 1.80194\leq y_1\leq 1.91472,$$} \noindent where $f$ is a function with $y_1$ as its variable. We go over the set \{$2\cos\frac{\pi}{7}$,$2\cos\frac{\pi}{8}$,$2\cos\frac{\pi}{9}$,$2\cos\frac{\pi}{10}$\} for possible value of  $h_1$ to substitute into the function $f(y_1)$, and have found that none of the function values conform to be form $2\cos\frac{\pi}{n}$.
			
			Besides, the Coxeter vector $\{2, 3, 2, 4, 2, \infty, 2, 7, 2, 2, \infty, 4, 2, 3, 2, 2, 2, 2, 5, 2, 2\}$ is one of the three instances with a singleton solution. The solution is in (6.4). It is evident that the solution of $y_1$ is not of the form $2\cos\frac{\pi}{n}$.
			 	\noindent
			 $$(6.4) ~~~
			 	\begin{cases}
			 		 x_1=\sqrt{\frac{3}{2}+\sqrt{5}+\sqrt{7+3 \sqrt{5}}}\\
			 		 x_2=\sqrt{\frac{1}{31} \left(-70+92 \sqrt{2}+53 \sqrt{5}-36 \sqrt{10}\right)}\\
			 		 \displaystyle y_1=\frac{2}{\sqrt{7+4 \sqrt{5}-2 \sqrt{2} \left(3+\sqrt{5}\right)}}\\
			 		 m=\frac{1}{2}\left(1+\sqrt{5}\right)\\p=\sqrt{2}
			 		
			 	\end{cases}$$
			 
			 \qedhere{\hfill\ensuremath{\square}}

			\end{example}
		
	This approach has been implemented in Mathematica, and we have discovered that fourteen of the admissible $4$-polytopes with $7$ facets can be realized as hyperbolic Coxeter polytopes. The statistics is presented in Table \ref{table:result}. The Coxeter vectors and related information are in Table \ref{table:resultall}--\ref{table:resultall3}, where volumes can be calculated via CoxIter \cite{Guglielmetti:2017} based on Proposition \ref{prop:volume}.  
	
	\begin{remark}
The distance data for non-incident facets can be very complicated. Extensive efforts were made to render the data in a more comprehensible format, as presented in Tables 13--14. While the information remains intricate, it is our best attempt to simplify it using \emph{Mathematica}; at least readers can identify the minimal number fields that contain the length hyperbolic cosines. The signature verification  is time-consuming, even in a cluster of servers. For readers interested in doing that, it might be better to solve for the length values based on the potential vectors provided in our Github (\emph{https://github.com/GeoTopChristy/HCPdm/tree/}
	
\noindent\emph{main/output/potential$\_$vectors$\_$ip47$\_$all}) using their preferred algebraic computation software. This might save time compared to manually entering the data provided here. While not ideal, a rapid check for zeros in the signatures can be achieved by employing the numerical results with a sufficiently high level of precision.

	\end{remark}
	
	\begin{table}[H]
	{\scriptsize
		\begin{tabular}{c|c|c|c|c|c}
			\Xcline{1-6}{1.2pt}
			label	&	\#SELCper	&	1st round	&	2nd round	&	3rd round	&	4th round	\\
			\Xcline{1-6}{1.2pt}
			1&13&N&13&13&13\\
			2&2&N&2&2&2\\
			3&47&N&1&1&1\\
			4&16&N&16&16&16\\
			5&18&N&2&1&1\\
			6&37&N&37&37&37\\
			7&303&52 [7]&5&2&2\\
			8&1541& 11 [6]&11&11&11\\
			9&548&134 [7]&134&134&134\\
			10&289&63 [7]&63  &4&4\\
			11&13&N&13&13&13\\
			12&1347&220 [7] -- 43 [6] -- 16 [5] -- 7 [4]&3&0&0\\
			13&10651&3776 [2]&31+10=41&8&8\\
			14&4992&966 [7] -- 275 [1]&72&8&8\\
			15&583&7 [7]&3+4=7&0&0\\
			16&81&N&81&81&81\\
			\Xcline{1-6}{1.2pt}
		\end{tabular}
	}
	
	\hspace*{0.5cm}
	\caption{In first-round column, the letter ``N" is used to indicate cases where the count of SELCper matrices is small, thereby bypassing this step. For all other entries in this column, the data is presented in the format $i~[j]$, where $i$ is the number of the cardinal of the union of out set and left set, and $j$ indicates the specific row and column excluded to form the $6\times 6$ minor. For example, for the polytope $P_{14}$, the set $d_{14}$ of disjoint pairs, as listed in Table \ref{table:disjoint}, is $\{\{0,6\}\}$. This means that as long as the first and seventh rows and columns are not included together in a minor, no length variable will be included.  We first exclude row 7 and column 7 to form a minor $M_7$ and use the equation $\det M_7=0$ along with related inequalities to filter the qualified SELCper matrices, resulting in 966 from the original 4992. Similarly, we consider the minor formed by deleting the first row and column, finally decreasing to 275.  In the second-round column, we report the number of results obtained after the second step. Specifically, for polytopes $P_{13}$ and $P_{15}$, we present the number in the form $i+j$, where $i$ is the  count of singleton solutions, and $j$ is the count of vague solutions. Vague solutions occur only in these two cases.}
	\label{table:Mathmaticaprocedure}
\end{table}

			\begin{table}[H]
				{\tiny
						\begin{tabular}{c|c|cccc|cccccc|ccc}
							\Xcline{1-15}{1.2pt}
							\# non-simple vertices &	\multicolumn{1}{c|}{3} & \multicolumn{4}{c|}{2} & \multicolumn{6}{c|}{1}& \multicolumn{3}{c}{0}\\
							\Xcline{1-15}{1.2pt}
							
							polytope labels&2&3&4&5&6&1&7&8&9&10&11&13&14&16\\
								\Xcline{1-15}{1.2pt}
													
							\# finite volume but not compact&	2&1&16&1&37&13&2&11&134 &4 &13&5&0&52\\	
							\hline
							\# compact  &N& N&N&N&N & N&N&N&N&N&N& 3&8&29\\
							\hline
							Roberts's list in dimension $4$ [Rob15]
&
							 \multirow{3}{*}{N}&\multirow{3}{*}{N}&\multirow{3}{*}{N}&\multirow{3}{*}{N}&\multirow{3}{*}{N}&\multirow{3}{*}{N}&\multirow{3}{*}{0}&\multirow{3}{*}{0}&\multirow{3}{*}{20}&\multirow{3}{*}{4}&\multirow{3}{*}{0}&\multirow{3}{*}{N}&\multirow{3}{*}{N}&\multirow{3}{*}{N}\\
							 (Hyperbolic Coxeter non-pyramidal $4$-polytope  &&&&&&&&&&&&&&\\
							   with $7$ facets and $1$ non-simple vertex) &&&&&&&&&&&&&&\\
							\hline
							
							Tumarkin's list in dimension $4$ [Tum07] [Tum0$4^{(2)}$]
&
							\multirow{3}{*}{N}&\multirow{3}{*}{N}&\multirow{3}{*}{N}&\multirow{3}{*}{N}&\multirow{3}{*}{N}&\multirow{3}{*}{13}&\multirow{3}{*}{N}&\multirow{3}{*}{N}&\multirow{3}{*}{N}&\multirow{3}{*}{N}&\multirow{3}{*}{N}&\multirow{3}{*}{3}&\multirow{3}{*}{8}&\multirow{3}{*}{29}\\
							(Compact and and pyramidal hyperbolic Coxeter  &&&&&&&&&&&&&&\\
							$4$-polytopes with $7$ facets) &&&&&&&&&&&&&&\\
							\Xcline{1-15}{1.2pt}

						\end{tabular}
					}
					
					\hspace*{0.5cm}
					\caption{Statistics of the result.  The symbol ``N" indicates that the combinatorial polytope is not included in the corresponding family. Specifically, for a compact hyperbolic Coxeter polytope to be considered, it must be simple, and thus a combinatorial polytope with non-simple vertices cannot possess the structure of a compact hyperbolic Coxeter polytope as shown on line $4$ and line $6$. Roberts considered non-pyramidal polytopes with one non-simple vertex, and therefore the symbol ``N" is used for other families with three, two, or zero non-simple vertices. 
						}
					\label{table:result}
				\end{table}

\newpage
				
				\section{Validation, Data Availability and Results}\label{section:validation}
\subsection{Validation}		
Hyperbolic Coxeter non-pyramidal $n$-polytopes with $n+3$ facets and only one non-simple vertex has been previously enumerated by Roberts. Our work includes Roberts's findings of $24$ polytopes in dimension $n=4$, as originally presented in \cite{Roberts:15} and depicted in Table \ref{table:resultall} with blue labels ($P_{9,7},\cdots,P_{9,16},P_{9,119},\cdots,P_{9,128},P_{10,1},\cdots,P_{10,4}$). Furthermore, we have discovered $140$ polytopes within this family that were lost from Roberts's list. Summary of these findings is displayed in Table \ref{table:result}.

Secondly, in \cite{Tumarkin: n3fv} and \cite{Tumarkin: n3}, Tumarkin classified all $13$ non-simple pyramidal Coxeter polytope and $40$ compact hyperbolic Coxeter $n$-polytopes with $n+3$ facets, respectively. Our findings coincide with Tumarkin's classification in dimension $4$, which are shown in Table \ref{table:resultall} with orange labels ($P_{1,1},\cdots,P_{1,13}$) and with  red labels ($P_{13,1},\cdots,P_{13,3}$,  $P_{14,1},\cdots,P_{14,8}$,  $P_{16,2},\cdots,P_{16,5}$,  $P_{16,7},\cdots P_{16,12},$ $P_{16,19},\cdots,P_{16,23},$ $P_{16,29},\cdots,P_{16,37},$ 

\noindent $P_{16,41},\cdots,P_{16,45}$), respectively. For statistic information, please refer to Table \ref{table:result}.
				
Thirdly, Felikson and Tumarkin studied finite-volume hyperbolic Coxeter polytopes with mutually intersecting facets in \cite{FT:08}. In particular, for compact hyperbolic Coxeter $4$-polytopes without disjoint facets, the possibilities are restricted to either a simplex or one of the seven Esselmann polytopes, each having exactly six facets. For a non-compact but finite-volume $4$-dimensional simple Coxeter polytope devoid of disjoint facets,  it must be either a simplex or a specific 4-polytope with exactly six facets. The $4$-polytope $P_{15}$, as presented in Table \ref{table: combinatoricSelction}, is simple, comprises  seven facets, and lacks any disjoint pair of facets. We show the absence of a hyperbolic Coxeter structure over this polytope, as detailed in Table \ref{table:result}, which is consistent with the results of Felikson and Tumarkin.

At last, the compactness and volume finiteness of the polytopes have been verified by CoxIter \cite{Guglielmetti:2017} \emph{https://github.com/rgugliel/CoxIter}.
						
\subsection{Data Availability}\label{dataaval}

The result data is completely stated here. The example codes and all the intermediate data is available at 
{\centering
\noindent\emph{ https://github.com/GeoTopChristy/HCPdm}.}

\subsection{Results}\label{dataaval}
The Coxeter vectors for the 331 hyperbolic Coxeter 4-polytopes with seven facets are presented in Tables \ref{table:resultall} -- \ref{table:resultall3}. These tables also provide documentation on their arithmeticity, the number of cusps, and the volume.  We use letters to alphabetically denote  the length values in the tables, rather than $-x_i$s, to make them easier to identify. The corresponding Coxeter diagrams are provided in the appendix.

\newgeometry{left=1cm,right=1cm,top=0.5cm,bottom=0.5cm}

\begin{landscape}

{\centering
	{\scriptsize
		\renewcommand\baselinestretch{1.5}\selectfont
\begin{longtable}{c|c|c|c|c|c|c|c|c|c|cc}

				\Xcline{1-12}{1.2pt}
				$i$ &	000000 & 11111 & 2222 & 333 & 44 & 5 &\multirow{2}{*}{\# cusps} & \multirow{2}{*}{$\frac{\text{volume}}{\pi^2}$}&\multirow{2}{*}{Arithmetic}&\multirow{2}{*}{ $a$} & \multirow{2}{*}{$b$}\\
				$j$ &	123456 & 23456 & 3456 & 456 & 56 & 6 &  &  & && \\
				\Xcline{1-12}{1.2pt}
				
{\color{orange}$P_{1,1}$	}&{\color{orange}	333222	}&{\color{orange}	2222$0$	}&{\color{orange}	22$0$2	}&{\color{orange}	$0$22}&{\color{orange}	22	}&{\color{orange}	2}&{\color{orange}1}&{\color{orange}$\frac{1}{144}$}&{\color{orange}A}&{\color{orange}}&{\color{orange}}\\
{\color{orange}$P_{1,2}$	}&{\color{orange}	333322	}&{\color{orange}2222$0$	}&{\color{orange}22$0$2	}&{\color{orange}	$0$22}&{\color{orange}	22}&{\color{orange}2}&{\color{orange}1}&{\color{orange}$\frac{1}{72}$}&{\color{orange}A}&{\color{orange}}&{\color{orange}}\\
{\color{orange}$P_{1,3}$	}&{\color{orange}333332}&{\color{orange}2222$0$}&{\color{orange}22$0$2}&{\color{orange}$0$22}&{\color{orange}22}&{\color{orange}2}&{\color{orange}1}&{\color{orange}$\frac{1}{36}$}&{\color{orange}A}&{\color{orange}}&{\color{orange}}\\
{\color{orange}$P_{1,4}$	}&{\color{orange}333333}&{\color{orange}2222$0$}&{\color{orange}22$0$2}&{\color{orange}$0$22}&{\color{orange}22}&{\color{orange}2}&{\color{orange}1}&{\color{orange}$\frac{1}{18}$}&{\color{orange}A}&{\color{orange}}&{\color{orange}}\\
{\color{orange}$P_{1,5}$	}&{\color{orange}	433222}&{\color{orange}2222$0$}&{\color{orange}22$0$2}&{\color{orange}$0$22}&{\color{orange}22}&{\color{orange}2}&{\color{orange}2}&{\color{orange}$\frac{1}{72}$}&{\color{orange}A}&{\color{orange}}&{\color{orange}}\\
{\color{orange}$P_{1,6}$}&{\color{orange}433223}&{\color{orange}2222$0$}&{\color{orange}22$0$2}&{\color{orange}$0$22}&{\color{orange}22}&{\color{orange}2}&{\color{orange}2}&{\color{orange}$\frac{1}{48}$}&{\color{orange}}&{\color{orange}}&{\color{orange}}\\
\end{longtable}}
}

{\centering
	{\tiny
		\renewcommand\baselinestretch{1.5}\selectfont
		\begin{longtable}{c|c|c|c}

			\Xcline{1-4}{1.2pt}

	&		a		&		b		&		c		\\
		\Xcline{1-4}{1.2pt}
$P_{16,1}$	&	$	\sqrt{\frac{2}{7}(5+4\sqrt{2})}	$	&	$	\frac{1}{7}(3+8\sqrt{2})	$	&	$	\sqrt{\frac{2}{7}(5+4\sqrt{2})}	$	\\
$P_{16,2}$	&	$	\frac{\sqrt{5+\sqrt{5}}}{2}	$	&	$	\sqrt{\frac{1}{11}(16+7\sqrt{5})}	$	&	$	\sqrt{\frac{1}{11}(10+3\sqrt{5})}	$	\\
$P_{16,3}$	&	$	\frac{1}{2}(1+\sqrt{5})	$	&	$	\frac{1}{2}(1+\sqrt{5})	$	&	$	\frac{\sqrt{3+\sqrt{5}}}{2}	$	\\
$P_{16,4}$	&	$	\frac{(21-53 \sqrt{2}+19 \sqrt{5}+17 \sqrt{10}) \sqrt{1037+731 \sqrt{2}+475 \sqrt{5}+315 \sqrt{10}}}{1364}	$	&	$	\sqrt{\frac{1}{11}(24+6\sqrt{2}+5\sqrt{5}+4\sqrt{10})}	$	&	$	\sqrt{\frac{1}{22}(19+13\sqrt{2}+9\sqrt{5}+5\sqrt{10})}	$	\\
$P_{16,5}$	&	$	\frac{1}{2}\sqrt{11+5\sqrt{5}}	$	&	$	3\sqrt{\frac{1}{31}(6+\sqrt{5})}	$	&	$	\sqrt{\frac{1}{31}(17+8\sqrt{5})}	$	\\
$P_{16,6}$	&	$	\sqrt{\frac{2}{19}(9+\sqrt{5})}	$	&	$	\sqrt{\frac{1}{209} \left(233+104 \sqrt{5}\right)}	$	&	$	\sqrt{\frac{2}{11}(7+\sqrt{5})}	$	\\
$P_{16,7}$	&	$	\sqrt{\frac{1}{38}(43+9\sqrt{5})}	$	&	$	\sqrt{\frac{1}{38}(39+17\sqrt{5})}	$	&	$	\frac{1}{2}\sqrt{\frac{1}{2}(9+\sqrt{5})}	$	\\
\end{longtable}}
}

{\centering
	{\tiny
		\renewcommand\baselinestretch{1.5}\selectfont
		\begin{longtable}{c|c}

			\Xcline{1-2}{1.2pt}

	$	a_{36}	$	&	$\displaystyle	\frac{\sqrt{w} \left(2 \sqrt{3} \left(-45-653 \sqrt{2}+228 \sqrt{5}+252 \sqrt{10}\right) \sqrt{u}+63041 \sqrt{5} \left(2+\sqrt{2}\right) \sqrt{\frac{1}{w}}-22 \left(119-33 \sqrt{2}+57 \sqrt{5}+63
		\sqrt{10}\right)\right)}{504328}+\frac{1}{8} \left(3 \sqrt{2}-2\right)$, where\\
	&$u=89+58 \sqrt{2}+41 \sqrt{5}+24 \sqrt{10}$,\\
&$w=2 \sqrt{30} \sqrt{u}+18 \sqrt{15} \sqrt{u}+14 \sqrt{6} \sqrt{u}+38 \sqrt{3} \sqrt{u}+294 \sqrt{10}+637 \sqrt{5}+650 \sqrt{2}+1489$.	\\
\hline
$	a_{37}	$	&	$\displaystyle	\frac{\sqrt{3} \left(1260-2103 \sqrt{2}+1074 \sqrt{5}+1193 \sqrt{10}\right) \sqrt{u}+11 \left(23221 \left(-\sqrt{2}+\sqrt{5}+\sqrt{10}\right) \sqrt{\frac{1}{w}}+\frac{69663}{\sqrt{w}}-786 \sqrt{10}-1357 \sqrt{5}-1842 \sqrt{2}-3019\right)}{1021724 \sqrt{2}\sqrt{\frac{1}{w}}}$, where\\
&$
u=89+58 \sqrt{2}+41 \sqrt{5}+24 \sqrt{10}$,\\

&$w=18 \sqrt{30} \sqrt{u}+74 \sqrt{15} \sqrt{u}+38 \sqrt{6} \sqrt{u}+166 \sqrt{3} \sqrt{u}+940 \sqrt{10}+1835 \sqrt{5}+2070 \sqrt{2}+4243$.	\\

\hline

$	a_{45}	$	&	$\displaystyle	\frac{\sqrt{w} \left(\sqrt{11} \left(1367+847 \sqrt{2}-9 \sqrt{5} \left(63+43 \sqrt{2}\right)\right) \sqrt{u}+2981 \left(\sqrt{2}+\sqrt{5}+\sqrt{10}\right) \sqrt{\frac{1}{w}}+2 \sqrt{11} \left(-383-227 \sqrt{2}+143 \sqrt{5}+89 \sqrt{10}\right)\right)}{11924}+\frac{1}{4}$, where\\
&$u=11+7 \sqrt{2}+5 \sqrt{5}+3 \sqrt{10}$,\\
\end{longtable}}
}

\end{landscape}
\restoregeometry

\begin{center}
	{\large \textbf{Appendix: Coxeter diagrams for hyperbolic Coxeter $4$-polytopes with $7$ facets.}}
\end{center}

\begin{figure}[h]
	\scalebox{0.47}[0.47]{\includegraphics {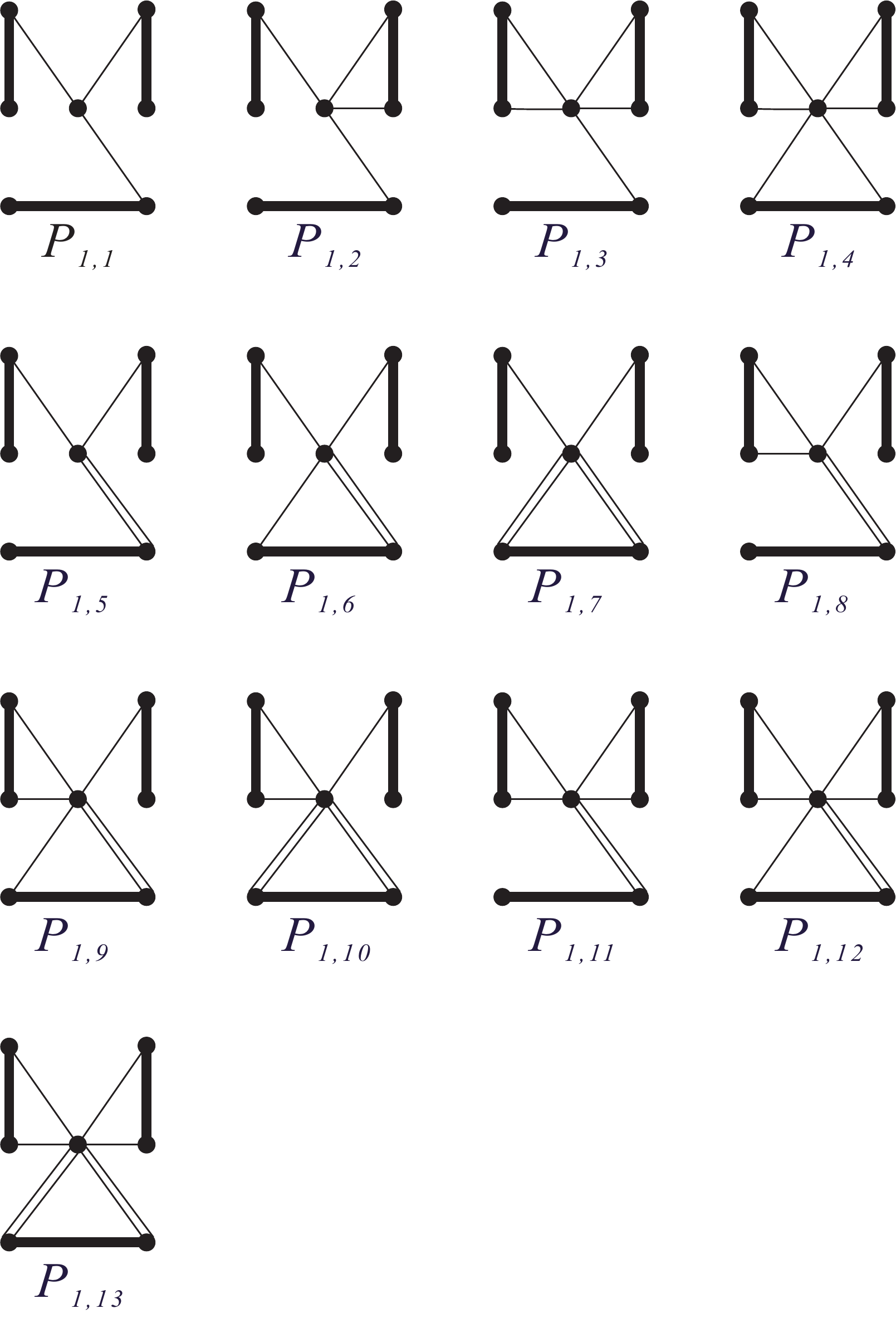}}
	\caption{The $13$ finite-volume hyperbolic Coxeter $4$-polytopes with $7$ facets over polytope   $P_{1}$.} \label{figure:p1}
\end{figure}

1. Coxeter diagrams for $P_{1}$ (Figure \ref{figure:p1}).
\vspace{1cm}

\begin{figure}[H]
	\scalebox{0.47}[0.47]{\includegraphics {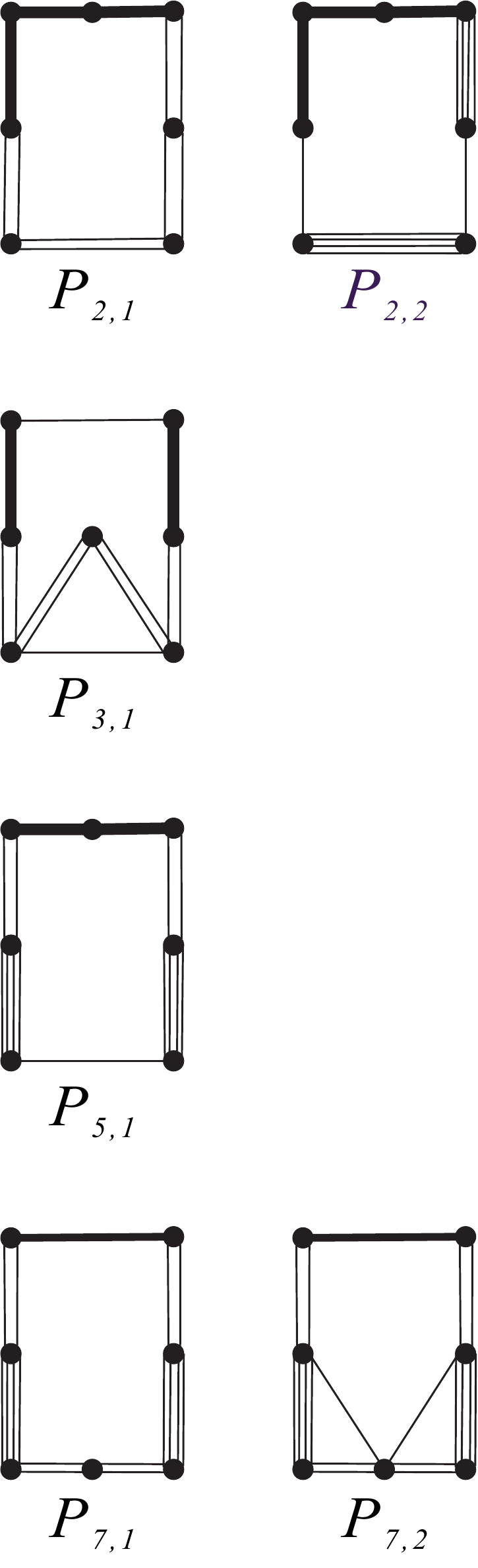} }
	\caption{The 2, 1, 1, 2 finite-volume hyperbolic Coxeter $4$-polytopes with $7$ facets over polytopes $P_2$, $P_3$, $P_5$, and $P_7$, respectively.} \label{figure:p2p3p5p7}
\end{figure}

\newpage 
2.Coxeter diagrams for $P_2$, $P_3$, $P_5$, and $P_7$ (Figure \ref{figure:p2p3p5p7}).

\begin{figure}[H]
	\scalebox{0.5}[0.47]{\includegraphics {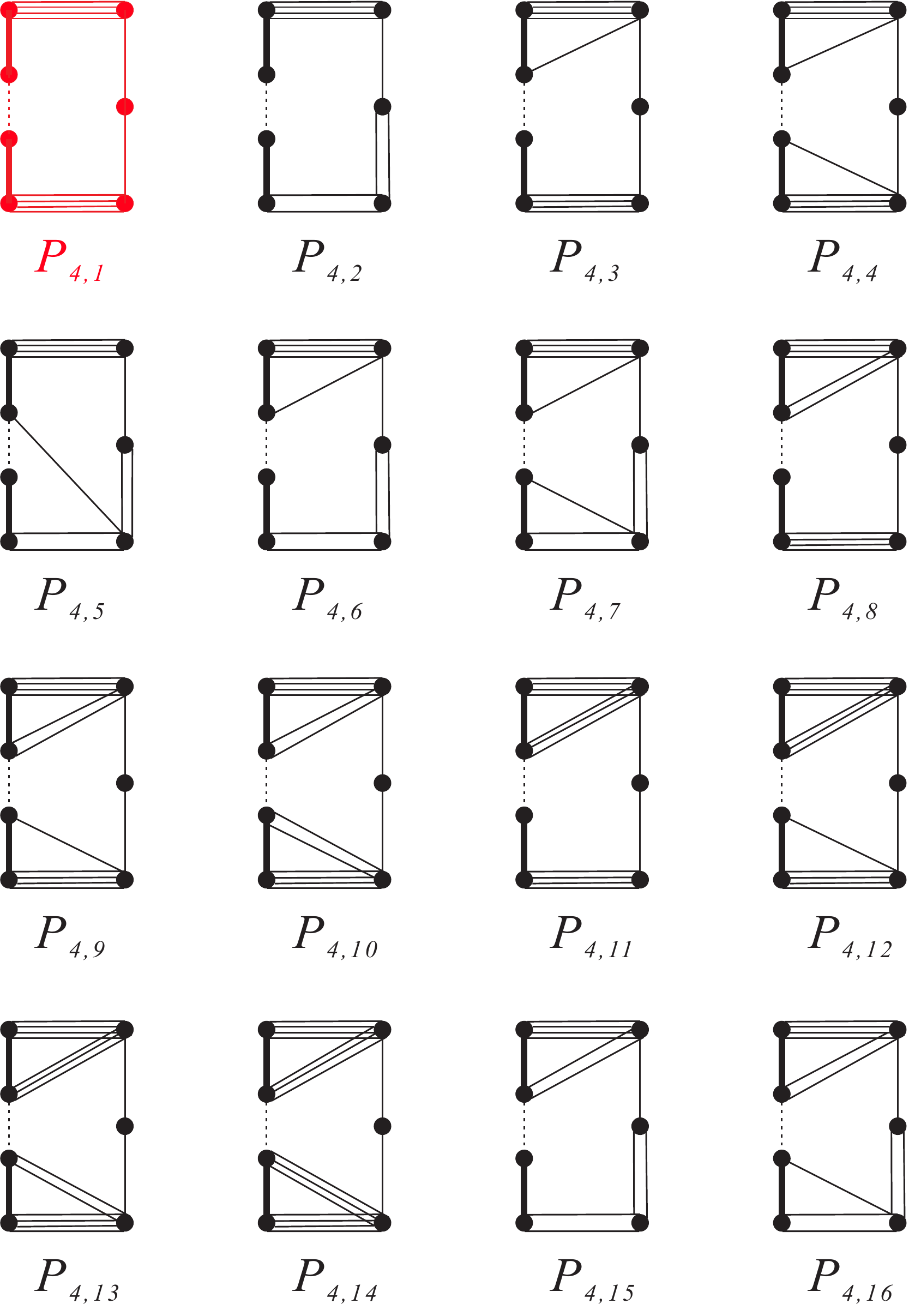}}
	\caption{The $16$ finite-volume hyperbolic Coxeter $4$-polytopes with $7$ facets over polytope $P_{4}$.} \label{figure:p41}
\end{figure}

\newpage
3. Coxeter diagrams for $P_{4}$ (Figure \ref{figure:p41}).

\begin{figure}[H]
	\scalebox{0.47}[0.47]{\includegraphics {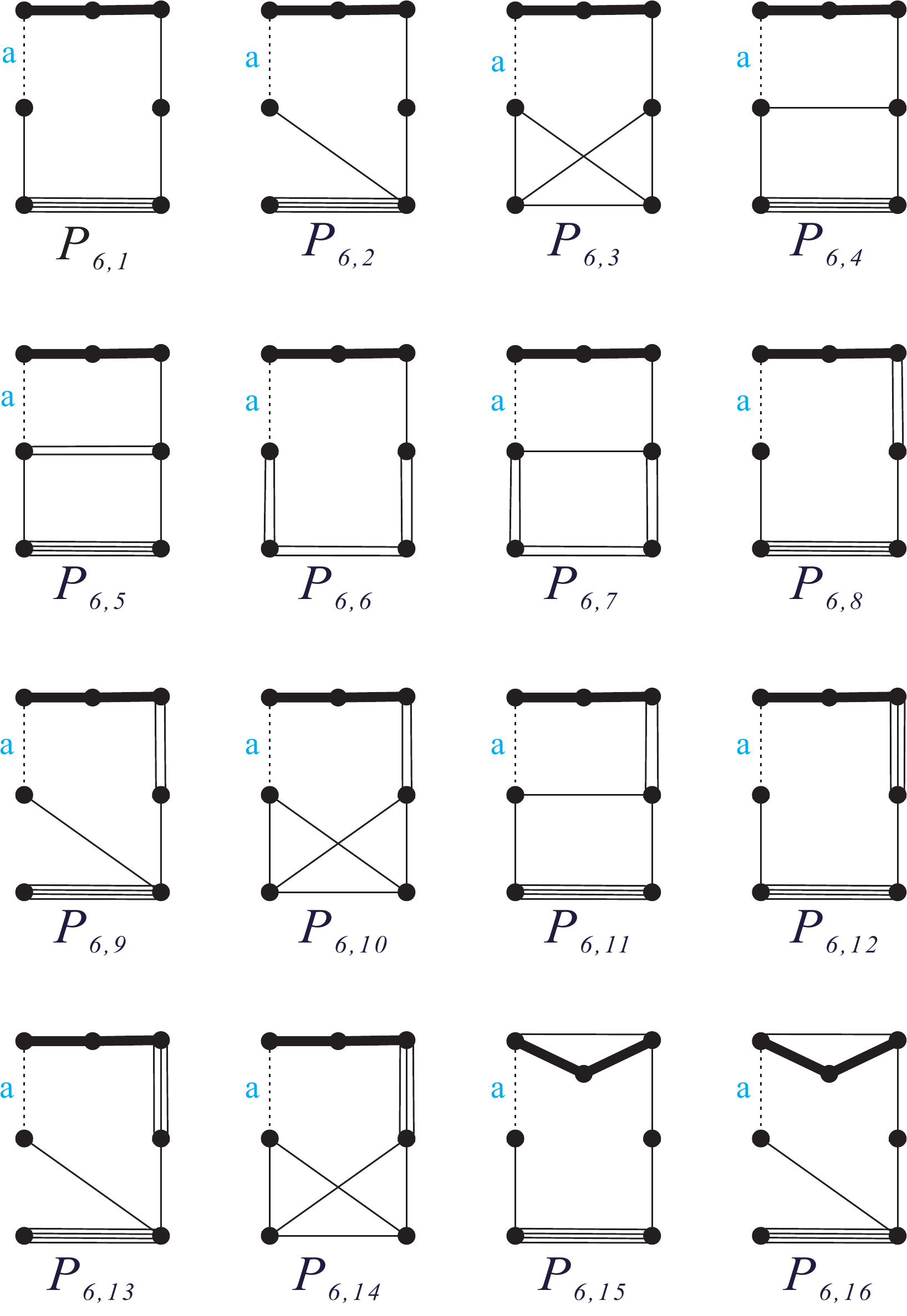}}
	\caption{The $37$ finite-volume hyperbolic Coxeter $4$-polytopes with $7$ facets over polytope $P_{6}$ (part 1).} \label{figure:p61}
\end{figure}

\begin{figure}[H]
	\scalebox{0.47}[0.47]{\includegraphics {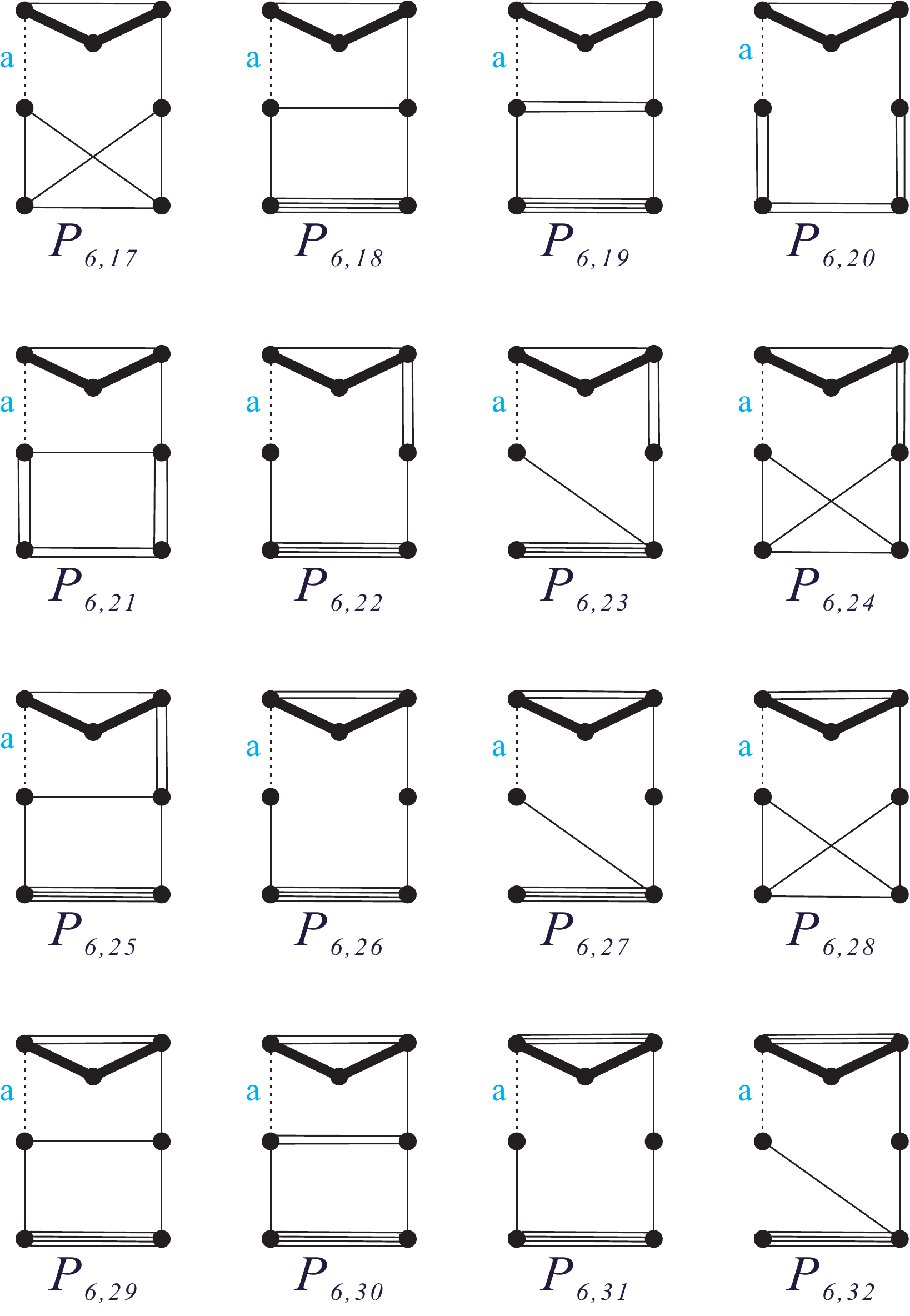}}
	\caption{The $37$ finite-volume hyperbolic Coxeter $4$-polytopes with $7$ facets over polytope $P_{6}$ (part 2).} \label{figure:p62}
\end{figure}

\begin{figure}[H]
	\scalebox{0.47}[0.47]{\includegraphics {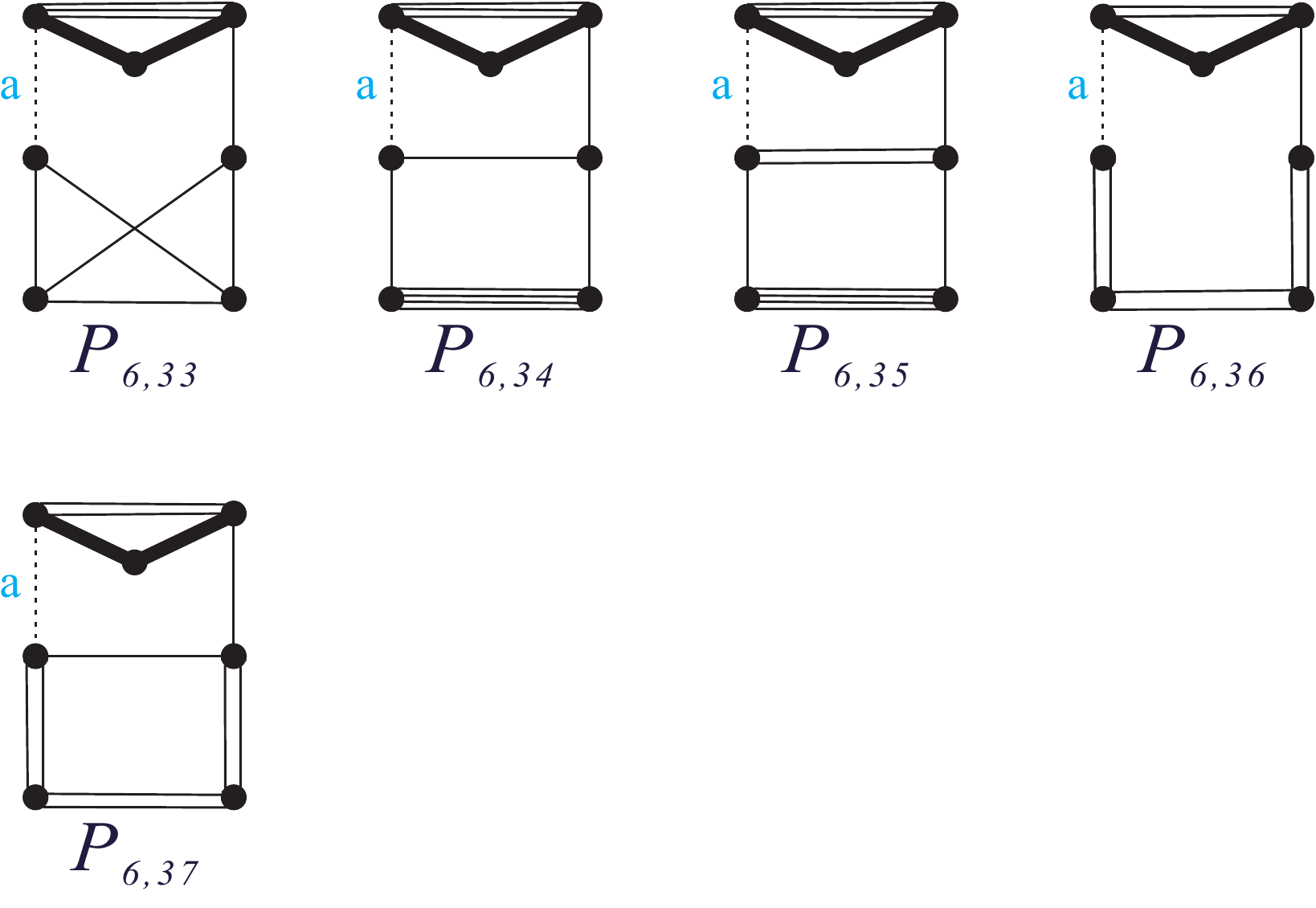}}
	\caption{The $37$ finite-volume hyperbolic Coxeter $4$-polytopes with $7$ facets over polytope $P_{6}$ (part 3).} \label{figure:p63}
\end{figure}

\newpage
4. Coxeter diagrams for $P_{6}$ (Figure \ref{figure:p61}--\ref{figure:p63}).

\newpage

\vspace{0.5cm}

\newpage

\vspace{0.5cm}

\newpage

\begin{figure}[H]
	\scalebox{0.47}[0.47]{\includegraphics {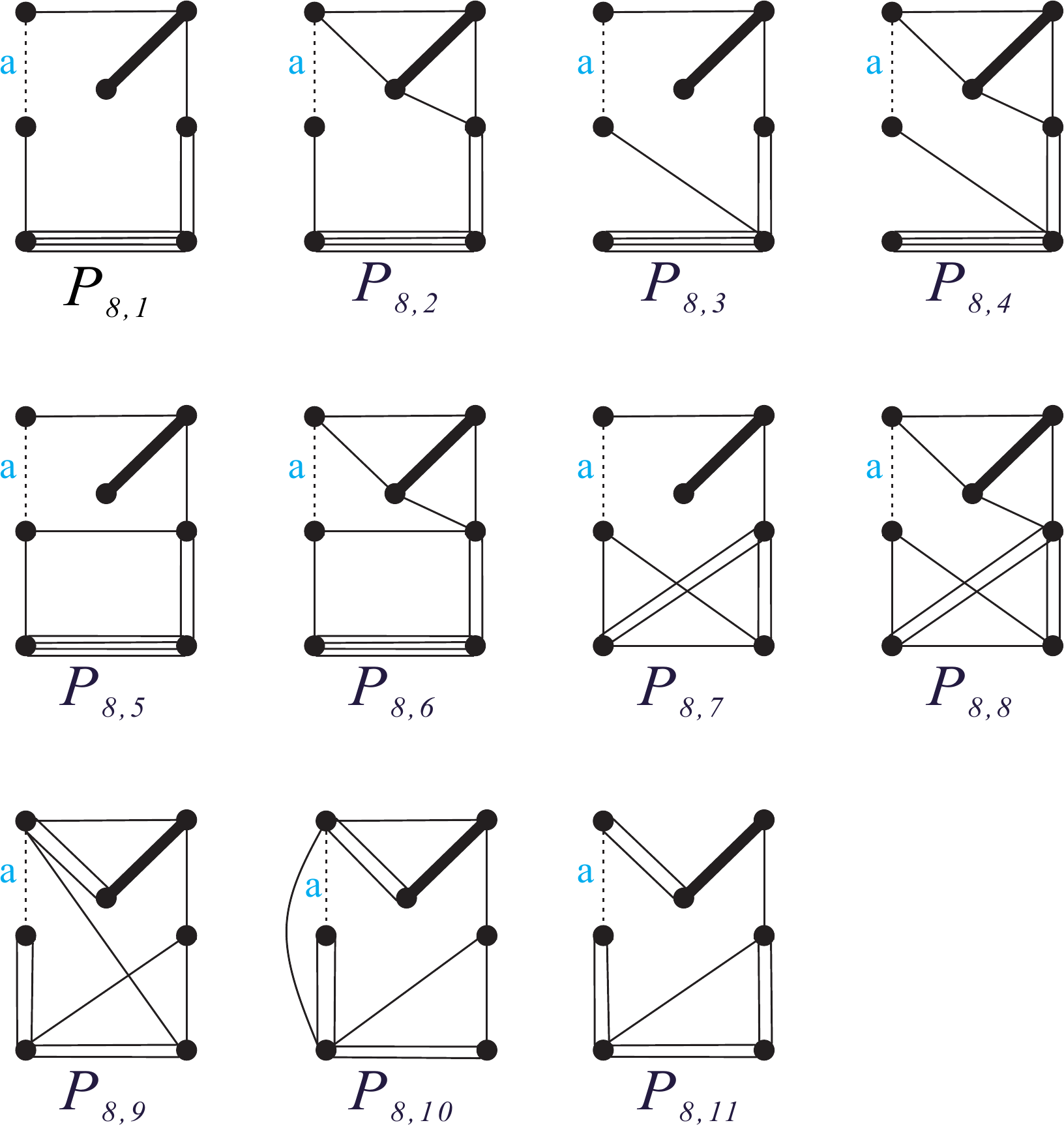}}
	\caption{There are $11$ finite-volume hyperbolic Coxeter $4$-polytopes with $7$ facets over polytope $P_{8}$.} \label{figure:p8}
\end{figure}

6. Coxeter diagrams for $P_{8}$ (Figure \ref{figure:p8}).

\newpage

\begin{figure}[H]
	\scalebox{0.47}[0.47]{\includegraphics {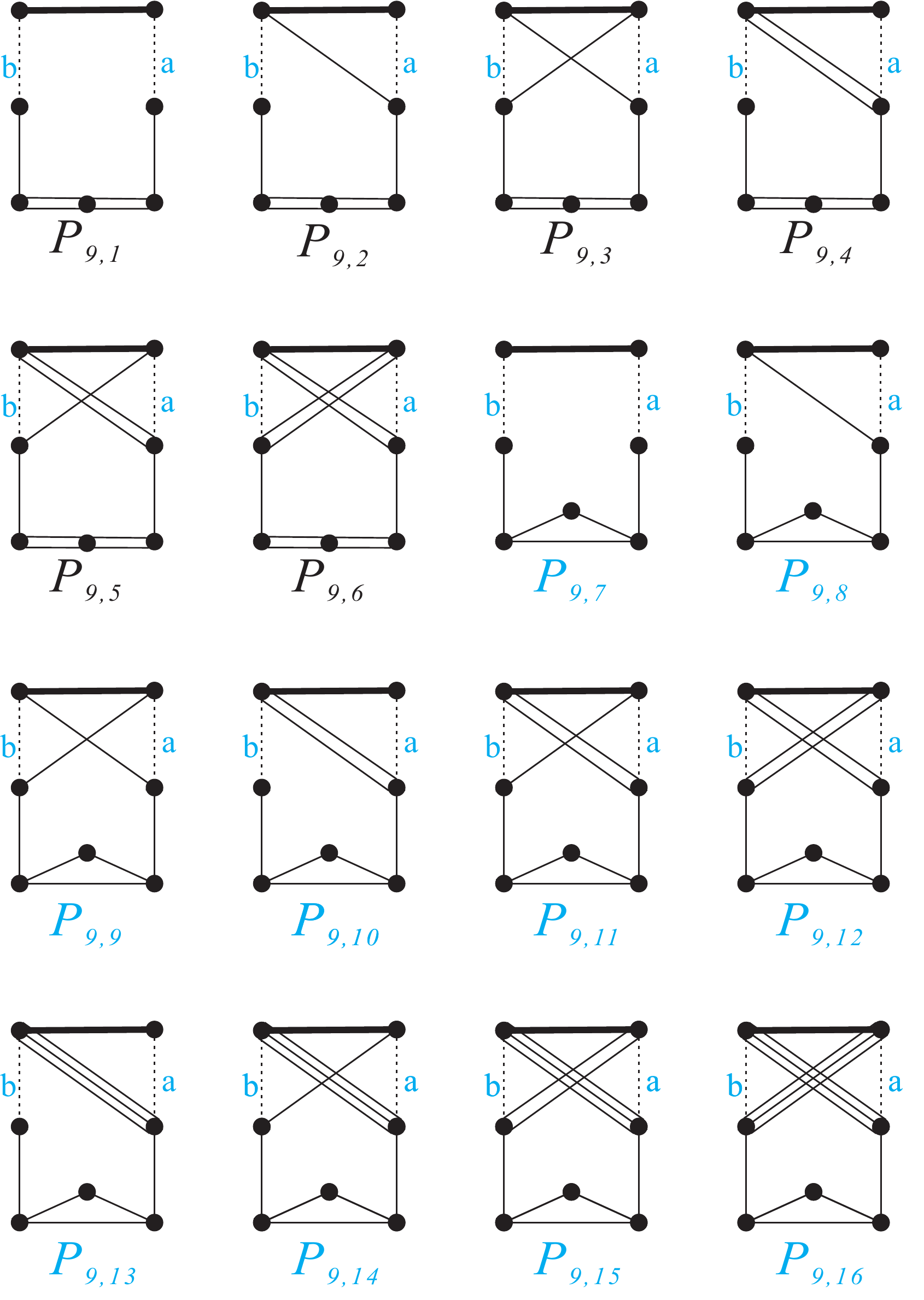}}
	\caption{The $134$ finite-volume hyperbolic Coxeter $4$-polytopes with $7$ facets over polytope $P_{9}$ (part 1).} \label{figure:p91}
\end{figure}

\begin{figure}[H]
	\scalebox{0.47}[0.47]{\includegraphics {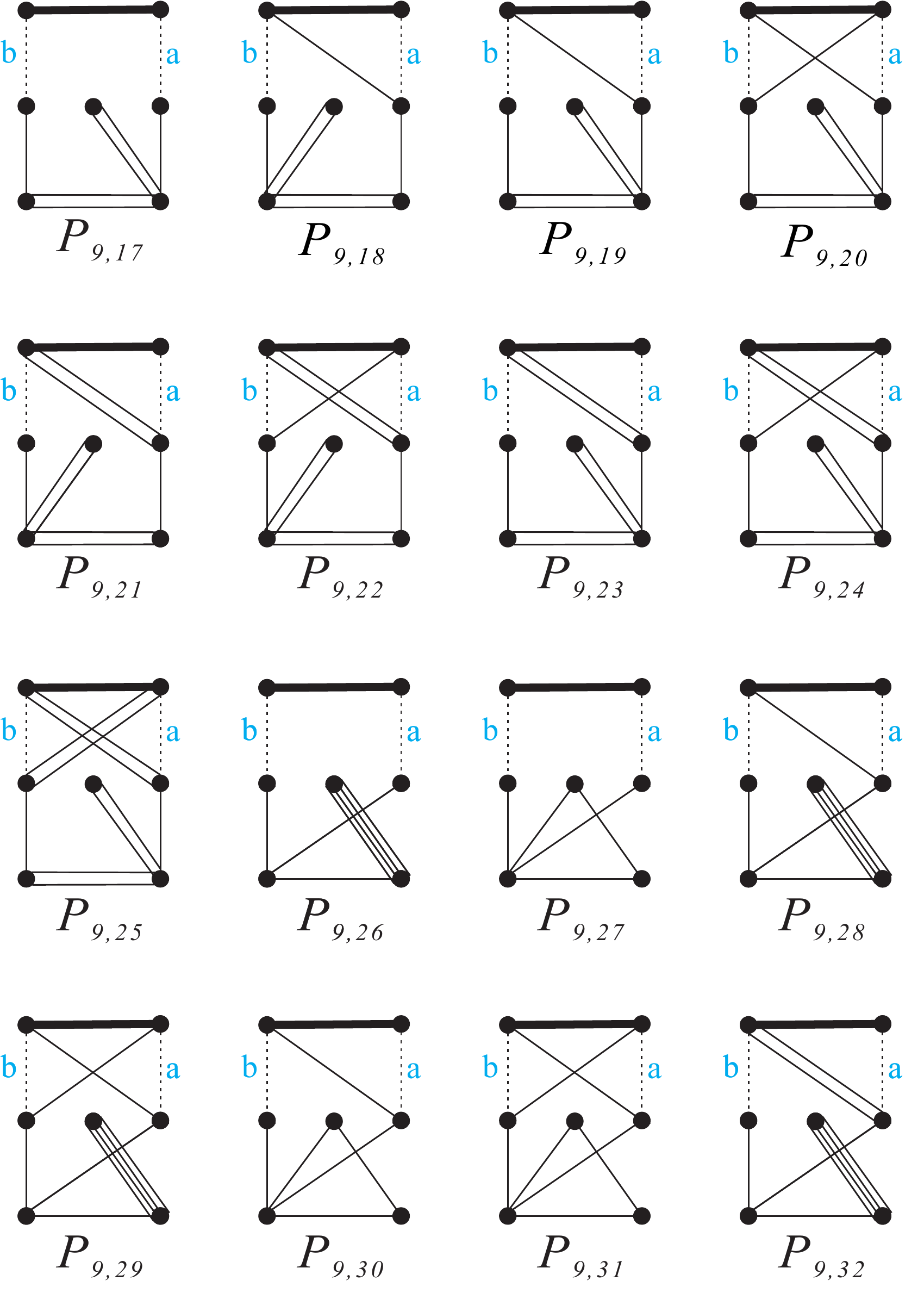}}
	\caption{The $134$ finite-volume hyperbolic Coxeter $4$-polytopes with $7$ facets over polytope $P_{9}$ (part 2).} \label{figure:p92}
\end{figure}

\begin{figure}[H]
	\scalebox{0.47}[0.47]{\includegraphics {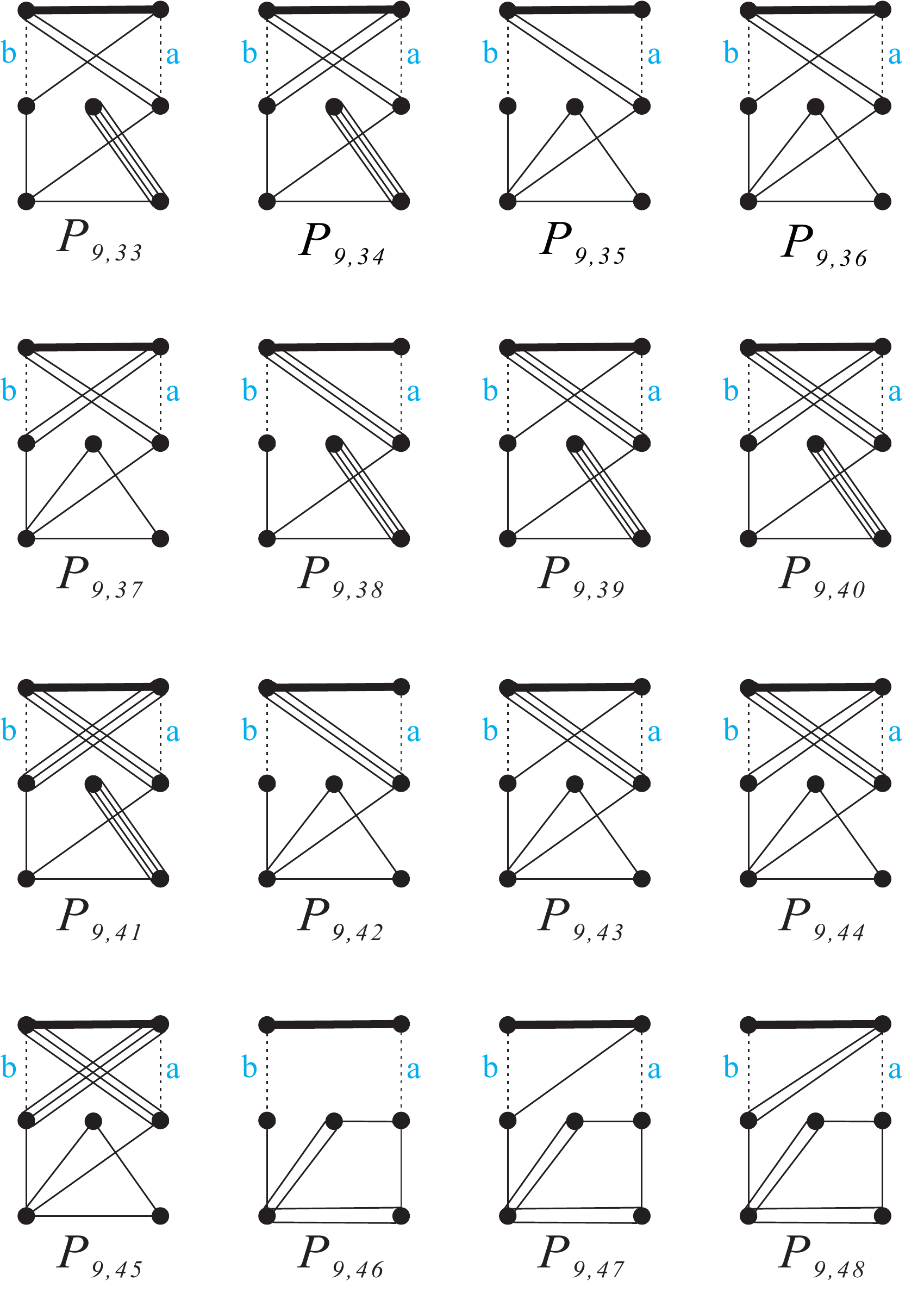}}
	\caption{The $134$ finite-volume hyperbolic Coxeter $4$-polytopes with $7$ facets over polytope $P_{9}$ (part 3).} \label{figure:p93}
\end{figure}

\begin{figure}[H]
	\scalebox{0.47}[0.47]{\includegraphics {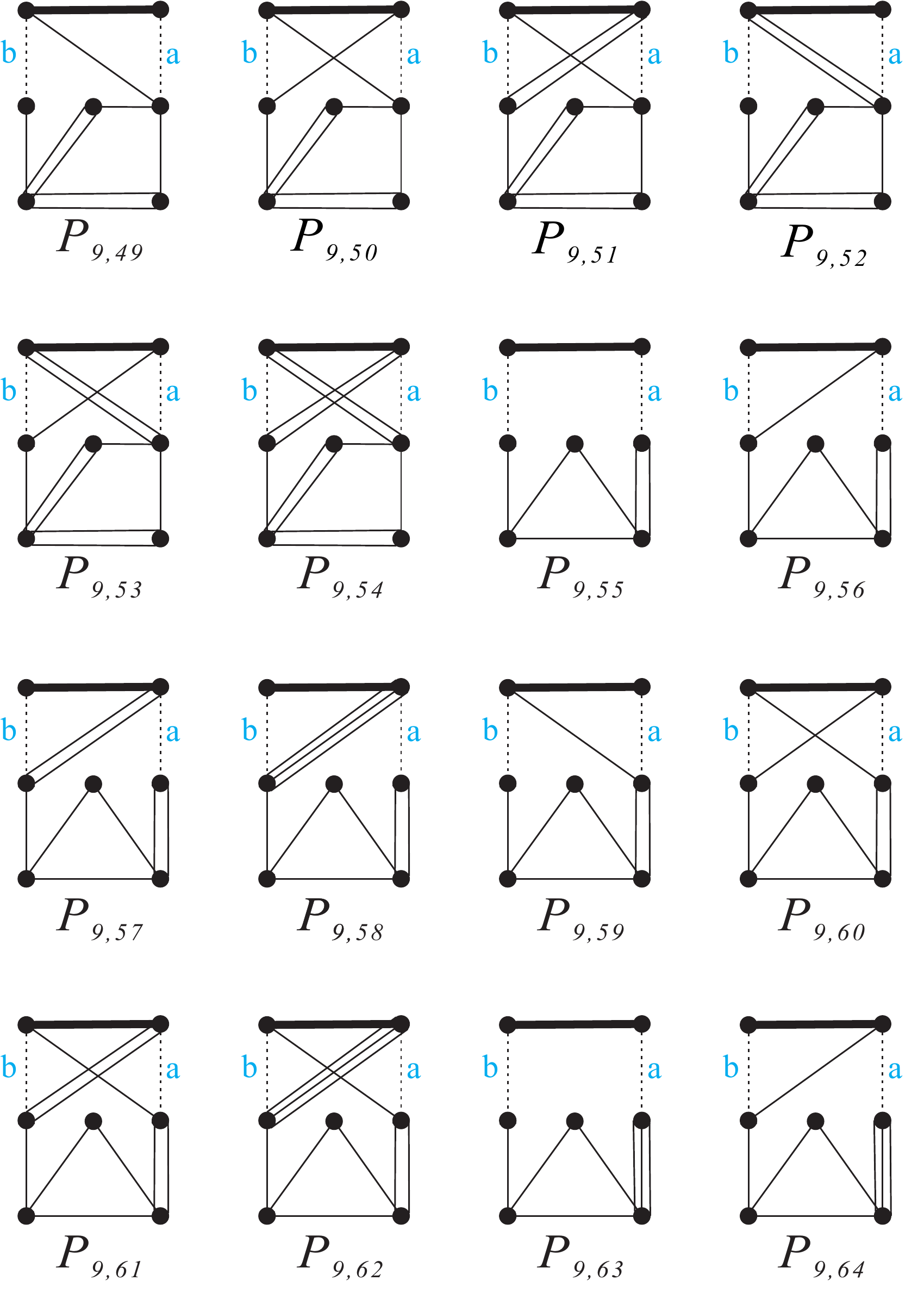}}
	\caption{The $134$ finite-volume hyperbolic Coxeter $4$-polytopes with $7$ facets over polytope $P_{9}$ (part 4).} \label{figure:p94}
\end{figure}

\begin{figure}[H]
	\scalebox{0.47}[0.47]{\includegraphics {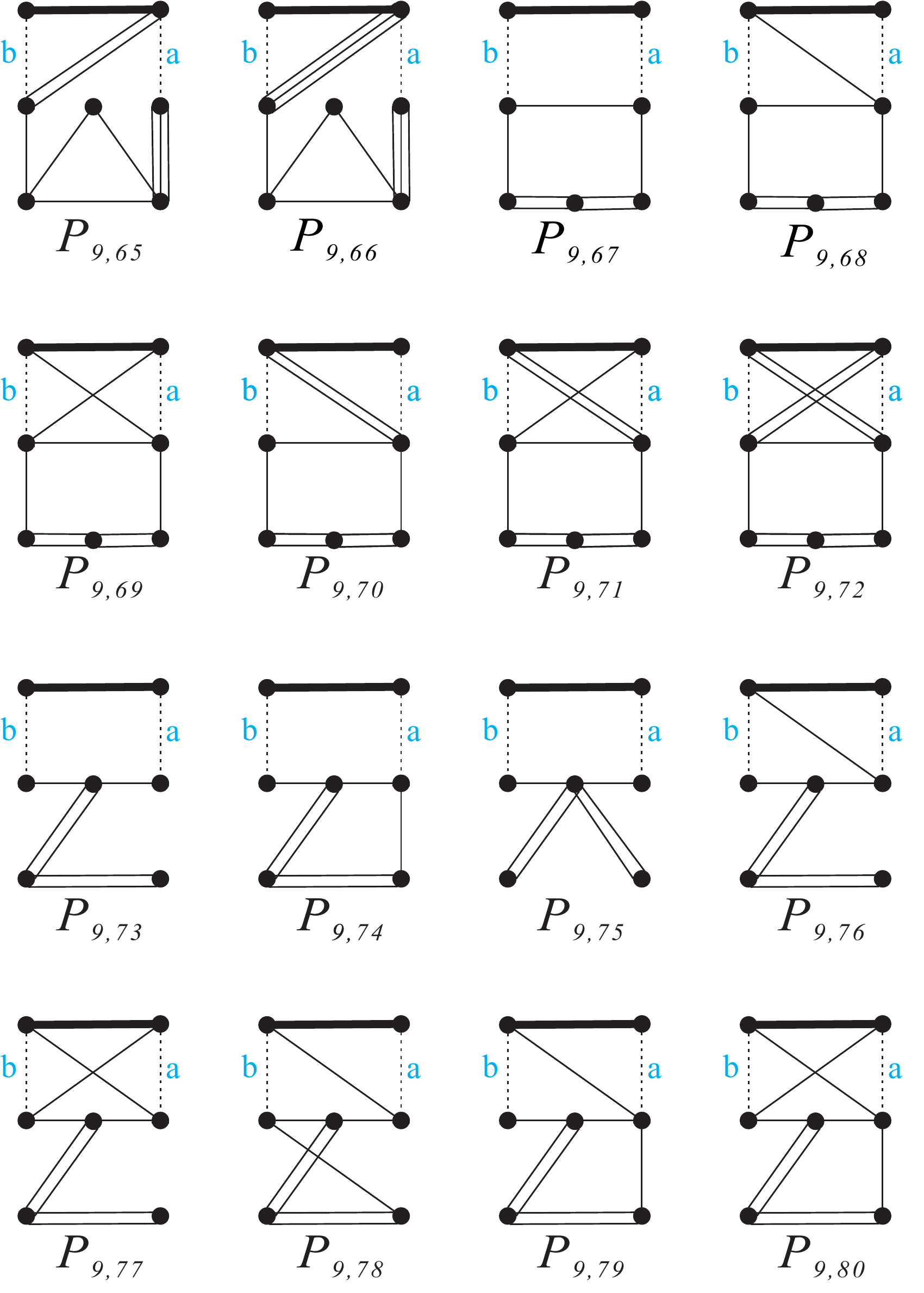}}
	\caption{The $134$ finite-volume hyperbolic Coxeter $4$-polytopes with $7$ facets over polytope $P_{9}$ (part 5).} \label{figure:p95}
\end{figure}

\begin{figure}[H]
	\scalebox{0.47}[0.47]{\includegraphics {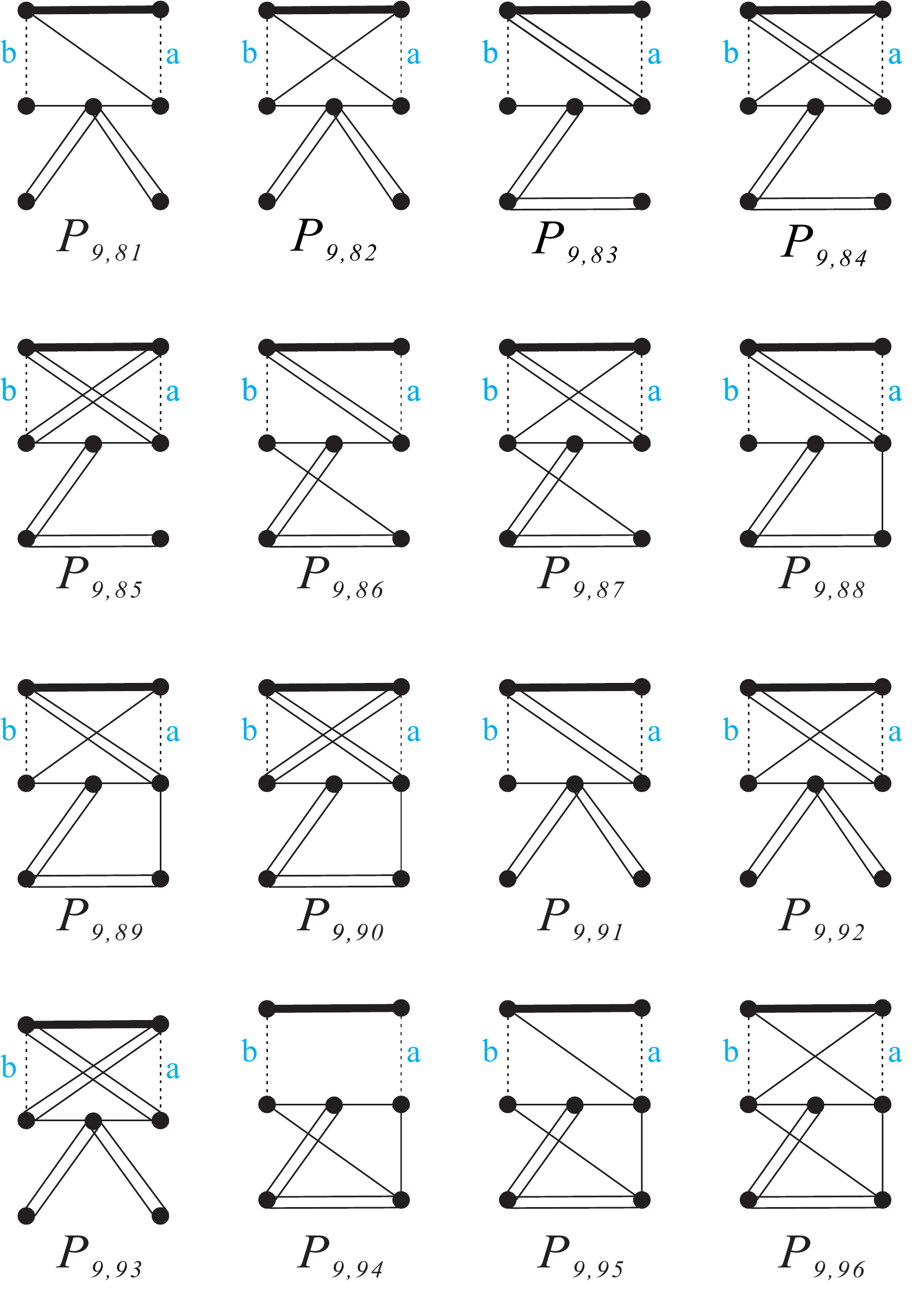}}
	\caption{The $134$ finite-volume hyperbolic Coxeter $4$-polytopes with $7$ facets over polytope $P_{9}$ (part 6).} \label{figure:p96}
\end{figure}

\begin{figure}[H]
	\scalebox{0.47}[0.47]{\includegraphics {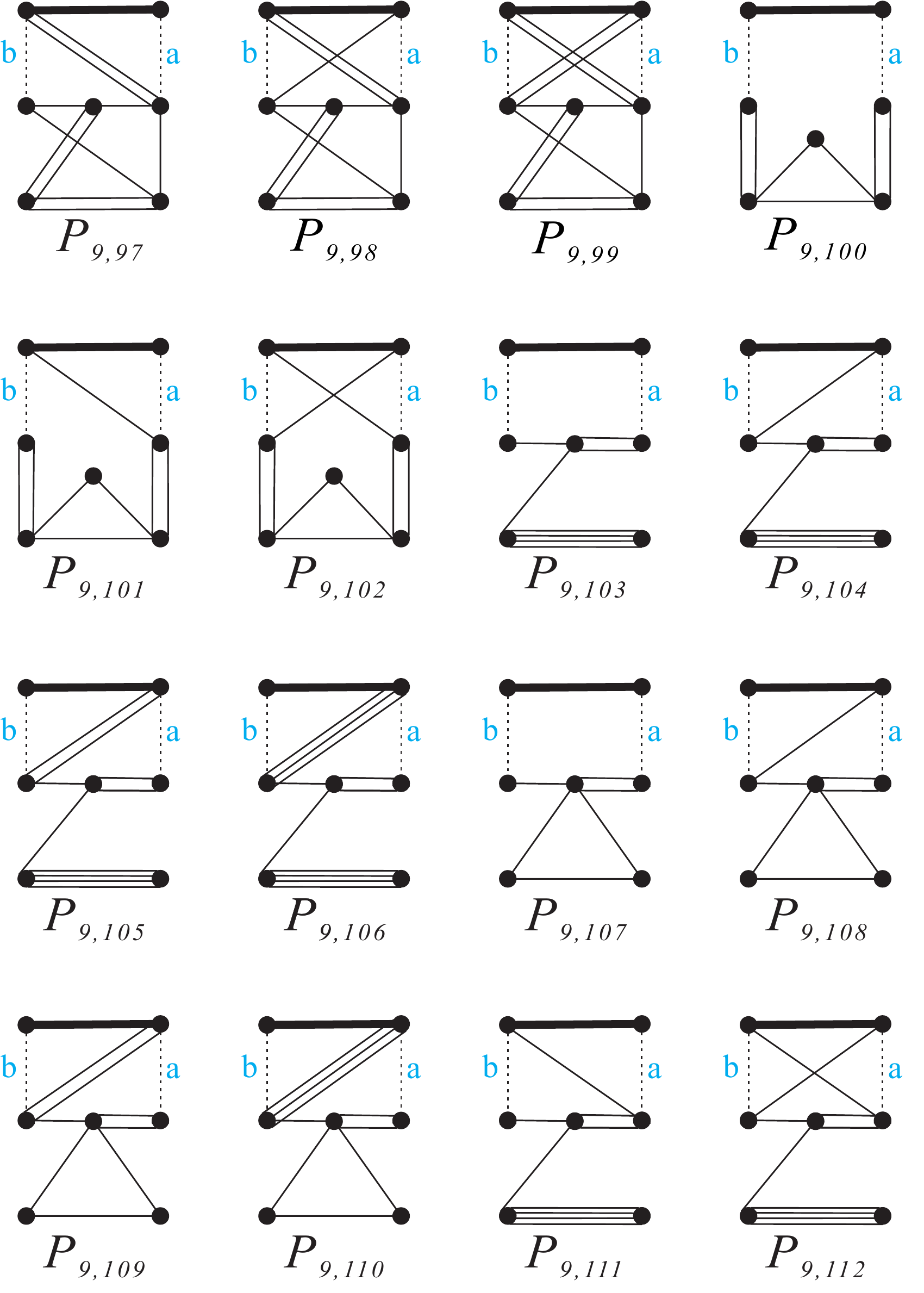}}
	\caption{The $134$ finite-volume hyperbolic Coxeter $4$-polytopes with $7$ facets over polytope $P_{9}$ (part 7).} \label{figure:p97}
\end{figure}

\begin{figure}[H]
	\scalebox{0.47}[0.47]{\includegraphics {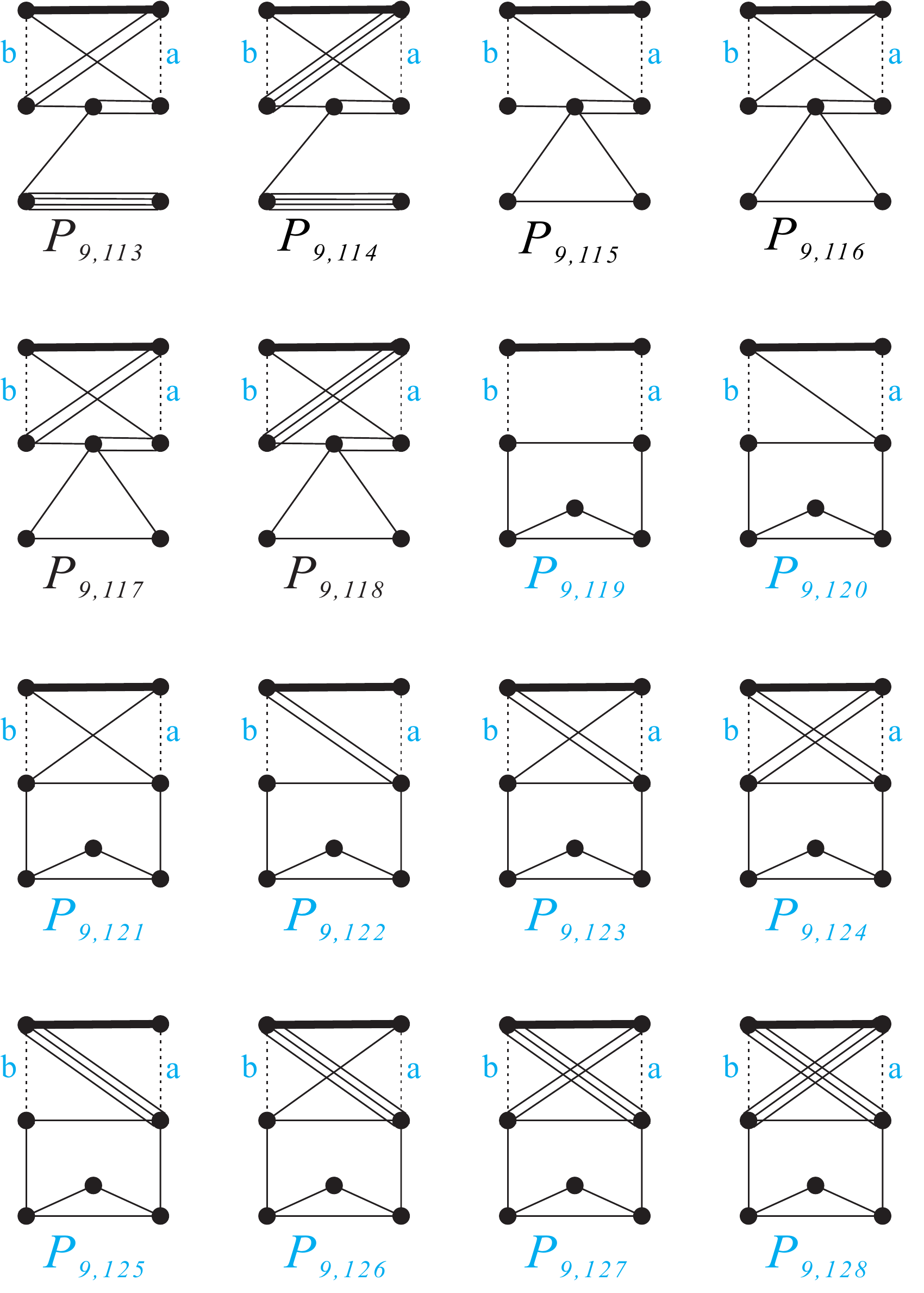}}
	\caption{The $134$ finite-volume hyperbolic Coxeter $4$-polytopes with $7$ facets over polytope $P_{9}$ (part 8).} \label{figure:p98}
\end{figure}

\begin{figure}[H]
	\scalebox{0.47}[0.47]{\includegraphics {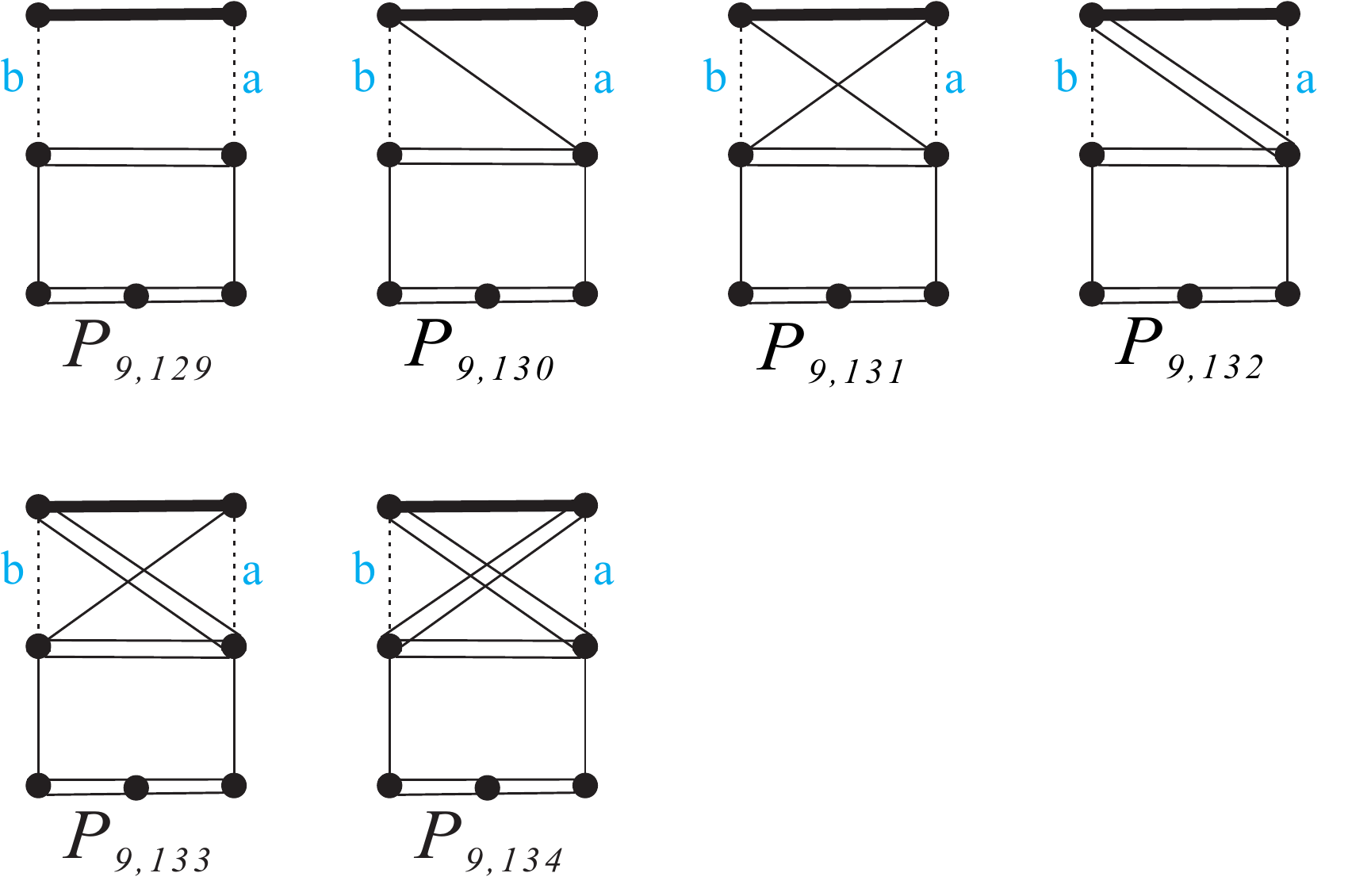}}
	\caption{The $134$ finite-volume hyperbolic Coxeter $4$-polytopes with $7$ facets over polytope $P_{9}$ (part 9).} \label{figure:p99}
\end{figure}

7. Coxeter diagrams for $P_{9}$ (Figure \ref{figure:p91}--\ref{figure:p99}).

\newpage

\newpage

\newpage

\newpage

\newpage

\newpage

\newpage

\newpage

\begin{figure}[H]
	\scalebox{0.47}[0.47]{\includegraphics {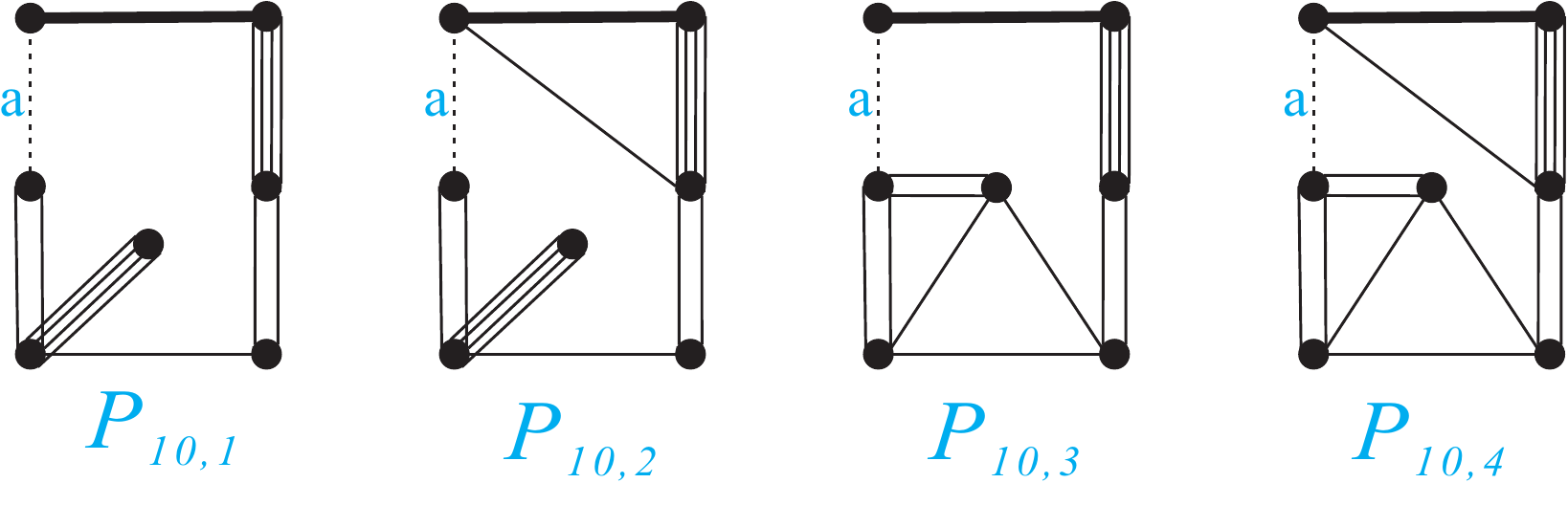}}
	\caption{The $4$ finite-volume hyperbolic Coxeter $4$-polytopes with $7$ facets over polytope $P_{10}$.} \label{figure:p10}
\end{figure}

8. Coxeter diagrams for $P_{10}$ (Figure \ref{figure:p10}).

\vspace{0.5cm}

\newpage

\begin{figure}[H]
	\scalebox{0.47}[0.47]{\includegraphics {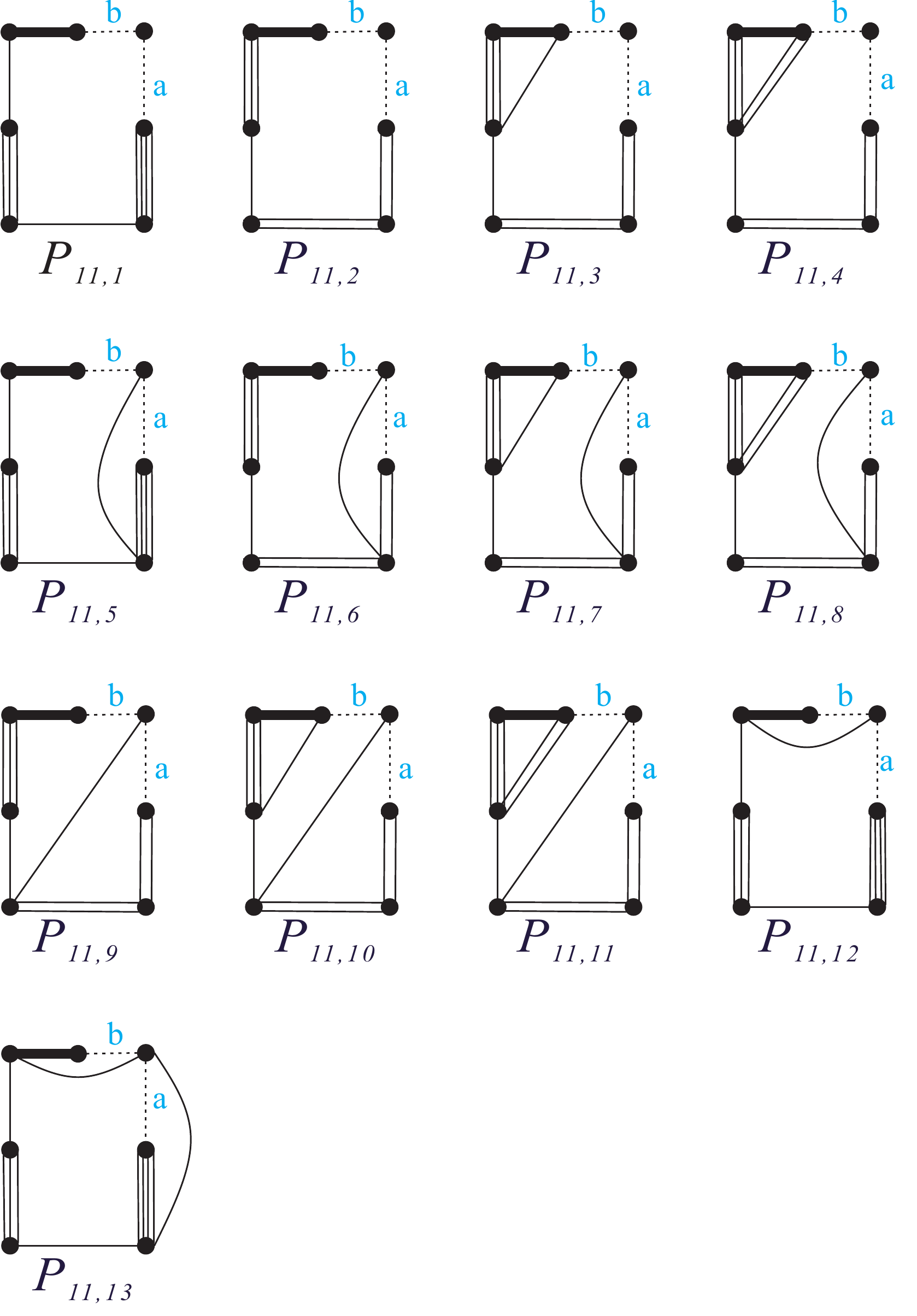}}
	\caption{The $13$ finite-volume hyperbolic Coxeter $4$-polytopes with $7$ facets over polytope $P_{11}$.} \label{figure:p11}
\end{figure}

9. Coxeter diagrams for $P_{11}$ (Figure \ref{figure:p11}).

\newpage

\begin{figure}[H]
	\scalebox{0.47}[0.47]{\includegraphics {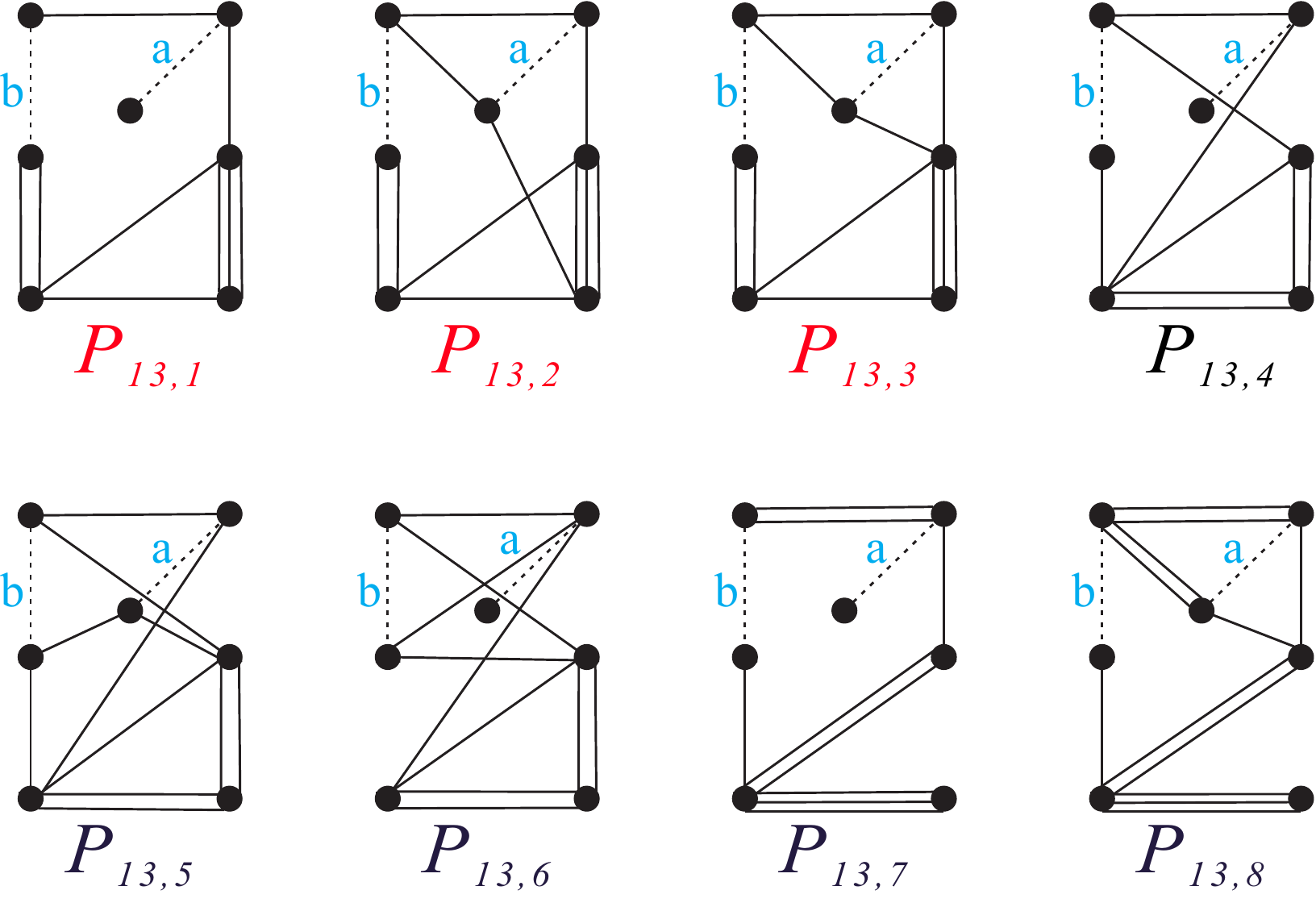}}
	\caption{The $8$ finite-volume hyperbolic Coxeter $4$-polytopes with $7$ facets over polytope $P_{13}$.} \label{figure:p13}
\end{figure}

10. Coxeter diagrams for $P_{13}$ (Figure \ref{figure:p13}).

\begin{figure}[H]
	\scalebox{0.47}[0.47]{\includegraphics {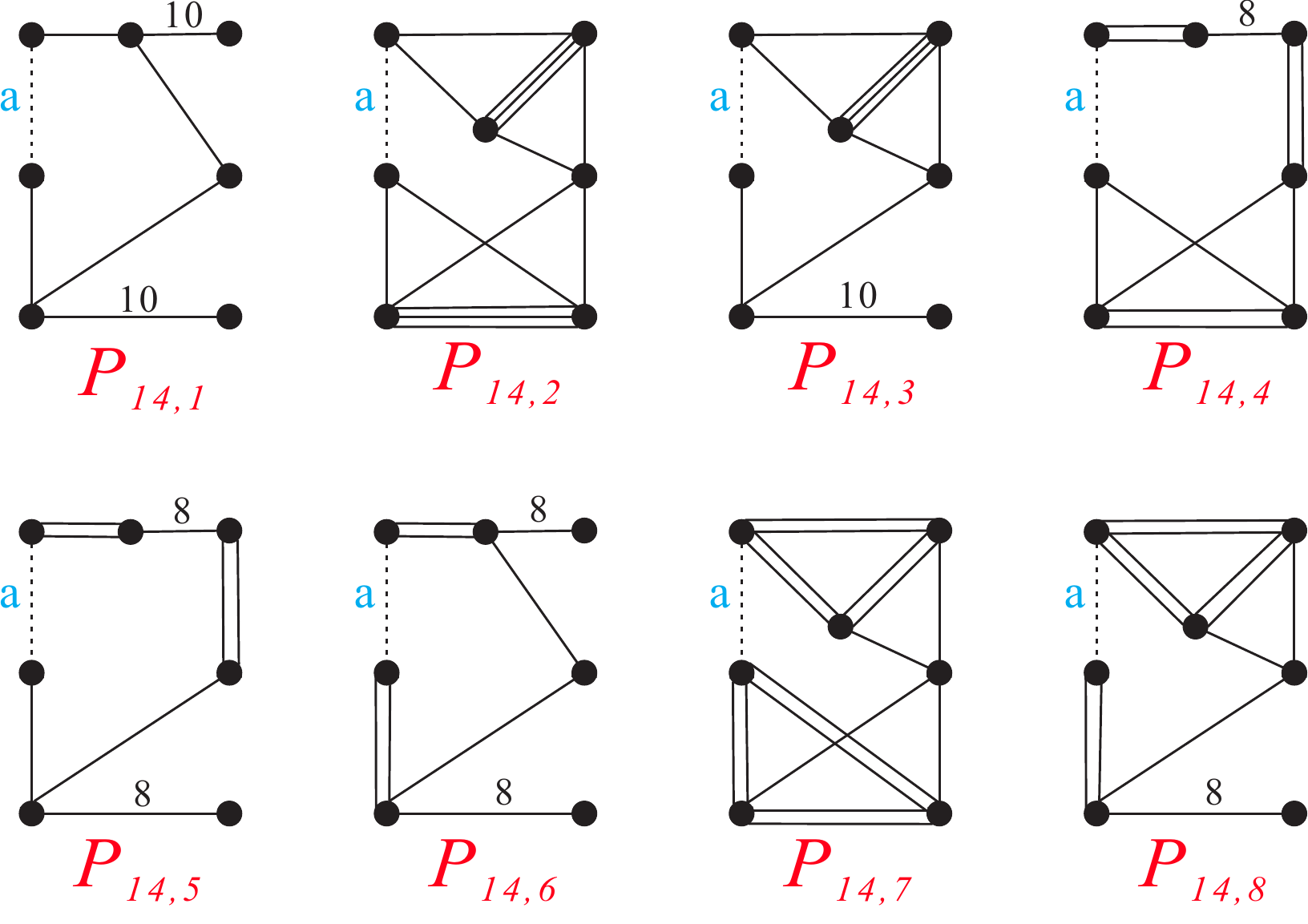}}
	\caption{The $8$ finite-volume hyperbolic Coxeter $4$-polytopes with $7$ facets over polytope $P_{14}$.} \label{figure:p14}
\end{figure}

11. Coxeter diagrams for $P_{14}$ (Figure \ref{figure:p14}).

\newpage

\begin{figure}[H]
	\scalebox{0.47}[0.47]{\includegraphics {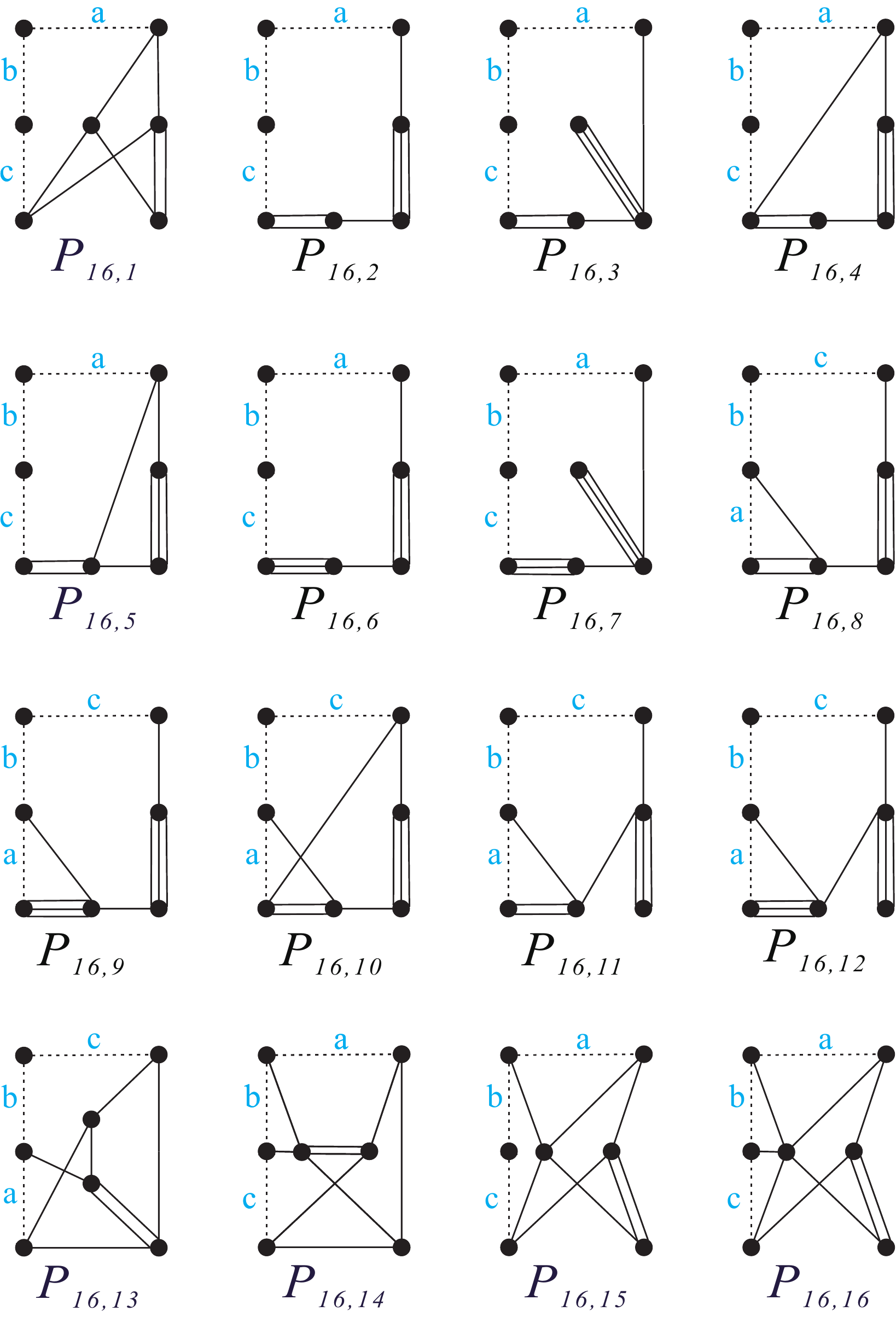}}
	\caption{The $81$ finite-volume hyperbolic Coxeter $4$-polytopes with $7$ facets over polytope $P_{16}$ (part 1).} \label{figure:p191}
\end{figure}

\begin{figure}[H]
	\scalebox{0.47}[0.47]{\includegraphics {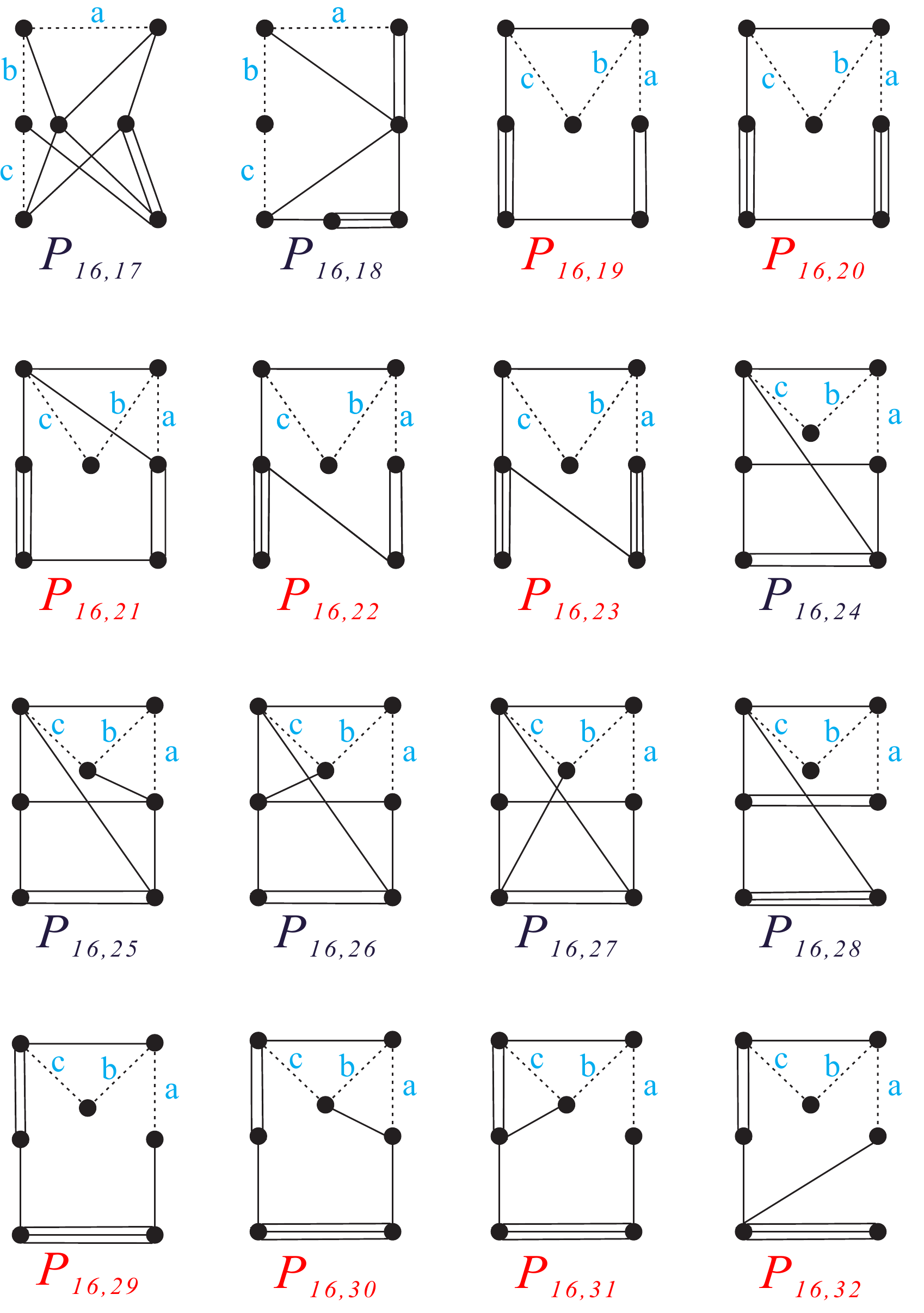}}
	\caption{The $81$ finite-volume hyperbolic Coxeter $4$-polytopes with $7$ facets over polytope $P_{16}$ (part 2).} \label{figure:p192}
\end{figure}

\begin{figure}[H]
	\scalebox{0.47}[0.47]{\includegraphics {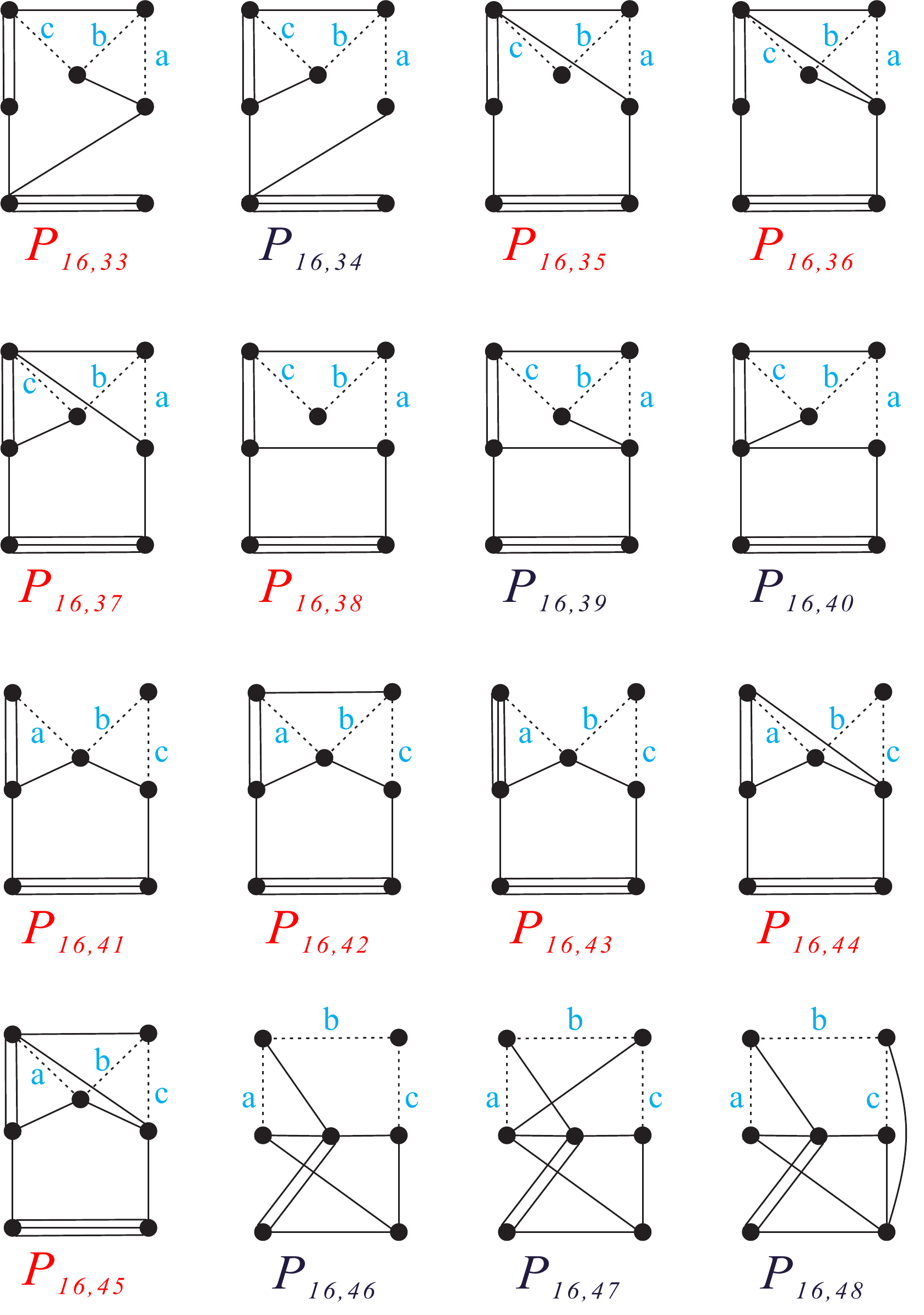}}
	\caption{The $81$ finite-volume hyperbolic Coxeter $4$-polytopes with $7$ facets over polytope $P_{16}$ (part 3).} \label{figure:p193}
\end{figure}

\begin{figure}[H]
	\scalebox{0.47}[0.47]{\includegraphics {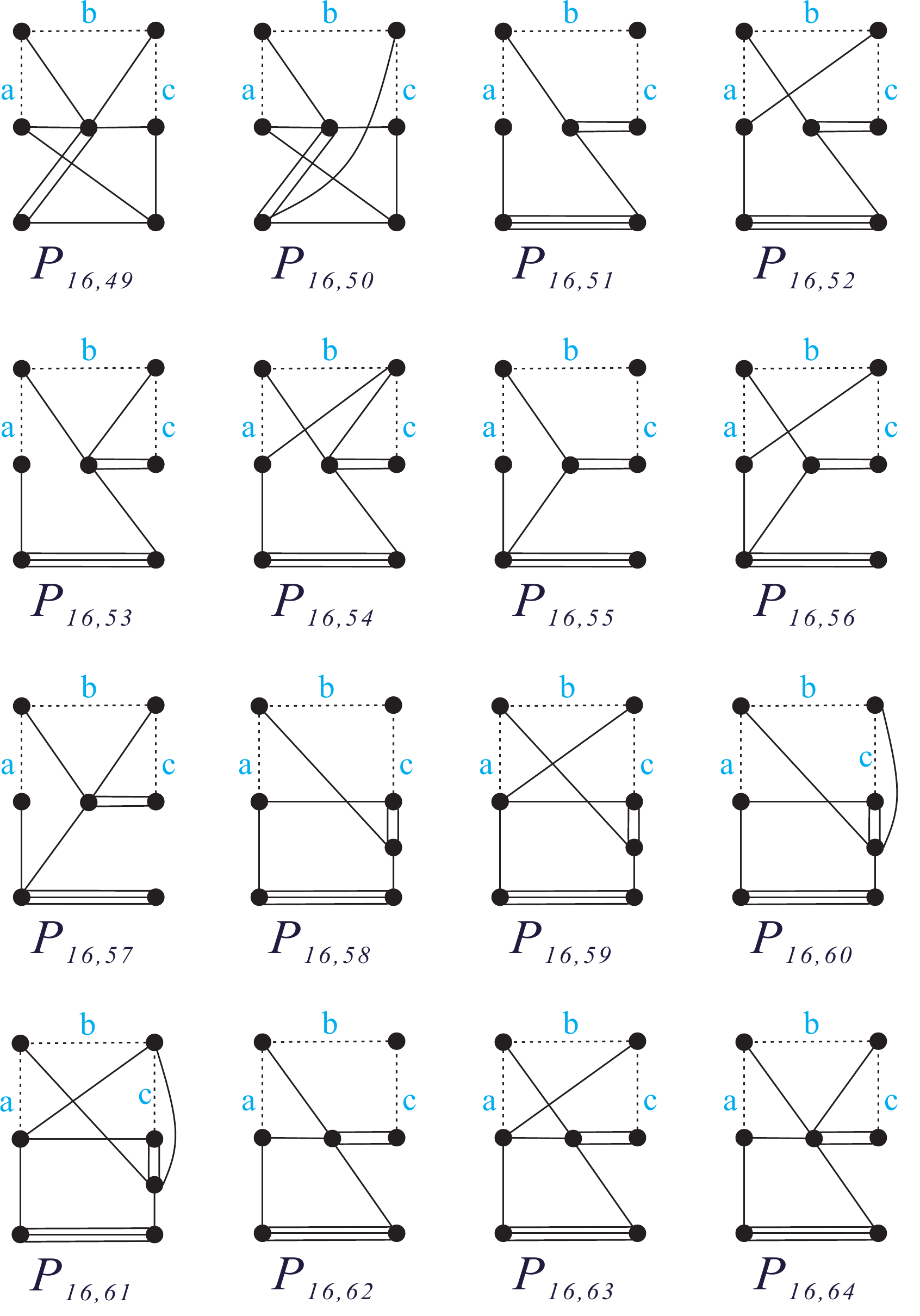}}
	\caption{The $81$ finite-volume hyperbolic Coxeter $4$-polytopes with $7$ facets over polytope $P_{16}$ (part 4).} \label{figure:p194}
\end{figure}

\begin{figure}[H]
	\scalebox{0.47}[0.47]{\includegraphics {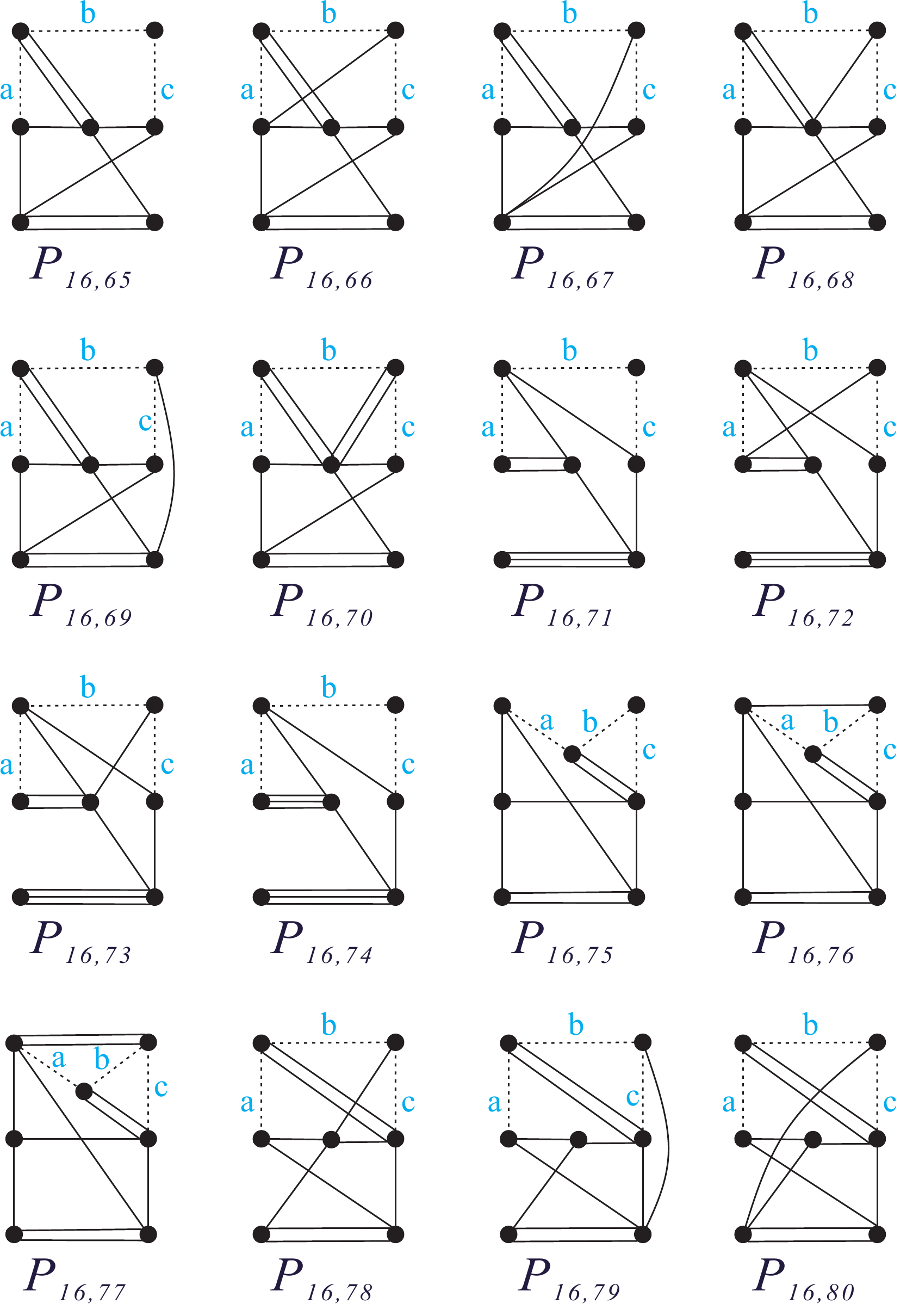}}
	\caption{The $81$ finite-volume hyperbolic Coxeter $4$-polytopes with $7$ facets over polytope $P_{16}$ (part 5).} \label{figure:p195}
\end{figure}

\begin{figure}[H]
	\scalebox{0.47}[0.47]{\includegraphics {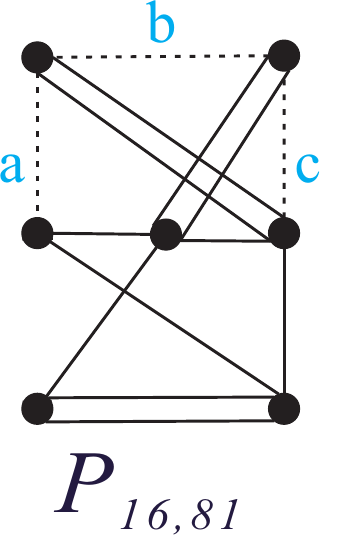}}
	\caption{The $81$ finite-volume hyperbolic Coxeter $4$-polytopes with $7$ facets over polytope $P_{16}$ (part 6).} \label{figure:p196}
\end{figure}

12. Coxeter diagrams for $P_{16}$ (Figure \ref{figure:p191}--\ref{figure:p196}).

\newpage

\newpage

\newpage

\newpage

\newpage

\end{document}